%% file: imp-imp-game.tex
\documentclass[11pt]{article}
\usepackage[english]{babel}
\usepackage{amsmath,amsthm,amssymb}
\usepackage{mathrsfs}
\usepackage{bbm}
\usepackage{amsfonts}
\usepackage[margin=0.8in]{geometry}
\usepackage{graphicx}
\usepackage{caption}
\usepackage{subcaption}
\usepackage{wasysym}
\usepackage{enumerate}
\usepackage{algorithm2e}
\usepackage{tabularx}
\usepackage{color}
\usepackage{wrapfig}
\usepackage{todonotes}

\newtheorem{thm}{Theorem}[section]
\newtheorem{ass}[thm]{Assumption}
\newtheorem{cor}[thm]{Corollary}
\newtheorem{lem}[thm]{Lemma}
\newtheorem{prop}[thm]{Proposition}
\newtheorem*{hyp*}{Hypothesis}
\theoremstyle{definition}
\newtheorem{defn}[thm]{Definition}

\theoremstyle{rem}
\newtheorem{rem}[thm]{Remark}

\numberwithin{equation}{section}
\newcommand{\R}{\mathbb R}
\newcommand{\eps}{\varepsilon}

\newcommand{\bbA}{\mathbb A}

\newcommand{\bbF}{\mathbb F}
\newcommand{\bbT}{\mathbb T}

\newcommand{\bbG}{\mathbb G}

\newcommand{\bbN}{\mathbb N}

\newcommand{\mcA}{\mathcal{A}}

\newcommand{\mcB}{\mathcal{B}}
\newcommand{\mcC}{\mathcal C}
\newcommand{\mcD}{\mathcal D}

\newcommand{\mcF}{\mathcal F}

\newcommand{\mcL}{\mathcal L}
\newcommand{\mcT}{\mathcal T}
\newcommand{\mcP}{\mathcal P}
\newcommand{\mcM}{\mathcal M}
\newcommand{\mcO}{\mathcal O}
\newcommand{\mcH}{\mathcal H}
\newcommand{\mcK}{\mathcal K}
\newcommand{\mcU}{\mathcal U}

\newcommand{\mcR}{\mathcal R}
\newcommand{\mcS}{\mathcal S}

\newcommand{\mcV}{\mathcal V}
\newcommand{\mcX}{\mathcal X}
\newcommand{\mcY}{\mathcal Y}
\newcommand{\mcZ}{\mathcal Z}

\newcommand{\E}{\mathbb{E}}
\newcommand{\Prob}{\mathbb{P}}

\newcommand{\esssup}{\mathop{\rm{ess}\sup}}
\newcommand{\essinf}{\mathop{\rm{ess}\inf}}
\newcommand{\argmax}{\mathop{\arg\max}}

\newcommand{\ett}{\mathbbm{1}}
\newcommand{\bigSET}{\mcD}

\newcommand{\cadlag}{c\`adl\`ag~}

\newcommand{\Pred}{\mathcal{P}}

\newcommand{\trace}{{\rm Tr}}

\newcommand{\jmp}{{\xi}}

\newcommand{\bigS}{\mathfrak S}

\newcommand{\infOP}{\mathcal M}
\newcommand{\supOP}{\mathcal W}

\newcommand{\ie}{\textit{i.e.\ }}
\newcommand{\eg}{\textit{e.g.\ }}

\begin{document}

\title{Zero-sum Stochastic Differential Games of Impulse Control with Random Intervention Costs}

\author{Magnus Perninge\footnote{M.\ Perninge is with the Department of Mathematics at M\"alardalen University,
Sweden. e-mail: magnus.perninge@lnu.se.}} %
\maketitle
% ----------------------------------------------------------------
\begin{abstract}
We consider a finite-horizon, zero-sum game in which both players control a stochastic differential equation by invoking impulses. We derive a control randomization formulation of the game and use the existence of a value for the randomized game to show that the upper and lower value functions of the original game coincide.

The main contribution of the present work is that we can allow intervention costs that are functions of the state as well as time, and that we do not need to impose any monotonicity assumptions on the involved coefficients.
\end{abstract}

% ----------------------------------------------------------------
\section{Introduction}
For $(t,x)\in [0,T]\times\R^d$ and impulse controls $u:=(\tau_j,\beta_j)_{j=1}^N$ and $\alpha:=(\eta_j,\theta_j)_{j=1}^M$ we let $X^{t,x;u,\alpha}$ solve the SDE with impulses driven by a $d$-dimensional Brownian motion, $W$,
\begin{align}\nonumber
X^{t,x;u,\alpha}_s&=x+\int_t^s a(r,X^{t,x;u,\alpha}_r)dr+\int_t^s\sigma(r,X^{t,x;u,\alpha}_r)dW_r+\sum_{j=1}^N\ett_{[\tau_j\leq s]}\jmp(\tau_j,X^{t,x;[u]_{j-1},\alpha}_{\tau_j},\beta_j)
\\
&\quad +\sum_{j=1}^M\ett_{[\eta_j\leq s]}\gamma(\eta_j,X^{t,x;u,[\alpha]_{j-1}}_{\eta_j},\theta_j),\label{ekv:fwd-sde}
\end{align}
where $N:=\max\{j:\tau_j< T\}$, $M:=\max\{j:\eta_j< T\}$ and $[\phi]_k$ is the truncation of $\phi=u,\alpha$ to the first $k$ interventions, so that \eg $[u]_k:=(\tau_j,\beta_j)_{j=1}^{k\wedge N}$. We then consider the zero-sum stochastic differential game of impulse control with upper (resp.~lower) value
\begin{align*}
  \bar v(t,x):=\esssup_{u^S\in\mcU^S_t}\essinf_{\alpha\in\mcA_t}J(t,x;u^S(\alpha),\alpha) \qquad (\text{resp.~}\underline v(t,x):=\essinf_{\alpha^S\in\mcA^S_t}\esssup_{u\in\mcU_t}J(t,x;u,\alpha^S(u))),
\end{align*}
where
\begin{align}\nonumber
  J(t,x;u,\alpha)&:=\E\Big[\psi(X^{t,x;u,\alpha}_T)+\int_t^{T}f(s,X^{t,x;u,\alpha}_s)ds - \sum_{j=1}^N\ell(\tau_j,X^{t,x;[u]_{j-1},\alpha},\beta_j)
  \\
  &\qquad + \sum_{j=1}^M\chi(\eta_j,X^{t,x;u,[\alpha]_{j-1}},\theta_j)\Big|\mcF_t\Big],\label{ekv:J-def}
\end{align}
and $\mcU^S_t$ (resp.~$\mcA^S_t$) is the set of non-anticipative maps from the set of impulse controls available to Player 1 at time $t\in[0,T]$, denoted $\mcU_t$, to the corresponding set of impulse controls for Player 2, denoted $\mcA_t$ (resp.~$\mcA_t\to\mcU_t$).

To achieve a fair game that does not reward preemption, we assume that for any $(t,x,b,e)\in [0,T] \times \R^d \times U \times A$, where $U$ (resp.~$A$) is a compact subset of $\R^d$ in which $\beta_j$ (resp.~$\alpha_j$) is valued, the jump functions satisfy the commutativity property
\begin{align}\label{ekv:jump-commute}
  \jmp(t,x+\gamma(t,x,e),b)=\gamma(t,x+\jmp(t,x,b),e)
\end{align}
whereas the intervention costs satisfy
\begin{align}\label{ekv:interv-cost-commute}
  \chi(t,x,e)=\chi(t,x+\jmp(t,x,b),e)\quad \text{and}\quad \ell(t,x,b)=\ell(t,x+\gamma(t,x,e),b).
\end{align}
These conditions are designed to rule out the importance of the order in which the players act and are crucial in proving that our game has a value, \ie that $\underline v \equiv \bar v$.

An example of a setting that matches the above assumptions is when Player $i$ ($i=1,2$) controls the state $X^{P_i}$ of her individual subsystem and there is a general system with state $X^G$ whose dynamics may depend on the state itself, but also on the vector $(X^{P_1},X^{P_2})$, whereas the jump function and the impulse cost of Player $i$ is a function of time, the chosen intervention and $(X^{G},X^{P_i})$.

We show that the game has a value by considering a dual formulation where the impulse control of the second player is randomized. This leads us to consider families of penalized reflected BSDEs driven by the Brownian motion, $W$, and an independent Poisson random measure $\mu$ on $[0,T]\times A$, with a compensator $\lambda$ that has full topological support. In particular, for each $k,n\in\bbN$ we assume that the family $(Y^{t,x,k,n},Z^{t,x,k,n},V^{t,x,k,n},K^{+,t,x,k,n}:(t,x)\in [0,T]\times\R^d)$ satisfies
\begin{align}\label{ekv:rbsde-pen}
  \begin{cases}
    Y^{t,x,k,n}_s=\psi(X^{t,x}_t)+\int_s^T f(r,X^{t,x}_r)dr-\int_s^T Z^{t,x,k,n}_r dW_r-\int_{s}^T\!\!\!\int_U V^{t,x,k,n}_{r-}(e)\mu(dr,de)
    \\
    \quad-n\int_s^T\!\!\!\int_A(V^{t,x,k,n}_{r}(e)+\chi(r,X^{t,x}_r,e))^-\lambda(de)ds+K^{+,t,x,k,n}_T-K^{t,x,k,n}_s,
    \\
    Y^{u,n}_s\geq \supOP v_{k-1,n}(s,X^{t,x}_s),\, \forall s\in [t,T],\quad \int_t^T\!\! \big(Y^{t,x,k,n}_s-\supOP v_{k-1,n}(s,X^{t,x}_s)\big)dK^{+,t,x,k,n}_s.
  \end{cases}
\end{align}
In \eqref{ekv:rbsde-pen}, the map $v_{k-1,n}:[0,T]\times\R^d\to\R$ is the (deterministic) function for which $v_{k-1,n}(t,x)=Y^{t,x,k-1,n}_t$, $\Prob$-a.s., the operator $\supOP$ is defined as $\supOP v(t,x):=\sup_{b\in U}\{v(t,x+\jmp(t,x,b))-\ell(t,x,b)\}$ and $X^{t,x}$ solves the standard SDE with jumps
\begin{align}
X^{t,x}_s&=x+\int_t^s a(r,X^{t,x}_r)dr+\int_t^s\sigma(r,X^{t,x}_r)dW_r + \int_t^s\!\!\!\int_A\gamma(r,X^{t,x}_{r-},e)\mu(dr,de),\quad\forall s\in [t,T].\label{ekv:fwd-sde-Markov}
\end{align}

A main result from the accompanying paper~\cite{dbl-obst-qvi} is that
\begin{align*}
\lim_{k\to\infty}\lim_{n\to\infty}  Y^{t,x,k,n}_t =\lim_{n\to\infty} \lim_{k\to\infty}  Y^{t,x,k,n}_t = v(t,x),\quad\forall (t,x)\in [0,T]\times\R^d,
\end{align*}
where the map $v:[0,T]\times\R^d\to\R$ is the unique viscosity solution to the double obstacle quasi-variational inequality
\begin{align}\label{ekv:var-ineq}
\begin{cases}
 \min\{v(t,x)-\supOP v(t,x), \max\{v(t,x)-\infOP v(t,x),-\frac{\partial}{\partial t}v(t,x)-\mcL v(t,x)\\
  \quad-f(t,x)\}\}=0,\quad\forall (t,x)\in[0,T)\times \R^d \\
  v(T,x)=\psi(x).
\end{cases}
\end{align}
The operator $\infOP$ in \eqref{ekv:var-ineq} is defined as $\infOP v(t,x):=\inf_{e\in A}\{v(t,x+\gamma(t,x,e))+\chi(t,x,e)\}$ and
\begin{align}
  \mcL:=\sum_{j=1}^d a_j(t,x)\frac{\partial}{\partial x_j}+\frac{1}{2}\sum_{i,j=1}^d (\sigma\sigma^\top(t,x))_{i,j}\frac{\partial^2}{\partial x_i\partial x_j}
\end{align}
is the infinitesimal generator related to the SDE
\begin{align*}
\check X^{t,x}_s=x + \int_t^s a(r,\check X^{t,x}_r)dr + \int_t^s \sigma(r,\check X^{t,x}_r)dW_r,\quad\forall s\in[t,T].
\end{align*}
We show that under the additional assumptions in \eqref{ekv:jump-commute}--\eqref{ekv:interv-cost-commute}, $v=\underline v=\bar v$. Exploiting the symmetry in our setup we deduce that $v$ is the unique solution to the max-min version of \eqref{ekv:var-ineq} as well. As a byproduct of our work, we thereby extend the results of \cite{dbl-obst-qvi} by proving that, when the driver $f$ only depends on $t$ and $x$ and \eqref{ekv:jump-commute}--\eqref{ekv:interv-cost-commute} holds, the unique viscosity solutions to the min-max and the max-min versions of the corresponding double obstacle QVIs agree.

\bigskip

\textbf{Related literature}  In the case when the impulse costs $\ell$ and $\chi$ are non-increasing in time and independent of the state variable, $x$, a zero-sum game of impulse control was resolved in \cite{Cosso13} by proving that the corresponding value function is indeed the unique viscosity solution to a double obstacle QVI. The results presented in \cite{Cosso13} extended previous results from \cite{TangHou07}, where one player was restricted to play a switching control. Since then, zero-sum games where both players implement switching controls have been considered in a more general framework in \cite{BollanSWG1} and \cite{Djehiche17}.

In \cite{dbl-obst-qvi}, existence of a unique viscosity solution, $v$, to \eqref{ekv:var-ineq} was proved by deriving a Feynman-Kac representation, relating $v$ to the solution of a BSDE with constrained jumps. The ensemble of approaches to solve stochastic optimal control problems that utilize BSDEs with constrained jumps is commonly referred to as \emph{control randomization}. The control randomization methodology was initiated by the seminal work \cite{Kharroubi2010}, which related the solution of a BSDE with constrained jumps to the unique viscosity solution of a (single-barrier) quasi-variational inequality (see also \cite{Bouchard09}). A significant breakthrough was achieved in \cite{Fuhrman15}, directly linking the value function of the randomized control problem to that of the original control problem. This eliminated the need for a Feynman-Kac representation, thus extending the framework to encompass stochastic systems with path-dependencies. Building on this foundation, subsequent advancements adapted the methodology to partial information settings \cite{Bandini18} and optimal switching problems \cite{Fuhrman2020}. More recently, \cite{imp-stop-game} further expanded the scope of control randomization by establishing a link between optimal stopping for BSDEs with constrained jumps and zero-sum games.\\

\textbf{Outline} In the following section, we introduce the set of assumptions used throughout the paper, define the basic notations, and state the main result. In Section~\ref{sec:game-trunc}, we demonstrate that the value functions can be chosen independently of $\omega$ and present truncations of the value functions. Section~\ref{sec:dual-game} introduces the dual game and establishes a connection between the primal and dual games, encapsulated in the central Proposition~\ref{prop:rel-rand-to-orig}, which translates the main result stated in Section~\ref{sec:prel}. Finally, Section~\ref{sec:proof} provides a detailed proof of Proposition~\ref{prop:rel-rand-to-orig}, thereby establishing the main result.%In the next section, we present the set of assumptions that hold throughout and state the main result, in addition to setting the basic notations. Then, in Section~\ref{sec:game-trunc} we explain that the value functions can be chosen independent of $\omega$ and introduce truncations of the value functions. In Section~\ref{sec:dual-game}, the dual game is introduced and the main result, formulated in Section~\ref{sec:prel}, is translated to a relation between the primal and the dual game in the central Proposition~\ref{prop:rel-rand-to-orig}. Finally, Section~\ref{sec:proof} provides a proof for Proposition~\ref{prop:rel-rand-to-orig} thus proving the main result.

\section{Problem formulation and main result\label{sec:prel}}

\subsection{Assumptions}
%\begin{ass}\label{ass:lambda}
% $\lambda$ is a finite positive measure on $(U,\mcB(U))$ with full topological support.
%\end{ass}

We assume that the coefficients of the forward SDE satisfy the following conditions:
\begin{ass}\label{ass:onSDE}
%For any $t \in [0,T]$, $b\in U$, $e\in A$ and $x,x'\in\R^d$ we have:
\begin{enumerate}[i)]
  \item\label{ass:onSDE-Gamma} The functions $\jmp:[0,T]\times\R^d\times U\to\R^d$ and $\gamma:[0,T]\times\R^d\times A\to\R^d$ are jointly continuous and satisfy the growth condition
  \begin{align}\label{ekv:imp-bound}
   |x+\jmp(t,x,b)| \vee |x+\gamma(t,x,e)|\leq K_{\gamma,\jmp}\vee |x|
  \end{align}
  for some constant $K_{\gamma,\jmp}>0$ and the commutativity property \eqref{ekv:jump-commute}. Moreover, $\gamma$ is Lipschitz continuous in $x$ uniformly in $(t,e)$, \ie
  \begin{align*}
    |\gamma(t,x',e)-\gamma(t,x,e)|\leq k_\gamma|x'-x|
  \end{align*}
  for all $(t,x,x',e)\in [0,T]\times \R^d\times\R^d\times A$.
  \item\label{ass:onSDE-a-sigma} The coefficients $a:[0,T]\times\R^d\to\R^{d}$ and $\sigma:[0,T]\times\R^d\to\R^{d\times d}$ are measurable functions that satisfy the growth condition
  \begin{align*}
    |a(t,x)|+|\sigma(t,x)|&\leq C(1+|x|)
  \end{align*}
  and the Lipschitz continuity
  \begin{align*}
    |a(t,x)-a(t,x')|+|\sigma(t,x)-\sigma(t,x')|&\leq k_{a,\sigma}|x'-x|
  \end{align*}
  for all $(t,x,x')\in [0,T]\times \R^d\times\R^d$.
\end{enumerate}
\end{ass}

Moreover, we make the following assumptions on the coefficients in the cost functional $J$:

\begin{ass}\label{ass:oncoeff}
For some $\rho\geq 1$ we have:
\begin{enumerate}[i)]
  \item\label{ass:oncoeff-f} The running cost/reward $f:[0,T]\times \R^{d}\to\R$ is jointly continuous and satisfies the growth condition
  \begin{align*}
    |f(t,x)|\leq C(1+|x|^\rho).
  \end{align*}
  \item The terminal value $\psi:\R^d\to\R$ is continuous, of polynomial growth, \ie
  \begin{align*}
    |\psi(x)|\leq C(1+|x|^\rho)
  \end{align*}
  and satisfies
  \begin{align}\label{ekv:ass@end}
    \sup_{b\in U}\big\{\psi(x+\jmp(T,x,b))-\ell(T,x,b)\big\}\leq \psi(x) \leq \inf_{e\in A}\big\{\psi(x+\gamma(T,x,e))+\chi(T,x,e)\big\}
  \end{align}
  for all $x\in \R^d$.
  \item\label{ass:oncoeff-ell} The intervention costs $\ell:[0,T]\times\R^d\times U\to [\delta,\infty)$ and $\chi:[0,T]\times\R^d\times A\to [\delta,\infty)$ (with $\delta>0$) are jointly continuous, of polynomial growth, \ie
      \begin{align*}
        \ell(t,x,b)+\chi(t,x,e)\leq C(1+|x|^\rho),
      \end{align*}
      and satisfy \eqref{ekv:interv-cost-commute}.
\end{enumerate}
\end{ass}

Moreover, we let $\supOP v(T,x)=\infOP v(T,x)=v(T,x)$. We also need the following assumption, which is a variant of the no-free-loop condition commonly used in optimal switching problems (see \eg assumption H3.(ii) in \cite{Morlais13} or assumption H4 in \cite{BollanSWG1}), and was used to prove the comparison result for viscosity solutions to \eqref{ekv:var-ineq} in \cite{dbl-obst-qvi}:

\begin{ass}\label{ass:nofreeloop}
For each $R>0$, there are functions $h_1,h_2:[0,T]\to (0,\infty)$ such that for each $(t,x_0)\in [0,T]\times \bar B_R(\R^d)$ it holds that, whenever $(b_i)_{i\in\bbN}$, $(e_i)_{i\in\bbN}$ and $(\iota_i)_{i\in \bbN}$ are sequences in $U$, $A$ and $\{1,2\}$, respectively, for which there exists a $\kappa\in\bbN$ such that $|x_{\kappa}-x_0|\leq h_1(t)$, where
\begin{align*}
  x_{j}=x_{j-1}+\ett_{[\iota_j=1]}\jmp(t,x_{j-1},b_j)+\ett_{[\iota_j=2]}\gamma(t,x_{j-1},e_j),\quad \text{for }j=1,2,\ldots,\kappa
\end{align*}
then
\begin{align*}
  |\sum_{j=1}^{\kappa}(\ett_{[\iota_j=1]}\ell(t,x_{j-1},b_j)-\ett_{[\iota_j=2]}\chi(t,x_{j-1},e_j))|\geq h_2(t).
\end{align*}
\end{ass}
The assumption is intended to prevent the players from becoming trapped in an infinite cycle of simultaneous impulses.

\subsection{Notations}
We let $(\Omega,\mcF,\Prob)$ be a complete probability space on which lives a $d$-dimensional Brownian motion $W$. We denote by $\bbF:=(\mcF_t)_{0\leq t\leq T}$ the augmented natural filtration generated by $W$ and for $t\in[0,T]$ we let $\bbF^t:=(\mcF^t_s)_{t\leq s\leq T}$ denote the augmented natural filtration generated by $(W_s-W_t:t\leq s\leq T)$.\\

\noindent Throughout, we will use the following notations:
\begin{itemize}
\item We let $\Pi^g$ denote the set of all functions $\varphi:[0,T]\times\R^d\to\R$ that are of polynomial growth in $x$, \ie there are constants $C,\rho>0$ such that $|\varphi(t,x)|\leq C(1+|x|^\rho)$ for all $(t,x)\in [0,T]\times\R^d$, and let $\Pi^g_c$ be the subset of jointly continuous functions.
  \item Given a filtration $\bbG$, we let $\Pred(\bbG)$ denote the $\sigma$-algebra of $\bbG$-predictably measurable subsets of $[0,T]\times \Omega$.
  \item We let $\mcT$ be the set of all $[0,T]$-valued $\bbF$-stopping times and for each $\eta\in\mcT$, we let $\mcT_\eta$ be the subset of stopping times $\tau$ such that $\tau\geq \eta$, $\Prob$-a.s.
  \item For $t\in [0,T]$, we let $\bar\mcU_t$ (resp.~$\bar\mcA_t$) be the set of all $u=(\tau_j,\beta_j)_{j=1}^N$ (resp.~$\alpha=(\eta_j,\theta_j)_{j=1}^M$), where $(\tau_j)_{j=1}^\infty$ (resp.~$(\eta_j)_{j=1}^\infty$) is a non-decreasing sequence of $\bbF$-stopping times in $\mcT_t$, $\beta_j$ (resp.~$\theta_j$) is a $\mcF_{\tau_j}$-measurable (resp.~$\mcF_{\eta_j}$-measurable) r.v.~taking values in the compact set $U$ (resp.~$A$) and $N:=\max\{j:\tau_j<T\}$ (resp.~$M:=\max\{j:\eta_j<T\}$) is $\Prob$-a.s.~finite.
\end{itemize}
Whenever $(u,\alpha)\in \bar\mcU_t\times \bar\mcA_t$, an iteration of standard results implies the existence of a unique solution, $X^{t,x;u,\alpha}$ to \eqref{ekv:fwd-sde}. Moreover, Lemma~\ref{app:lem:SDEmoment} implies that the solution has moments of all orders and, thus, belongs to $\mcS^{p}_{t}$ for all $p\geq 1$. Since the cumulative intervention costs for Player 1 (resp.~Player 2),
\begin{align*}
  \Lambda^{t,x;u,\alpha}_s:=\sum_{j=1}^\infty\ett_{[\tau_j < s]}\ell(\tau_j,X^{t,x,[u]_{j-1},\alpha}_{\tau_j},\beta_j),\quad \forall s\in [0,T]
\end{align*}
(resp.
\begin{align*}
  \Theta^{t,x;u,\alpha}_s:=\sum_{j=1}^\infty\ett_{[\eta_j < s]}\chi(\eta_j,X^{t,x,u,[\alpha]_{j-1}}_{\eta_j},\theta_j),\quad \forall s\in [0,T]),
\end{align*}
may still be unbounded, we need to impose further restrictions on the control sets for the cost/reward functional $J(t,x;u,\alpha)$ to be well defined. We thus introduce the sets of impulse controls:
\begin{itemize}
  \item For $t\in [0,T]$, we let $\mcU_t$ (resp.~$\mcA_t$) be the set of all $u\in\bar\mcU_t$ (resp.~$\alpha\in\bar\mcA_t$) such that $\Lambda^{t,x;u,\bar\alpha}_T\in L^2(\Omega,\mcF_T,\Prob)$ for each $x\in\R^d$ and $\bar\alpha\in\bar\mcA_t$ (resp.~$\Theta^{t,x;\bar u,\alpha}_T\in L^2(\Omega,\mcF_T,\Prob)$ for each $x\in\R^d$ and $\bar u\in\bar\mcA_t$).
  \item For each $k\in\bbN$, we let $\mcU^{k}_t$ (resp.~$\mcA^{k}_t$) be the subset of $\mcU_t$ (resp.~$\mcA_t$) consisting of all impulse controls for which $N\leq k$ (resp.~$M\leq k$), $\Prob$-a.s.
  \item For $u\in\mcU$, we let $[u]_{j}:=(\tau_i,\beta_i)_{i=1}^{N\wedge j}$. Moreover, we introduce $N(s):=\max\{j\geq 0:\tau_j\leq s\}$ and let $u_s:=[u]_{N(s)}$. For $\alpha\in\mcA$, we let $[\alpha]_{j}:=(\eta_i,\theta_i)_{i=1}^{M\wedge j}$.
  \item We define the concatenation of impulse controls $u_1=(\tau^1_j,\beta^1_j)_{j=1}^{N_1}\in\mcU$ and $u_2=(\tau^2_j,\beta^2_j)_{j=1}^{N_2}\in\mcU$ at $t\in [0,T]$ as
      \begin{align*}
        u_1\oplus_t u_2:=((\tau^1_j,\beta^1_j),\ldots,(\tau^1_{N_1(t)},\beta^1_{N_1(t)}),(\tau^2_{N_2(t)+1},\beta^2_{N_2(t)+1}),\ldots, (\tau^2_{N_2},\beta^2_{N_2}))
      \end{align*}
      and similarly on $\mcA$.
\end{itemize}

Next, we define the set of \emph{non-anticipative impulse control strategies}.
\begin{defn}
For $t\in [0,T]$, we denote by $\mcU^S_t$ (resp.~$\mcA^{S}_t$) the set of all strategies $u^S:\mcA\to\mcU_t$ (resp.~$\alpha^S:\mcU\to\mcA_t$) such that $(u^S(\alpha))_r=(u^S(\tilde\alpha))_r$ whenever $\alpha,\tilde\alpha\in\mcA$ agree on $[0,r]$, \ie $\alpha_r=\tilde\alpha_r$, $\Prob$-a.s. (resp. $(\alpha^S(u))_r=(\alpha^S(\tilde u))_r$ whenever $u,\tilde u\in\mcU$ agree on $[0,r]$). Furthermore, for $k\in\bbN$, we let $\mcU^{S,k}$ (resp.~$\mcA^{S,k}_t$) be the subset of $\mcU^S_t$ (resp.~$\mcA^{S}_t$) containing all maps $u^S:\mcA\to\mcU^k_t$ (resp.~$\alpha^S:\mcU\to\mcA^k_t$).
\end{defn}

We also mention that, unless otherwise specified, all inequalities involving random variables are assumed to hold $\Prob$-a.s.

\subsection{Viscosity solutions}

We define the upper, $v^*$, and lower, $v_*$, semi-continuous envelope of a function $v:[0,T]\times\R^d\to\R$ as
\begin{align*}
v^*(t,x):=\limsup_{(t',x')\to(t,x),\,t'<T}v(t',x')\quad {\rm and}\quad v_*(t,x):=\liminf_{(t',x')\to(t,x),\,t'<T}v(t',x').
\end{align*}
Next, we introduce the limiting parabolic superjet $\bar J^+v$ and subjet $\bar J^-v$.
\begin{defn}\label{def:jets}
Subjets and superjets
\begin{enumerate}[i)]
\item For a l.s.c. (resp. u.s.c.) function $v : [0, T]\times \R^d \to \R$, the parabolic subjet, denote by $J^-v(t, x)$, (resp.~the parabolic superjet, $J^+v(t, x)$) of $v$ at $(t, x) \in [0, T]\times \R^d$, is defined as the set of triples $(p, q,M) \in \R \times\R^d \times \mathbb S^d$ satisfying
    \begin{align*}
      v(t', x') &\geq (\text{resp.~}\leq)\: v(t, x) + p(t'- t)+ < q, x' - x > + \tfrac{1}{2} < x' - x,M(x' - x) >
      \\
      &\quad+o(|t' - t| + |x' - x|^2)
    \end{align*}
    for all $(t',x')\in[0,T]\times\R^d$, where $\mathbb S^d$ is the set of symmetric real matrices of dimension $d\times d$.
\item For a l.s.c. (resp. u.s.c.) function $v : [0, T]\times \R^d \to \R$ we denote by $\bar J^-v(t, x)$ the parabolic limiting
subjet (resp. $\bar J^+v(t, x)$ the parabolic limiting superjet) of $v$ at $(t, x) \in [0, T]\times \R^d$, defined as the set of triples $(p, q,M) \in \R \times\R^d \times \mathbb S^d$ such that:
\begin{align*}
  (p, q,M) = \lim_{n\to\infty} (p_n, q_n,M_n),\quad (t, x) = \lim_{n\to\infty}(t_n, x_n)
\end{align*}
for some sequence $(t_n,x_n,p_n,q_n,M_n)_{n\geq 1}$ with $(p_n, q_n,M_n) \in J^-v(t_n, x_n)$ (resp. $(p_n, q_n,M_n) \in J^+v(t_n, x_n)$) for all $n\geq 1$ and $v(t, x) = \lim_{n\to\infty}v(t_n, x_n)$.
\end{enumerate}
\end{defn}

We now give the definition of a viscosity solution for the QVI in \eqref{ekv:var-ineq}. (see also pp. 9--10 of \cite{UsersGuide}).
\begin{defn}\label{def:visc-sol-jets}
Let $v$ be a locally bounded function from $[0,T]\times \R^d$ to $\R$. Then,
\begin{enumerate}[a)]
  \item It is referred to as a viscosity supersolution (resp. subsolution) to \eqref{ekv:var-ineq} if it is l.s.c.~(resp. u.s.c.) and satisfies:
  \begin{enumerate}[i)]
    \item $v(T,x)\geq \psi(x)$ (resp. $v(T,x)\leq \psi(x)$)
    \item For any $(t,x)\in [0,T)\times\R^d$ and $(p,q,X)\in \bar J^- v(t,x)$ (resp. $\bar J^+ v(t,x)$) we have
    \begin{align*}
      \min\big\{&v(t,x)-\supOP v(t,x),\max\{v(t,x)-\infOP v(t,x),-p-q^\top a(t,x)
      \\
      &-\frac{1}{2}\trace(\sigma\sigma^\top(t,x)X)-f(t,x)\}\big\}\geq 0\quad (\text{resp. }\leq 0)
    \end{align*}
  \end{enumerate}
  \item It is called a viscosity solution to \eqref{ekv:var-ineq} if $v_*$ is a supersolution and $v^*$ is a subsolution.
\end{enumerate}
\end{defn}

We will sometimes use the following alternative definition of viscosity supersolutions (resp. subsolutions):
\begin{defn}\label{def:visc-sol-dom}
  A l.s.c.~(resp. u.s.c.) function $v$ is a viscosity supersolution (resp.~subsolution) to \eqref{ekv:var-ineq} if $v(T,x)\leq \psi(x)$ (resp. $\geq \psi(x)$) and whenever $\varphi\in C^{1,2}([0,T]\times\R^d\to\R)$ is such that $\varphi(t,x)=v(t,x)$ and $\varphi-v$ has a local maximum (resp. minimum) at $(t,x)$, then
  \begin{align*}
    \min\big\{&v(t,x)-\supOP v(t,x),\max\{v(t,x)-\infOP v(t,x),-\varphi_t(t,x)-\mcL\varphi(t,x)
    \\
    &-f(t,x)\}\big\}\geq 0\quad(\text{resp. }\leq 0).
  \end{align*}
\end{defn}

%\begin{rem}\label{rem:mcM-monotone}
%Let $u,v:[0,T]\times \R^d\to\R$ be locally bounded functions. We remark that $\supOP$ and $\infOP$ are both monotone (if $u\leq v$ pointwise, then $\supOP u\leq \supOP v$ and $\infOP u\leq \infOP v$). Moreover, $\supOP(u_*)$ and $\infOP(u_*)$ (resp. $\supOP(u^*)$ and $\infOP(u^*)$) are l.s.c.~(resp. u.s.c.).% and $(\supOP u)_*\leq
%
%In particular, $\supOP v$ and $\infOP v$ are both jointly continuous whenever $v$ is.\\
%\end{rem}

\subsection{Main result\label{sec:zsg}}
The main result of the paper is the following:

\begin{thm}\label{thm:main}
The game has a value and $v=\underline v=\bar v$ belongs to $\Pi^g_c$ and is the unique (in $\Pi^g$) viscosity solution to \eqref{ekv:var-ineq}.
\end{thm}

Notable is that the setup is completely symmetric and from the proof of Theorem~\ref{thm:main} we find that $v$ is also a viscosity solution to
\begin{align}\label{ekv:var-ineq-max-min}
\begin{cases}
 \max\{v(t,x)-\infOP v(t,x),\min\{v(t,x)-\supOP v(t,x), -\frac{\partial}{\partial t}v(t,x)-\mcL v(t,x)\\
  \quad-f(t,x)\}\}=0,\quad\forall (t,x)\in[0,T)\times \R^d \\
  v(T,x)=\psi(x).
\end{cases}
\end{align}
In particular, this lead us to the following corollary result:

\begin{cor}
The unique viscosity solutions to \eqref{ekv:var-ineq} and \eqref{ekv:var-ineq-max-min} coincide.
\end{cor}

\section{Truncations of the game\label{sec:game-trunc}}
Our approach to prove Theorem~\ref{thm:main} is based on a truncation of the game, where we restrict the number of interventions that the players are allowed to make. Before we turn to these truncations we prove that the value functions are members of the set $\Pi^g$:

\begin{lem}
The value functions $\underline v$ and $\bar v$ can be chosen deterministic. Moreover, there is a $C>0$ such that
\begin{align}\label{ekv:vfs-bound}
|\underline v(t,x)| + |\bar v(t,x)| \leq C(1+|x|^\rho).
\end{align}
\end{lem}

\noindent\emph{Proof.} That the upper and lower value functions for a zero-sum game of classical controls can be chosen deterministic was proved in \cite{Buckdahn08} (see Proposition 4.1 therein). Appealing to an analogous argument shows that the same property holds for the value functions in the present setting as well. We turn to the matter of showing polynomial growth and have
\begin{align}
\inf_{\alpha\in\mcA}J(t,x;\emptyset,\alpha)\leq \underline v(t,x),\,\bar v(t,x)\leq \sup_{u\in\mcU}J(t,x;u,\emptyset),
\end{align}
where $\emptyset$ denotes the impulse control without any interventions. Hence,
\begin{align*}
|\underline v(t,x)| + |\bar v(t,x)| \leq 2\sup_{(u,\alpha)\in\mcU\times\mcA}\E\Big[|\psi(X^{t,x;u,\alpha}_T)|+\int_t^{T}|f(s,X^{t,x;u,\alpha}_s)|ds\Big]\leq C(1+|x|^\rho)
\end{align*}
by the polynomial growth assumptions on $\psi$ and $f$ and Lemma~\ref{app:lem:SDEmoment}.\qed\\

Next, we develop useful approximations of the value functions by their counterparts with a truncated number of interventions. Similarly to the above, there are for each $k,l\in\bbN$, maps $\underline v^{k,l},\bar v^{k,l}\in\Pi^g$ such that
\begin{align}\label{ekv:lvf-k-def}
\underline v^{k,l}(t,x)=\essinf_{\alpha^S\in\mcA^{S,l}_t}\,\esssup_{u\in\mcU_t^k}J(t,x;u,\alpha^S(u))
\end{align}
and
\begin{align}\label{ekv:uvf-k-def}
\bar v^{k,l}(t,x)=\esssup_{u^S\in\mcU^{S,k}_t}\,\essinf_{\alpha\in\mcA^l_t}J(t,x;u^S(\alpha),\alpha)
\end{align}
for all $(t,x)\in[0,T]\times\R^d$. In addition, we let $\underline v^{\infty,l}(t,x)$, $\underline v^{k,\infty}(t,x)$, $\bar v^{\infty,l}(t,x)$ and $\bar v^{k,\infty}(t,x)$ (all elements of $\Pi^g$) denote the corresponding values when only one of the players has a restriction on the number of interventions; for example,
\begin{align*}
\bar v^{\infty,l}(t,x)=\esssup_{u^S\in\mcU^{S}_t}\,\essinf_{\alpha\in\mcA^l_t}J(t,x;u^S(\alpha),\alpha).
\end{align*}
Clearly, we have
\begin{align}\label{ekv:uvf-bnd}
\bar v^{k,\infty}(t,x)\leq \bar v(t,x)\leq \bar v^{\infty,l}(t,x)
\end{align}
and
\begin{align}\label{ekv:lvf-bnd}
\underline v^{k,\infty}(t,x)\leq \underline v(t,x)\leq \underline v^{\infty,l}(t,x).
\end{align}
Moreover, the value functions are related through the following convergence result:

\begin{prop}\label{prop:trunk-vfs-approx}
For each $l\in\bbN$, there is a $C>0$ such that
\begin{align}\label{ekv:lvf-approx}
\underline v^{\infty,l}(t,x)-\underline v^{k,l}(t,x) +\bar v^{\infty,l}(t,x)-\bar v^{k,l}(t,x) \leq \frac{C}{\sqrt{k}}(1+|x|^{3\rho/2})
\end{align}
for all $k\geq 0$. Similarly, for each $k\in\bbN$, there is a $C>0$ such that
\begin{align}\label{ekv:uvf-approx}
\underline v^{k,l}(t,x)-\underline v^{k,\infty}(t,x) +\bar v^{k,l}(t,x)-\bar v^{k,\infty}(t,x) \leq \frac{C}{\sqrt{l}}(1+|x|^{3\rho/2})
\end{align}
for all $l\geq 0$.
\end{prop}

\noindent\emph{Proof.} We only prove \eqref{ekv:lvf-approx} as \eqref{ekv:uvf-approx} follows by an identical argument. For any $\alpha^S\in\mcA^{S,l}$ and $u\in\mcU_t$, we have
\begin{align*}
J(t,x;u,\alpha^S(u))&\leq \E\Big[ |\psi(X^{t,x;u,\alpha^S(u)})+|\int_t^{T} |f(r,X^{t,x;u,\alpha^S(u)}_r)|dr+ \sum_{j=1}^M\chi(\eta_j,X^{t,x;u,[\alpha]_{j-1}},\theta_j) +N\delta\Big|\mcF_t\Big]
\\
&\leq C(1+|x|^\rho)-\delta\E\big[N\big|\mcF_t\big],
\end{align*}
where the constant $C>0$ may depend on $l$ but can be chosen independent of $\alpha^S$ and $u$. There is thus a $C_1>0$ (that does not depend on $\alpha^S$ and $u$) such that the control $u\in\mcU_t$ is dominated by the control $\emptyset\in\mcU_t$ (\ie the impulse control that makes no interventions in $[0,T)$) whenever
\begin{align*}
\E\big[N\big|\mcF_t\big] > C_1(1+|x|^\rho).
\end{align*}
We let $\mcU^*_t$ be the subset of all $u=(\tau_j,\beta_j)_{j=1}^N$ such that
\begin{align*}
\E\big[N\big|\mcF_t\big] \leq C_1(1+|x|^\rho)
\end{align*}
and have that
\begin{align*}
  \esssup_{u\in\mcU_t}J(t,x;u,\alpha^S(u))=\esssup_{u\in\mcU^*_t}J(t,x;u,\alpha^S(u))
\end{align*}
for any $\alpha^S\in\mcA^{S,l}_t$. Hence,
\begin{align*}
&\esssup_{u\in\mcU_t}J(t,x;u,\alpha^S(u))- \esssup_{u\in\mcU^k_t}J(t,x;u,\alpha^S(u))\leq \esssup_{u\in\mcU^*_t}(J(t,x;u,\alpha^S(u))-J(t,x;[u]_k,\alpha^S([u]_k)))
\\
&\leq \esssup_{u\in\mcU^*_t}\E\Big[ \ett_{[N>k]}\big(|\psi(X^{t,x;u,\alpha^S(u)}_T)|+|\psi(X^{t,x;[u]_k,\alpha^S([u]_k)}_T)|+
\\
&\quad \int_t^T(|f(r,X^{t,x;u,\alpha^S(u)}_r)|+|f(r,X^{t,x;[u]_k,\alpha^S([u]_k)}_r)|)dr\big)\Big|\mcF_t\Big]
\\
&\leq C\esssup_{u\in\mcU^*_t}\E\Big[\ett_{[N>k]} (1+\sup_{s\in[t,T]}|X^{t,x;u,\alpha^S(u)}|^{2\rho}+\sup_{s\in[t,T]}|X^{t,x;[u]_k,\alpha^S([u]_k)}|^{2\rho})\Big|\mcF_t\Big]
\\
&\leq C\esssup_{u\in\mcU^*_t}\E\big[\ett_{[N>k]}\big]^{1/2} \E\Big[1+\sup_{s\in[t,T]}|X^{t,x;u,\alpha^S(u)}|^\rho+\sup_{s\in[t,T]}|X^{t,x;[u]_k,\alpha^S([u]_k)}|^\rho\Big|\mcF_t\Big]^{1/2}
\\
& \leq C\sqrt{\frac{C_1(1+|x|^\rho)}{k}}(1+|x|^\rho),
\end{align*}
where we use Lemma~\ref{app:lem:SDEmoment} along with the bound on $\E\big[N\big|\mcF_t\big]$ that holds for all $u\in \mcU^*_t$ to reach the last inequality. Since $\alpha^S\in\mcA^{S,l}_t$ was arbitrary this proves the existence of a $C>0$ such that \eqref{ekv:lvf-approx} holds for the first term. We may address the second term similarly by noting that
\begin{align*}
  \bar v^{\infty,l}(t,x)=\esssup_{u^S\in\mcU^{S,*}_t}\,\essinf_{\alpha\in\mcA^l_t}J(t,x;u^S(\alpha),\alpha),
\end{align*}
where $\mcU^{S,*}_t$ is the subset of all $u^S=(\tau^S_j,\beta^S_j)_{j=1}^{N^S}\in\mcU^{S}_t$ for which
\begin{align*}
\E\big[N^S(\alpha)\big|\mcF_t\big] \leq C_1(1+|x|^\rho),
\end{align*}
for all $\alpha\in\mcA^l$. Then,
\begin{align*}
&\esssup_{u^S\in\mcU^{S}_t}\,\essinf_{\alpha\in\mcA^l_t}J(t,x;u^S(\alpha),\alpha)- \esssup_{u^S\in\mcU^{S,k}_t}\,\essinf_{\alpha\in\mcA^l_t}J(t,x;u^S(\alpha),\alpha)
\\
&\leq \esssup_{u\in\mcU^*_t}\esssup_{\alpha\in\mcA^l_t}(J(t,x;u,\alpha)-J(t,x;[u]_k,\alpha))
\end{align*}
and the result follows by repeating the above steps.\qed\\

%\begin{rem}\label{rem:SDEflow}
%Letting $\Lambda(u,\tilde u):=[0,T]\setminus\bigcup_{j=1}^k[\underline\eta_j,\bar\eta_j)$ with $\underline\eta_j:=\eta_j\wedge\tilde\eta_j$ and $\bar\eta_j:=\eta_j\vee\tilde\eta_j$ for $j=1,\ldots,k$, we note that $\sup_{s\in \Lambda(u,\tilde u)}|X^{t,u}_s-X^{t,\tilde u}_s|\leq\bold d[(T,X^{t,u}),(T,X^{t,\tilde u})]$. Lemma~\ref{app:lem:SDEflow} thus gives an effective bound on the expected difference between the paths of the state process under different controls.
%\end{rem}

%%%%%%%%%%%%%%%%%%%%%%%%%%%%%%%%%%%%%%%%%%%%%%%%%%%%%%%%%%%%%%%%%%%%%%%%%%%%%%%%%%%%%%%%%%%%%%%%%%%%%%%%%%%%%%%%%%%%%%%%%%%%%%%%%%%%%%%%%%%%%%%%%%
%%%%%%%%%%%%%%%%%%%%%%%%%%%%%%%%%%%%%%%%%%%%%%%%%%%%%%%%%%%%%%%%%%%%%%%%%%%%%%%%%%%%%%%%%%%%%%%%%%%%%%%%%%%%%%%%%%%%%%%%%%%%%%%%%%%%%%%%%%%%%%%%%%
%%%%%%%%%%%%%%%%%%%%%%%%%%%%%%%%%%%%%%%%%%%%%%%%%%%%%%%%%%%%%%%%%%%%%%%%%%%%%%%%%%%%%%%%%%%%%%%%%%%%%%%%%%%%%%%%%%%%%%%%%%%%%%%%%%%%%%%%%%%%%%%%%%

\subsection{A discretization of the control sets}
Our approach requires a discretization of the control sets and we introduce the following sets:
\begin{defn}\label{def:discretization}
For $\eps>0$:
\begin{itemize}
  \item We let $n^\eps\geq 0$ be the smallest integer such that $2^{-n^\eps}T\leq\eps$, set $n_{\bbT}^\eps:=2^{n^\eps}+1$ and introduce the discrete set $\bbT^\eps:=\{t^\eps_i:t^\eps_i=(i-1)2^{-\iota}T,i=1,\ldots,n_{\bbT}^\eps\}$, a discritization of $[0,T]$ with step-size $\Delta t^\eps:=2^{-n^\eps}T$. For $s\in [0,T]$, we let $\bbT^\eps_s:=\bbT^\eps\cap[s,T]$.
  \item We let $(A^\eps_{i})_{i=1}^{n^\eps_A}$ be a Borel-partition of $A$ such that each $A^\eps_i$ has non-empty interior in $A$ and a diameter that does not exceed $\eps$ and let $(e^\eps_i)_{i=1}^{n^\eps_A}$ be a sequence with $e^\eps_i\in \text{int}\,A^\eps_i$ and denote by $\bar A^\eps:=\{e^\eps_1,\ldots,e^\eps_{n^\eps_A}\}$ the corresponding discritization of $A$.
  \item We let $(U^\eps_{i})_{i=1}^{n^\eps_U}$ be a Borel-partition of $U$ such that each $U^\eps_i$ has non-empty interior in $U$ and a diameter that does not exceed $\eps$ and let $(b^\eps_i)_{i=1}^{n^\eps_U}$ be a sequence with $b^\eps_i\in \text{int}\,U^\eps_i$ and denote by $\bar U^\eps:=\{b^\eps_1,\ldots,b^\eps_{n^\eps_U}\}$ the corresponding discritization of $U$.
\end{itemize}
\end{defn}
For $\eps>0$ and $k\in\bbN$, we then let
\begin{align*}
\mcA^{k,\eps}_t:=\{\alpha\in\mcA^{k}_t: \eta_j \in \bbT^\eps,\:\theta_j\in  \bar A^\eps,\:\Prob\text{-a.s.},\text{ for }j=1,\ldots,k\}
\end{align*}
(resp.
\begin{align*}
\mcU^{k,\eps}_t:=\{u\in\mcU^{k}_t: \tau_j \in \bbT^\eps,\:\beta_j\in  \bar U^\eps,\:\Prob\text{-a.s.},\text{ for }j=1,\ldots,k\}).
\end{align*}
For $\alpha\in\mcA^{k}_t$ (resp. $u\in\mcU^{k}_t$), we let $\pi^{\eps}_A(\alpha)=(\eta^\eps_j,\theta^\eps_j)_{j=1}^k:=\big(\pi_1^\eps(\eta_j), \pi^{A,\eps}_2(\theta_j)\big)_{j=1}^k$ (resp. $\pi^{\eps}_U(u):=\big(\pi_1^\eps(\tau_j), \pi^{U,\eps}_2(\beta_j)\big)_{j=1}^k$) with $\pi_1^\eps(s):= \min\{r\in \bbT^\eps:r\geq s\}$ and $\pi^{A,\eps}_2(e):=\sum_{i=1}^{n^\eps_A}e^\eps_i\ett_{A^\eps_i}(e)$ (resp. $\pi^{U,\eps}_2(b):=\sum_{i=1}^{n^\eps_U}b^\eps_i\ett_{U^\eps_i}(b)$) and note that $\pi^{\eps}_A(\alpha)\in\mcA^{k,\eps}_t$ (resp. $\pi^{\eps}_U(u)\in\mcU^{k,\eps}_t$).

\begin{defn}
For $t\in[0,T]$, $\eps>0$ and $k\in\bbN$, we let $\mcA^{S,k,\eps}_t$ (resp.~$\mcU^{S,k,\eps}_t$) be the subset of $\mcA^{S,k}_t$ (resp.~$\mcU^{S,k}_t$) of strategies $\alpha^S:\mcU_t\to\mcA^{k,\eps}_t$ (resp.~$u^S:\mcA_t\to\mcU^{k,\eps}_t$).

%Moreover, we let $\mcA^{S,k,\eps}_t$ (resp.~$\mcU^{S,k,\eps}_t$) bet the subset of $\bar\mcA^{S,k,\eps}_t$ (resp.~$\bar\mcU^{S,k,\eps}_t$) such that on the set $\{\omega:\pi_U^\eps(u)_{t^\eps_i}=\pi_U^\eps(\tilde u)_{t^\eps_i}\}$ (resp.~$\{\omega:\pi_A^\eps(\alpha)_{t^\eps_i}=\pi_A^\eps(\tilde \alpha)_{t^\eps_i}\}$) we have $(\alpha^S(u))_{t^\eps_i}=(\alpha^S(\tilde u))_{t^\eps_i}$ (resp.~$(u^S(\alpha))_{t^\eps_i}=(u^S(\tilde \alpha))_{t^\eps_i}$) outside of a $\Prob$-null set for $i=1,\ldots,n_\bbT^\eps$. (Introducera nedan istället!!!!)
\end{defn}

For $\eps > 0$ and $k,l\in\bbN$, we introduce the discretized value functions $\underline v^{k,l}_{\eps},\bar v^{k,l}_{\eps}\in\Pi^g$ defined as
\begin{align*}
  \underline v^{k,l}_{\eps}(t,x)&:=\essinf_{\alpha^S\in\mcA^{S,l,\eps}_t}\,\esssup_{u\in\mcU^{k,\eps}_t}J(t,x;u,\alpha^S(u))
\end{align*}
and
\begin{align*}
  \bar v^{k,l}_{\eps}(t,x)&:=\esssup_{u^S\in\mcU^{S,k,\eps}_t}\,\essinf_{\alpha\in\mcA^{l,\eps}_t}J(t,x;u^S(\alpha),\alpha).
\end{align*}
%using the convention that $\bbT^0:=[0,T]$, $U^0:=U$ and $A^0:=A$.

The main result of the present section is:
\begin{prop}\label{prop:disc-vfs}
For each $k,l\in\bbN$ we have
\begin{align}\label{ekv:disc-vfs}
  \lim_{\eps\searrow 0}\underline v^{k,l}_{\eps}(t,x)=\underline v^{k,l}(t,x)\leq \bar v^{k,l}(t,x)=\lim_{\eps\searrow 0}\bar v^{k,l}_{\eps}(t,x)
\end{align}
for all $(t,x)\in [0,T]\times\R^d$.
\end{prop}

The proof is distributed over a sequence of lemmas where we start by considering the value function
\begin{align*}
  \underline v^{k,l,\eps,0}(t,x)&:=\essinf_{\alpha^S\in\mcA^{S,l,\eps}_t}\,\esssup_{u\in\mcU^{k}_t}J(t,x;u,\alpha^S(u)),
\end{align*}
in which the discretization only affects the strategy of Player 2. Clearly, we have $\underline v^{k,l,\eps,0}(t,x)\geq \underline v^{k,l}(t,x)$. On the other hand,
\begin{align*}
  \underline v^{k,l,\eps,0}(t,x)=\essinf_{\alpha^S\in\mcA^{S,l}_t}\,\esssup_{u\in\mcU^{k}_t}J(t,x;u,\pi^\eps_A(\alpha^S(u)))
\end{align*}
and by Lemma~\ref{app:lem:SDEflow} it follows that $\lim_{\eps\searrow 0}\underline v^{k,l,\eps,0}(t,x) = \underline v^{k,l}(t,x)$.

\begin{lem}
For each $(t,x)\in [0,T]\times\R^d$ and $k,l\in\bbN$, it holds that
\begin{align*}
  \lim_{\eps\searrow 0}\underline v^{k,l}_{\eps}(t,x)=\lim_{\eps\searrow 0}\underline v^{k,l,\eps,0}(t,x)
\end{align*}
\end{lem}

\noindent\emph{Proof.} By definition $\underline v^{k,l}_{\eps}(t,x)\leq \underline v^{k,l,\eps,0}(t,x)$. On the other hand,
\begin{align*}
  \underline v^{k,l,\eps,0}(t,x)-\underline v^{k,l}_{\eps}(t,x)&=\essinf_{\alpha^S\in\mcA^{S,l,\eps}_t}\, \esssup_{u\in\mcU^k_t}J(t,x;u,\alpha^S(u)) - \essinf_{\alpha^S\in\mcA^{S,l,\eps}_t}\,\esssup_{u\in\mcU^{k,\eps}_t}J(t,x;u,\alpha^S(u))
  \\
  &\leq\essinf_{\alpha^S\in\mcA^{S,l,\eps}_t}\, \esssup_{u\in\mcU^k_t}J(t,x;u,\alpha^S(\pi^\eps_U(u))) - \essinf_{\alpha^S\in\mcA^{S,l,\eps}_t}\,\esssup_{u\in\mcU^{k,\eps}_t}J(t,x;u,\alpha^S(u))
  \\
  &\leq\esssup_{\alpha^S\in\mcA^{S,l,\eps}_t}\, \esssup_{u\in\mcU^k_t}\{J(t,x;u,\alpha^S(\pi^\eps_U(u))) - J(t,x;\pi^\eps_U(u),\alpha^S(\pi^\eps_U(u)))\}
  \\
  &\leq\esssup_{\alpha\in\mcA^{l,\eps}_t} \,\esssup_{u\in\mcU^k_t}|J(t,x;u,\alpha) - J(t,x;\pi^\eps_U(u),\alpha)|,
\end{align*}
where the last term tends to 0 in $L^1$ as $\eps\searrow 0$ by Lemma~\ref{app:lem:SDEflow}.\qed\\

The above lemma gives the first equality in \eqref{ekv:disc-vfs}. To derive the middle inequality we need the following lemma:

\begin{lem}\label{lem:DynP-eps}
The map $\underline v^{k,l}_{\eps}:[0,T]\times\R^d\to\R$ satisfies the dynamic programming relation
\begin{align}\label{ekv:dynP-eps}
  \begin{cases}
    \underline v^{k,l}_{\eps}(T,x)=\psi(x),\forall x\in\R^d,
    \\
    \underline v^{k,l}_{\eps}(t^\eps_i,x)
    =\supOP^k_\eps\infOP^l_\eps\underline v^{\cdot,\cdot}_{\eps}(t^\eps_i+,x) ,\quad\forall (i,x)\in \{1,\ldots,n^\eps_\bbT-1\}\times\R^d
    \\
    \underline v^{k,l}_{\eps}(t,x)
    =\E\Big[\int_{t}^{\pi^\eps_1(t)}f(s,X^{t,x,\emptyset,\emptyset}_s)ds + \underline v^{k,l}_{\eps}(\pi^\eps_1(t),X^{t,x,\emptyset,\emptyset}_{\pi^\eps_1(t)})\Big]\quad\forall (t,x)\in ([0,T)\setminus\bbT^\eps)\times\R^d,
  \end{cases}
\end{align}
where for any sequence $w=(w^{k,l}:[0,T]\times\R^d\to\R)_{0\leq k,l}$, the operators $\supOP^k_\eps$ and $\infOP^l_\eps$ are defined recursively through
\begin{align*}
  \infOP^l_\eps w^{k,\cdot}(t,x):=w^{k,l}(t,x)\wedge \inf_{e\in \bar A^\eps}\{\infOP^{l-1}_\eps w^{k,\cdot}(t,x+\gamma(t,x,e))+\chi(t,x,e)\}
\end{align*}
and
\begin{align*}
  \supOP^k_\eps w^{\cdot,l}(t,x):=w^{k,l}(t,x)\vee \sup_{b\in \bar U^\eps}\{\supOP^{k-1}_\eps w^{\cdot,l}(t,x+\jmp(t,x,b))-\ell(t,x,e)\}
\end{align*}
for $(t,x)\in [0,T]\times\R^d$ and $k,l\geq 1$, with $\infOP^0_\eps w^{k,\cdot}(t,x):= w^{k,0}(t,x)$ and $\supOP^0_\eps w^{\cdot,l}(t,x):= w^{0,l}(t,x)$.
\end{lem}

\noindent\emph{Proof.} Let $w=(w^{k,l}:[0,T]\times\R^d\to\R)_{0\leq k,l}$ be the unique solution to the recursion \eqref{ekv:dynP-eps} and note that for each $k,l\in\bbN$, the map $w^{k,l}:(t^\eps_i,t^\eps_{i+1}]\times\R^d\to \R$ is continuous for $i=1,\ldots,n_\bbT^\eps$. Moreover, the map $x\to w^{k,l}(t^\eps_i+,x)$ is continuous and, therefor, the map $x\mapsto \mcM^l\underline v^{k,\cdot}_{\eps}(t^\eps_{i}+,x)$ is continuous for each $k,l\in\bbN$ and $i=1,\ldots,n_\bbT^\eps$. The proof is based on induction and we assume that $w^{k,l}(t,x) = \underline v^{k,l}_{\eps}(t,x)$ for all $(t,x)\in [t^\eps_{i'+1},T]\times \R^d$ for some $i'\in\{1,3,\ldots,n^\eps_\bbT-1\}$. We pick arbitrary $x\in \R^d$ and prove that $w^{k,l}(t^\eps_{i'},x)=\underline v^{k,l}_{\eps}(t^\eps_{i'},x)$ over the two steps below. The result for arbitrary $t\in (t^\eps_{i'},t^\eps_{i'+1})$ then follows by a simpler argument. To ease notation we borrow from \cite{LiPeng09} and introduce the family of operators $(G^{t,x;u,\alpha}_{s,s'}:0\leq t\leq s\leq s'\leq T, (x,u,\alpha)\in\R^d\times\mcU_s\times\mcA_s)$, where $G^{t,x;u,\alpha}_{s,s'}:L^1(\mcF_{s'},\Prob)\to L^1(\mcF_s,\Prob)$ is defined as:
\begin{align*}
  G^{t,x;u,\alpha}_{s,s'}[\vartheta]&=\E\Big[\vartheta+\int_s^{s'}f(r,X^{t,x;u,\alpha}_r)dr - \sum_{j=1}^N \ell(\tau_j,X^{t,x;[u]_{j-1},\alpha},\beta_j) + \sum_{j=1}^M\chi(\eta_j,X^{t,x;u,[\alpha]_{j-1}},\theta_j)\Big|\mcF_s\Big].
\end{align*}

\textbf{Step 1:} We first prove that $w^{k,l}(t^\eps_{i'},x)\leq \underline v^{k,l}_{\eps}(t^\eps_{i'},x)$, where we start off by noting that (with $\mcU^k_{\{t\}}$, $\mcA^l_{\{t\}}$, etcetera, denoting the corresponding set of impulse controls where impulses are only allowed to take place at time $t$)
\begin{align*}
  w^{k,l}(t^\eps_{i'},x)&=\essinf_{\alpha^S\in\mcA^{S,l,\eps}_{\{t^\eps_{i'}\}}}\esssup_{u\in\mcU^{k,\eps}_{\{t^\eps_{i'}\}}} \E\Big[\int_{{t^\eps_{i'}}}^{t^\eps_{i'+1}}f(s,X^{t^\eps_{i'},x,u,\alpha^S(u)}_s)ds - \sum_{j=1}^N\ell(t^\eps_{i'},X^{t^\eps_{i'},x;[u]_{j-1},\alpha^S(u)}_s,\beta_j)
  \\
  &\quad + \sum_{j=1}^{M^S(u)}\jmp(t^\eps_{i'},X^{t^\eps_{i'},x;u,[\alpha^S(u)]_{j-1}},\theta^S_j(u)) + w^{k-N,l-M^S(u)}(t^\eps_{i'+1},X^{t^\eps_{i'},x,u,\alpha^{S}(u)}_{t^\eps_{i'+1}})\Big|\mcF_{t^\eps_{i'}}\Big]
  \\
  &=\essinf_{\alpha^S\in\mcA^{S,l,\eps}_{\{t^\eps_{i'}\}}}\esssup_{u\in\mcU^{k,\eps}_{\{t^\eps_{i'}\}}} G_{t^\eps_{i'},t^\eps_{i'+1}}^{t^\eps_{i'},x;u,\alpha^S(u)}[ w^{k-N,l-M^S(u)}(t^\eps_{i'+1},X^{t^\eps_{i'},x,u,\alpha^{S}(u)}_{t^\eps_{i'+1}})].
\end{align*}
We pick an arbitrary $\alpha^S\in\mcA^{S,l,\eps}_{t^\eps_{i'}}$ and note that we can define the restriction, $\alpha_1^S=(\eta^S_{1,j},\theta^S_{1,j})_{j=1}^{M^S_1}$, of $\alpha^S$ to $\mcA^{S,l,\eps}_{\{t^\eps_{i'}\}}$ as
\begin{align*}
  \alpha^S_1(u_1):=\alpha^S(u_1)_{t^\eps_{i'}},\qquad \forall u_1\in\mcU_{\{t^\eps_{i'}\}}.
\end{align*}
%Moreover, we can define the restriction of $\alpha^S$ to $\mcA_{t+h}$ as $\alpha^S_2(u)=\alpha(u)|_{[t+h,T]}$.
We fix $\varrho>0$ and have by stability under pasting\footnote{We can paste together two controls $u_1,u_2\in\mcU^k_s$ on sets $B_1\in\mcF_s$ and $B_2=B_1^c$ by setting $u=\ett_{B_1}u_1+\ett_{B_2}u_2\in\mcU^k_s$ and get by linearity of the conditional expectation that $G_{s,r}^{t^\eps_{i'},x;u,\alpha}[\vartheta]=\ett_{B_1}G_{s,r}^{t^\eps_{i'},x;u_1,\alpha}[\vartheta]+\ett_{B_2}G_{s,r}^{t^\eps_{i'},x;u_2,\alpha}[\vartheta]$.}
that there is a $u^\varrho_{1}=(\tau^{1,\varrho}_j,\beta^{1,\varrho}_j)_{1\leq j\leq N^\varrho_1}\in \mcU^k_{\{t^\eps_{i'}\}}$ such that
\begin{align}\nonumber
w^{k,l}(t^\eps_{i'},x)&\leq\esssup_{u\in\mcU^k_{\{t^\eps_{i'}\}}}G_{t^\eps_{i'},t^\eps_{i'+1}}^{t^\eps_{i'},x;u,\alpha_1(u)} [w^{k-N,l-M^S_1(u)}(t^\eps_{i'+1},X^{t^\eps_{i'},x;u_1,\alpha_1(u_1)}_{t^\eps_{i'+1}})]
\\
&\leq G_{t^\eps_{i'},t^\eps_{i'+1}}^{t^\eps_{i'},x;u^\varrho_1,\alpha_1(u^\varrho_1)}[w^{k-N^\varrho_1,l-M^S_1(u^\varrho_1)}(t^\eps_{i'+1}, X^{t^\eps_{i'},x;u^\varrho_1,\alpha^S_1(u^\varrho_1)}_{t^\eps_{i'}})]+\varrho.\label{ekv:w1}
\end{align}
Now, given $u^\varrho_{1}$ we can define the restriction, $\alpha^S_2$, of $\alpha^S$ to $\mcA_{t^\eps_{i'+1}}$ as
\begin{align*}
  \alpha^S_2(u_2):=(\eta^S_j(u_1^\varrho\oplus_{t^\eps_{i'}} u_2),\theta^S_j(u_1^\varrho\oplus_{t^\eps_{i'}} u_2))_{j=M^S_1(u_1^\varrho)+1}^{M^S(u_1^\varrho\oplus_{t^\eps_{i'}} u_2)},\qquad \forall u_2\in\mcU_{t^\eps_{i'}}.
\end{align*}
We let $(\mcO_\varpi)_{\varpi\geq 0}\subset\mcB(\R^n)$ be a partition of $\R^n$ with $\mcO_0:=\{x\in\R^d:|x|\geq K\}$ and $K>0$ chosen such that
\begin{align*}
  \Prob[\{\omega\in\Omega:|X^{t^\eps_{i'},x,\alpha,u}_{t^\eps_{i'+1}}|\geq K\}]\leq \varrho
\end{align*}
for all $(u,\alpha)\in \mcU_t\times\mcA_t$ and $(\mcO_\varpi)_{\varpi\geq 1}$ chosen such that
\begin{align*}
  |w^{k',l'}(t^\eps_{i'+1},x)-w^{k',l'}(t^\eps_{i'+1},x')|\leq \varrho
\end{align*}
and
\begin{align*}
  |J(t^\eps_{i'+1},x,u,\alpha)-J(t^\eps_{i'+1},x',u,\alpha))|\leq \varrho
\end{align*}
for all $\varpi\geq 1$, $k'\in \{0,\ldots,k\}$, $l'\in \{0,\ldots,l\}$, $(u,\alpha)\in \mcU^k_t\times\mcA^l_t$ and $x,x'\in\mcO_\varpi$. We pick $x_\varpi\in\mcO_\varpi$ and have by our induction assumption and the same pasting property as above that there is for each $\varpi\geq 1$, $k'\in \{0,\ldots,k\}$ and $l'\in \{0,\ldots,l\}$, a $u^\varrho_{2,\varpi,k',l'}\in\mcU_{t^\eps_{i'+1}}^{k'}$ such that
\begin{align*}
  w^{k',l'}(t^\eps_{i'+1},x_\varpi) &\leq J(t^\eps_{i'+1},x_\varpi,u^\varrho_{2,\varpi,k',l'},\alpha^S_2(u^\varrho_{2,\varpi,k',l'}))+\varrho.
\end{align*}
Consequently, whenever $X^{t^\eps_{i'},x;u^\varrho_1,\alpha_1^S(u^\varrho_1)}_{t^\eps_{i'+1}}\notin \mcO_0$, we have
\begin{align*}
&w^{k-N^\varrho_1,l-M^S_1(u^\eps_1)}(t^\eps_{i'+1},X^{t^\eps_{i'},x;u^\varrho_1,\alpha_1^S(u^\varrho_1)}_{t^\eps_{i'+1}})
\\
&\leq \sum_{\varpi\geq 1}\ett_{[X^{t^\eps_{i'},x;u^\varrho_1,\alpha_1^S(u^\varrho_1)}_{t^\eps_{i'+1}}\in \mcO_\varpi]}\sum_{k'=0}^k\ett_{[k-N^\varrho_1=k']}\sum_{l'=0}^l\ett_{[l-M^S_1(u^\varrho_1)=l']}w^{k',l'}_{\eps}(t^\eps_{i'+1},x_\varpi)+\varrho
\\
&\leq\sum_{\varpi\geq 1}\ett_{[X^{t^\eps_{i'},x;u^\varrho_1,\alpha_1^S(u^\varrho_1)}_{t^\eps_{i'+1}}\in \mcO_\varpi]}\sum_{k'=0}^k\ett_{[k-N^\varrho_1=k']}\sum_{l'=0}^l\ett_{[l-M^S_1(u^\varrho_1)=l']} J(t^\eps_{i'+1},x_\varpi,u^\varrho_{2,\varpi,k',l'},\alpha^S_2(u^\varrho_{2,\varpi,k',l'}))+2\varrho
\\
&\leq\sum_{\varpi\geq 1}\ett_{[X^{t^\eps_{i'},x;u^\varrho_1,\alpha_1^S(u^\varrho_1)}_{t^\eps_{i'+1}}\in \mcO_\varpi]}\sum_{k'=0}^k\ett_{[k-N^\varrho_1=k']}\sum_{l'=0}^l\ett_{[l-M^S_1(u^\varrho_1)=l']} J(t^\eps_{i'+1},X^{t^\eps_{i'},x;u_1^\varrho,\alpha_1^S(u_1^\varrho)}_{t^\eps_{i'+1}},u^\varrho_{2,\varpi,k',l'},\alpha^S(u^\varrho))+3\varrho,
\end{align*}
with
\begin{align*}
  u^\varrho:=u_1^\varrho\oplus_{t^\eps_{i'}}\sum_{\varpi\geq 1}\ett_{[X^{t^\eps_{i'},x;u_1^\varrho,\alpha_1^S(u_1^\varrho)}_{t^\eps_{i'+1}}\in\mcO_\varpi]}\ \sum_{k'=0}^k\ett_{[k-N^\varrho_1=k']}\sum_{l'=0}^l\ett_{[l-M^S_1(u^\varrho_1)=l']}u_{2,\varpi,k',l'}.
\end{align*}
Plugging this into \eqref{ekv:w1} we get
\begin{align*}
w^{k,l}(t^\eps_{i'},x) &\leq G_{t^\eps_{i'},t^\eps_{i'+1}}^{t^\eps_{i'},x;u^\varrho_1,\alpha_1^S(u^\varrho_1)}[\sum_{\varpi\geq 1}\ett_{[X^{t^\eps_{i'},x;u^\varrho_1,\alpha_1^S(u^\varrho_1)}_{t^\eps_{i'+1}}\in \mcO_\varpi]}\sum_{k'=0}^k\ett_{[k-N^\varrho_1=k']}\sum_{l'=0}^l\ett_{[l-M^S_1(u^\varrho_1)=l']}
\\
&\quad \cdot J(t^\eps_{i'+1},X^{t^\eps_{i'},x;u_1^\varrho,\alpha_1^S(u_1^\varrho)}_{t^\eps_{i'+1}},u^\varrho_{2,\varpi,k',l'},\alpha^S(u^\varrho))\\
&\quad+\ett_{[X^{t^\eps_{i'},x;u^\varrho_1,\alpha_1^S(u^\varrho_1)}_{t^\eps_{i'+1}}\in \mcO_0]}J(t^\eps_{i'+1},X^{t^\eps_{i'},x;u_1^\varrho,\alpha_1^S(u_1^\varrho)}_{t^\eps_{i'+1}},\emptyset,\alpha^S(u^\varrho_1))] +4\varrho
\\
&= J(t^\eps_{i'},x;u^\varrho,\alpha^S(u^\varrho))+C\sqrt{\varrho}+4\varrho
\\
&\leq \esssup_{u\in\mcU^k}J(t^\eps_{i'},x;u,\alpha^S(u))+C(\sqrt{\varrho}+\varrho),
\end{align*}
Now, as this holds for all $\alpha^S\in \mcA^{S,l,\eps}_t$ we conclude that $w^{k,l}(t^\eps_{i'},x)\leq \underline v^{k,l}_{\eps}(t^\eps_{i'},x)+C(\sqrt{\varrho}+\varrho)$, but $\varrho>0$ was arbitrary and the result follows.\\

\textbf{Step 2:} We show that $w^{k,l}(t^\eps_{i'},x)\geq \underline v^{k,l}_{\eps}(t^\eps_{i'},x)$. We again fix a $\varrho>0$ and let $(\mcO_\varpi)_{\varpi\geq 0}$ be defined as above. We pick an $x_\varpi\in\mcO_\varpi$ for each $\varpi\geq 1$ and note that by our induction assumption, there is a $\alpha_{2,\varpi,k',l'}^S\in\mcA^{S,l,\eps}_{t^\eps_{i'+1}}$ (see \cite{Buckdahn08} Lemma 4.5) such that
\begin{align*}
  w^{k',l'}_{\eps}(t^\eps_{i'+1},x_\varpi)\geq J(t^\eps_{i'+1},x_\varpi;u_2,\alpha^S_{2,\varpi,k',l'}(u_2))-\varrho,
\end{align*}
for all $u_2\in\mcU^{k,\eps}_{t^\eps_{i'+1}}$, $\varpi\geq 1$, $k'\in \{0,\ldots,k\}$ and $l'\in \{0,\ldots,l\}$. Moreover, there is an $\alpha^S_1=(\eta^S_{1,j},\theta^S_{1,j})_{j=1}^{M^S_1}\in\mcA^{S,l,\eps}_{\{t^\eps_{i'}\}}$ such that
\begin{align*}
  w^{k,l}_{\eps}(t^\eps_{i'},x)\geq G_{t^\eps_{i'},t^\eps_{i'+1}}^{t^\eps_{i'},x;u_1,\alpha^S_1(u_1)}[w^{k-N_1,l-M^S_1(u_1)} (t^\eps_{i'+1},X^{t^\eps_{i'},x;u_1,\alpha^S_1(u_1)}_{t^\eps_{i'}})]-\varrho,
\end{align*}
for all $u_1=(\tau_{1,j},\beta_{1,j})_{j=1}^{N_1}\in\mcU^{k}_{\{t^\eps_{i'}\}}$. Now, each $u=(\tau_j,\beta_j)_{j=1}^N\in\mcU^k_t$ can be uniquely decomposed as $u=u_1\oplus_{t^\eps_{i'}} u_2$ with $u_1\in\mcU^{k}_{\{t^\eps_{i'}\}}$ (with $N_1:=\max\{j\geq 0:\tau_j\leq t^\eps_{i'}\}$ interventions) and $u_2\in\mcU^{k',\eps}_{t^\eps_{i'+1}}$ (with first intervention at $\tau^2_1\geq t^\eps_{i'+1}$). Then,
\begin{align*}
  w^{k,l}_{\eps}(t^\eps_{i'},x)&\geq G_{t^\eps_{i'},t^\eps_{i'+1}}^{t^\eps_{i'},x;u_1,\alpha^S_1(u_1)}[w^{k-N_1,l-M^S_1(u_1)}(t^\eps_{i'+1},X^{t^\eps_{i'}, x;u_1,\alpha^S_1(u_1)}_{t^\eps_{i'+1}})]-\varrho
  \\
  &=G_{t^\eps_{i'},t^\eps_{i'+1}}^{t^\eps_{i'},x;u_1,\alpha^S_1(u_1)}[\sum_{k'=0}^k\ett_{[k-N_1=k']}\sum_{l'=0}^l\ett_{[l-M^S_1(u_1)=l']}  w^{k',l'}_{\eps}(t^\eps_{i'+1},X^{t^\eps_{i'},x;u_1,\alpha^S_1(u_1)}_{t^\eps_{i'+1}})]-\varrho
  \\
  &\geq G_{t^\eps_{i'},t^\eps_{i'+1}}^{t^\eps_{i'},x;u_1,\alpha^S_1(u_1)}[\sum_{k'=0}^k\ett_{[k-N_1=k']}\sum_{l'=0}^l\ett_{[l-M^S_1(u_1)=l']} \big(\sum_{\varpi\geq 1}\ett_{[X^{t^\eps_{i'},x;u_1,\alpha^S_1(u_1)}_{t^\eps_{i'+1}}\in\mcO_\varpi]}  w^{k',l'}_{\eps}(t^\eps_{i'+1},x_\varpi)
  \\
  &\quad +\ett_{[X^{t^\eps_{i'},x;u_1,\alpha^S_1(u_1)}_{t^\eps_{i'+1}}\in\mcO_0]}  w^{k',l'}(t^\eps_{i'+1},X^{t^\eps_{i'},x;u_1,\alpha^S_1(u_1)}_{t^\eps_{i'+1}})\big)]-2\varrho
  \\
  &\geq G_{t^\eps_{i'},t^\eps_{i'+1}}^{t^\eps_{i'},x;u_1,\alpha^S_1(u_1)}[\sum_{k'=0}^k\ett_{[k-N_1=k']}\sum_{l'=0}^l\ett_{[l-M^S_1(u_1)=l']}  \sum_{\varpi\geq 1}\ett_{[X^{t^\eps_{i'},x;u_1,\alpha^S_1(u_1)}_{t^\eps_{i'+1}}\in\mcO_\varpi]}
  \\
  &\quad \cdot J(t^\eps_{i'+1},x_\varpi;u_2,\alpha^S_{2,\varpi,k',l'}(u_2))]-C(\varrho+\sqrt{\varrho})
  \\
  &\geq G_{t^\eps_{i'},t^\eps_{i'+1}}^{t^\eps_{i'},x;u_1,\alpha^S_1(u_1)}[\sum_{k'=0}^k\ett_{[k-N_1=k']}\sum_{l'=0}^l\ett_{[l-M^S_1(u_1)=l']}  \sum_{\varpi\geq 1}\ett_{[X^{t^\eps_{i'},x;u_1,\alpha^S_1(u_1)}_{t^\eps_{i'+1}}\in\mcO_\varpi]}
  \\
  &\quad\cdot J(t^\eps_{i'},X^{t^\eps_{i'},x;u_1,\alpha^S(u_1)}_{t^\eps_{i'+1}};u_2,\alpha^S_{2,\varpi,k',l'}(u))]-C(\varrho+\sqrt{\varrho})
  \\
  &=J(t^\eps_{i'},x;u,\alpha_1^S(u_1)\oplus_{t^\eps_{i'}}\alpha_2^S(u))-C(\varrho+\sqrt{\varrho}),
\end{align*}
with
\begin{align*}
  \alpha_2^S(u):=\sum_{k'=0}^k\ett_{[k-N_1=k']}\sum_{l'=0}^l\ett_{[l-M^S_1(u_1)=l']}\sum_{\varpi\geq 1}\ett_{[X^{t^\eps_{i'},x;u_1,\alpha^S_1(u_1)}_{t^\eps_{i'}}\in\mcO_\varpi]}\alpha^S_{2,\varpi,k',l'}(u_2).
\end{align*}
Since $u\in\mcU^{k,\eps}_{t^\eps_{i'}}$ was arbitrary and $u\mapsto \alpha^S(u):=u\mapsto \alpha_1^S(u_1)\oplus_{t^\eps_{i'}}\alpha_2^S(u)\in\mcA^{S,l,\eps}_{t^\eps_{i'}}$, we conclude that $w^{k,l}(t^\eps_{i'},x)\geq \underline v^{k,l}_{\eps}(t^\eps_{i'},x)-C(\varrho+\sqrt{\varrho})$, where $C>0$ does not depend on $\varrho>0$ which in turn was arbitrary and the result follows.\qed\\

\noindent\emph{Proof of Proposition~\ref{prop:disc-vfs}.} Similarly to the above, we find that $\lim_{\eps\searrow 0}\bar v^{k,l}_{\eps}(t,x)=\bar v^{k,l}(t,x)$ for each $(t,x)\in [0,T]\times\R^d$ and that $\bar v^{k,l}_{\eps}$ satisfies the dynamic programming relation
\begin{align*}
  \begin{cases}
    \bar v^{k,l}_{\eps}(T,x)=\psi(x),\forall x\in\R^d,
    \\
    \bar v^{k,l}_{\eps}(t^\eps_i,x)
    =\supOP^k_\eps\infOP^l_\eps \bar v^{\cdot,\cdot}_{\eps}(t^\eps_i+,x),\quad\forall (i,x)\in \{1,\ldots,n^\eps_\bbT-1\}\times\R^d
    \\
    \bar v^{k,l}_{\eps}(t,x)
    =\E\Big[\int_{t}^{\pi^\eps_1(t)}f(s,X^{t,x,\emptyset,\emptyset}_s)ds + \bar v^{k,l}_{\eps}(\pi^\eps_1(t),X^{t,x,\emptyset,\emptyset}_{\pi^\eps_1(t)})\Big]\quad\forall (t,x)\in ([0,T)\setminus\bbT^\eps)\times\R^d.
  \end{cases}
\end{align*}
But then using an induction argument, we find that $\underline v^{k,l}_{\eps}(t,x)\leq\bar v^{k,l}_{\eps}(t,x)$ for all $(t,x)\in [0,T]\times\R^d$ and all $\eps>0$. Taking the limit as $\eps\searrow 0$ now gives that $\underline v^{k,l}(t,x)\leq\bar v^{k,l}(t,x)$.\qed\\

Combining the above proposition with \eqref{ekv:uvf-bnd}-\eqref{ekv:lvf-bnd} and Proposition~\ref{prop:trunk-vfs-approx} we arrive at the following result:

\begin{cor}\label{cor:vf-bnds}
For each $(t,x)\in [0,T]\times\R^d$ and $k,l\in\bbN$, we have $\underline v^{k,\infty}\leq \underline v(t,x)\wedge \bar v(t,x)$ and $\underline v(t,x) \vee \bar v(t,x)\leq \bar v^{\infty,l}$.
\end{cor}

\section{The dual game\label{sec:dual-game}}
A successful approach to represent the value function for various types of control problems (including those with path-dependencies) has been to consider a weak formulation where the auxiliary probability space is endowed with an independent Poisson random measure that is used to represent the control. Optimization is then carried out by altering the probability measure to modify the compensator of the random measure, so that (in the limit and on an intuitive level) the path of the corresponding Poisson jump processes have probabilistic characteristics that mimic those of an optimal control.

In the novel work presented in \cite{Kharroubi2010}, it was shown that minimal solutions to BSDEs driven by a Brownian motion and an independent Poisson random measure with constrained jumps satisfy a Feynman-Kac formula, where the corresponding partial differential equation (PDE) is a quasi-variational inequality (QVI). Moreover, when the driver is independent of the solution to the BSDE it was shown that the state process, $Y$, is the value function of an impulse control problem over randomized impulse controls. Later, \cite{KharroubiPham15} extended the approach to provide a Feynman-Kac representation of Hamilton-Jacobi-Bellman integro-partial differential equations (HJB-IPDEs).

The first result that directly linked the value function of the randomized control problem to that of the original control problem, without going through a Feynman-Kac representation, appeared in \cite{Fuhrman15} that dealt with a stochastic optimal control problem with path-dependencies over classical controls. Subsequently, this result was extended to the partial information setting in \cite{Bandini18} and optimal switching problems in \cite{Fuhrman2020}.

Our approach to show that the game introduced in the previous section has a value goes through control randomization, where we represent the impulse control of the minimizer, \ie the control denoted $\alpha$, by the sequence $(\sigma_j,\zeta_j)$ that appears in the Dirac sum formulation of the random measure, $\mu=\sum_{j\geq 1}\delta_{(\sigma_j,\zeta_j)}$, and then control the integrand in the Dol\'eans-Dade exponential appearing in a Girsanov transformation applied to the probability measure $\Prob$, effectively changing the probability distribution of $\mu$.

\subsection{Additional definitions and notation}
Let $\bigSET^\mcR$ be the set of all admissible setups
\begin{align*}
  \mathbb A=(\Omega,\mcF,\bbF^{B,\mu},\Prob,B,\mu),
\end{align*}
for the dual problem, where $(\Omega,\mcF,\bbF,\Prob)$ is a complete probability space that holds a $d$-dimensional standard Brownian motion $B$ and an independent Poisson random measure $\mu$ on $[0,T]\times A$ with compensator $\lambda$, the filtration $\bbF^{B,\mu}:=(\mcF^{B,\mu}_t)_{0\leq t\leq T}$ is the right-continuous, completed natural filtration generated by $B$ and $\mu$.

We add to the above notations by invoking the following:
\begin{itemize}
  \item We let $\mcT^{B,\mu}$ be the set of all $[0,T]$-valued $\bbF^{B,\mu}$-stopping times and for each $\eta\in\mcT^{B,\mu}$, we let $\mcT^{B,\mu}_\eta$ be the subset of stopping times $\tau$ such that $\tau\geq \eta$, $\Prob$-a.s.% Moreover, we let $\mcT^t$ (resp. $\mcT^t_\eta$) be the corresponding subsets of $\bbF^t$-stopping times, with $\tau\geq t$ (resp. $\tau\geq \eta$), $\Prob$-a.s.
  \item For $p\geq 1$, $t\in [0,T]$ and $\tau\in\mcT_t$, we let $\mcS^{p}_{[t,\tau]}(\bbA)$ be the set of all $\R$-valued, $\bbF^{B,\mu}$-progressively measurable, \cadlag processes $(Z_s: s\in [t,\tau])$ for which $\|Z\|_{\mcS^p_{[t,\tau]}(\bbA)}:=\E\big[\sup_{s\in[t,\tau]} |Z_s|^p\big]<\infty$. Moreover, we let $\mcK^p_{[t,\tau]}(\bbA)$ be the subset of $\bbF^{B,\mu}$-predictably measurable and non-decreasing processes with $Z_t=0$. Whenever $\tau=T$ we use the notations $\mcS^p_{t}(\bbA)$ and $\mcK^p_{t}(\bbA)$, respectively.
  \item We let $\mcH^{p}_{[t,\tau]}(B,\bbA)$ denote the set of all $\R^d$-valued $\mcP(\bbF^{B,\mu})$-measurable processes $(Z_s: s\in[t,\tau])$ such that $\|Z\|_{\mcH^p_{[t,\tau]}(B,\bbA)}:=\E\big[\big(\int_t^\tau |Z_s|^2 ds\big)^{p/2}\big]^{1/p}<\infty$. Furthermore, we set $\mcH^{p}_{t}(B,\bbA):=\mcH^{p}_{[t,T]}(B,\bbA)$.
  \item We let $\mcH^{p}_{[t,\tau]}(\mu,\bbA)$ denote the set of all $\R$-valued, $\mcP(\bbF^{B,\mu}) \otimes\mcB(A)$-measurable maps $(Z_s(e): (s,e)\in[t,\tau]\times A)$ such that $\|Z\|_{\mcH^p_{[t,\tau]}(\mu,\bbA)}:=\E\big[\big(\int_t^\tau \!\!\int_A |Z_s(e)|^2 \lambda(de)ds\big)^{p/2}\big]^{1/p}<\infty$ and set $\mcH^{p}_{t}(\mu,\bbA)=\mcH^{p}_{[t,T]}(\mu,\bbA)$.
  \item We let $\bigS^{p}(\bbA)$ be the set of all maps $Z:\cup_{t\in[0,T]}[t,T]\times \Omega\times \{t\}\times \R^d\to \R:(s,\omega,t,x)\mapsto Z^{t,x}_s(\omega)$ such that $Z^{t,x}\in \mcS^p_t(\bbA)$ and there is a $v\in\Pi^g_c$ such that $v(t,x)=Z^{t,x}_t$, $\Prob$-a.s., for all $(t,x)\in [0,T]\times\R^n$. Moreover, given $Z\in\bigS^p(\bbA)$, we let $(t,x)\mapsto Z(t,x)$ denote this deterministic function, so that $ Z(t,x)=v(t,x)$ for all $(t,x)\in[0,T]\times \R^d$.
  \item We define $\mcU^{\mcR,k}_t(\bbA)$ (resp.~$\mcU^{\mcR,k,\eps}_t(\bbA)$) analogously to $\mcU^{k}_t$ (resp.~$\mcU^{k,\eps}_t$) using the filtration $\bbF^{B,\mu}$ instead of $\bbF$. To define $\mcU^{\mcR}_t(\bbA)$ we also replace the square integrability constraint on the total cost of interventions and assume that, instead, $\Lambda^{t,x;u^\mcR}_T\in L^2(\Omega,\mcF,\Prob)$, where\footnote{For the definition of $X^{t,x;u^\mcR}$ see \eqref{ekv:fwd-sde-dual}.}
    \begin{align*}
      \Lambda^{t,x;u^\mcR}_s:=\sum_{j=1}^\infty\ett_{[\tau_j < s]}\ell(\tau_j,X^{t,x,[u^\mcR]_{j-1}}_{\tau_j},\beta_j),\quad \forall s\in [0,T],
    \end{align*}
    whenever $u^\mcR\in\mcU^{\mcR}_t(\bbA)$.
\end{itemize}
In the dual game, the minimizer will not directly choose the intervention times and the corresponding interventions. Instead, the minimizer influences the probability distribution for the jumps by altering the compensation of the Poisson random measure. We thus let $\mcV(\bbA)$ be the set of all $\Pred(\bbF^{B,\mu})\otimes\mcB(A)$-measurable, bounded maps $\nu=\nu_t(\omega,e):[0,T]\times\Omega\times U\to [0,\infty)$. For each $(t,x)\in[0,T]$, $u^\mcR=(\tau^\mcR_j,\beta^\mcR_j)_{j=1}^{N^\mcR}\in\mcU^{\mcR}_t(\bbA)$ and $\nu\in\mcV(\bbA)$ we then define the corresponding reward/cost functional
\begin{align}\nonumber
   J^{\mcR}(t,x;u^\mcR,\nu,\bbA)&:=\E^\nu\Big[\psi(X_T^{t,x;u^\mcR})+\int_t^{T}f(s,X_s^{t,x;u^\mcR})ds - \sum_{j=1}^{N^\mcR}\ell(\tau^\mcR_j,X^{t,x;[u^\mcR]_{j-1}}_{\tau^\mcR_j},\beta^\mcR_j)
   \\
   &\quad+ \int_t^T\int_A\chi(s,X_{s-}^{t,x;u^\mcR},e)\mu(ds,de)\Big|\mcF^{B,\mu}_t\Big],\label{ekv:obj-fun-dual}
\end{align}
where $X^{t,x;u^\mcR}$ is the unique solution to the SDE
\begin{align}\nonumber
X^{t,x;u^\mcR}_s&=x+\int_t^s a(r,X^{t,x;u^\mcR}_r)dr+\int_t^s\sigma(r,X^{t,x;u^\mcR}_r)dB_r+\sum_{j=1}^N\ett_{[\tau_j\leq s]}\jmp(\tau^\mcR_j,X^{t,x;[u^\mcR]_{j-1}}_{\tau^\mcR_j},\beta^\mcR_j)
\\
&\quad +\int_t^s\gamma(r,X^{t,x;u^\mcR}_{r-},e)d\mu(dr,de),\quad\forall s\in [t,T]\label{ekv:fwd-sde-dual}
\end{align}
and $\E^\nu$ is expectation with respect to the probability measure $\Prob^\nu$ on $(\Omega,\mcF)$ defined by $d\Prob^\nu:=\kappa^\nu d\Prob$ with
\begin{align*}
\kappa^{\nu}_s:=\exp\Big(\int_{0}^s\int_A(1-\nu_r(e))\lambda(de)dr\Big)\prod_{\sigma_j\leq s}\nu_{\sigma_j}(\zeta_j).
\end{align*}
We define the set of non-anticipative density strategies in the following natural way:
\begin{defn}
We let $\mcV^S(\bbA)$ be the set of maps $\nu^S:\mcU^{\mcR}(\bbA)\to \mcV(\bbA)$ such that for any $u^\mcR,\tilde u^\mcR\in\mcU^{\mcR}(\bbA)$ it holds that $\nu^S[u^\mcR]=\nu^S[\tilde u^\mcR]$, $d\Prob\otimes ds\otimes d\lambda$-a.e.~on $\Omega\times[0,t]\times A$ whenever $u^\mcR_t=\tilde u^\mcR_t$, $\Prob$-a.s.% Moreover, we let $\mcV^S_{\inf>0}(\bbA)$ be the subset of $\mcV^S(\bbA)$ containing all maps $\nu^S:\mcU^{\mcR}(\bbA)\to \mcV_{\inf>0}(\bbA)$.
\end{defn}
We then consider the value function
\begin{align*}
  v^{\mcR}(t,x):=\essinf_{\nu^S\in\mcV^{S}(\bbA)}\esssup_{u^\mcR\in\mcU^{\mcR}_t(\bbA)}J^{\mcR}(t,x;u^\mcR,\nu^S(u^\mcR),\bbA).
\end{align*}

\subsection{A truncation}
As in Section~\ref{sec:game-trunc}, we apply a truncation to the control sets. In the randomized setting it turns out to be more effective to put an upper bound on the map $\nu$, rather that to directly limit the number of interventions. For each $n\in\bbN$, we thus introduce the set\footnote{To simplify notation we will suppress the dependence of $\bbA$ through the remined of this section.} $\mcV^n$ as the subset of $\mcV$ containing all maps $\nu=\nu_t(\omega,e):[0,T]\times\Omega\times A\to [0,n]$ and let $\mcV^n_{\inf>0}:=\{\nu\in\mcV^n:\inf\nu>0\}$. In addition, we let $\mcV_{\inf>0}^{S,n}$ be the subset of $\mcV^S_{\inf>0}$ with all maps $\nu^S:\mcU^{\mcR}\to\mcV^{n}_{\inf>0}$.

For $(t,x)\in[0,T]\times\R^d$, we introduce the double sequence $(Y^{t,x,k,n})_{k,n\in\bbN}$, where for each $n\in\bbN$, the process $Y^{t,x,0,n}\in\mcS^2_t$ is the first component in the unique solution to the standard BSDE with jumps (recall the definition of $X^{t,x}$ in \eqref{ekv:fwd-sde-Markov})
\begin{align}\nonumber
    Y^{t,x,0,n}_s&=\psi(X^{t,x}_T)+\int_s^T f(r,X^{t,x}_r)dr-n\int_t^T\!\!\int_A(V^{t,x,0,n}_r(e)+\chi(r,X^{t,x}_{r},e))^-\lambda(de)dr
    \\
    &\quad -\int_s^T Z^{t,x,0,n}_r dB_r-\int_{s}^T\!\!\!\int_A V^{t,x,0,n}_{r}(e)\mu(dr,de),\quad \forall s\in [t,T],\label{ekv:rbsde-pen-0}
\end{align}
whereas $Y^{t,x,k,n}$ for $k\geq 1$ and $n\in\bbN$, is defined recursively by letting $(Y^{t,x,k,n},Z^{t,x,k,n},V^{t,x,k,n},K^{+,t,x,k,n})\in\mcS^2_t\times\mcH^2_t(B)\times\mcH^2_t(\mu)\times\mcA^2_t$ be the unique solution to the reflected BSDE with jumps
\begin{align}\label{ekv:rbsde-pen-trunk}
  \begin{cases}
    Y^{t,x,k,n}_s=\psi(X^{t,x}_T)+\int_s^T f(r,X^{t,x}_r)dr-n\int_t^T\!\!\int_A(V^{t,x,k,n}_r(e)+\chi(r,X^{t,x}_{r},e))^-\lambda(de)dr
    \\
    \quad -\int_s^T Z^{t,x,k,n}_r dB_r-\int_{s}^T\!\!\!\int_A V^{t,x,k,n}_{r}(e)\mu(dr,de)+K^{+,t,x,k,n}_T-K^{+,t,x,k,n}_s,\quad \forall s\in [t,T],
    \\
    Y^{t,x,k,n}_s\geq \supOP Y^{k-1,n}(s,X^{t,x}_s),\, \forall s\in [t,T] \quad\text{and}\quad\int_t^T\!\! \big(Y^{t,x,k,n}_s-\supOP Y^{k-1,n}(s,X^{t,x}_s)\big)dK^{+,t,x,k,n}_s
  \end{cases}
\end{align}
where $Y^{k-1,n}(t,x):=Y^{t,x,k-1,n}_t$. Moreover, we set $K^{-,t,x,k,n}_t:=n\int_0^t\!\!\int_A(V^{t,x,k,n}_s(e)+\chi(s,X^{t,x}_{s},e))^-\lambda(de)ds$. Existence and uniqueness of solutions to this system of BSDEs follow from repeated use of Proposition~\ref{prop:rbsde-jmp}.

%The following proposition shows that, when truncating the densities, the dual game has a value. The proof is rather standard and is therefore relegated to Appendix~\ref{app:dual-trunc}.

The following proposition shows that, when truncating the densities, the dual game has a value. Since the proof is fairly standard, it is deferred to Appendix~\ref{app:dual-trunc}.

\begin{prop}\label{prop:trunk-game}
For each $k,n\in\bbN$, there is a $v^{\mcR,k,n}\in\Pi^g_c$ (that does not depend on $\bbA\in\bigSET^\mcR$) such that for each $(t,x)\in[0,T]\times\R^d$, we have
\begin{align}\label{ekv:trunk-game-value}
  v^{\mcR,k,n}(t,x)=Y^{t,x,k,n}_{t}=\essinf_{\nu\in\mcV^{S,n}}\esssup_{u^\mcR\in\mcU^{\mcR,k}_{t}}J^{\mcR}(t,x;u^\mcR,\nu^S(u^\mcR)) = \esssup_{u^\mcR\in\mcU^{\mcR,k}_{t}}\essinf_{\nu\in \mcV^{n}}J^{\mcR}(t,x;u^\mcR,\nu).
\end{align}
\end{prop}

\begin{rem}\label{rem:inf-g-0}
For later use it will be convenient to note that we have the equivalent representation
\begin{align*}
  v^{\mcR,k,n}(t,x)=\essinf_{\nu\in\mcV_{\inf>0}^{S,n}}\esssup_{u^\mcR\in\mcU^{\mcR,k}_{t}}J^{\mcR}(t,x;u^\mcR,\nu^S(u^\mcR)) = \esssup_{u^\mcR\in\mcU^{\mcR,k}_{t}}\essinf_{\nu\in \mcV_{\inf>0}^{n}}J^{\mcR}(t,x;u^\mcR,\nu).
\end{align*}
\end{rem}

\begin{cor}\label{cor:Ykn-ineq}
$Y^{t,x,k,n+1}\leq Y^{t,x,k,n}\leq Y^{t,x,k+1,n}$
\end{cor}

\noindent\emph{Proof.} The inequalities are immediate from the representation of $Y^{t,x,k,n}$ in \eqref{ekv:trunk-game-value} as $\mcV^{t,x,n}\subset\mcV^{t,x,n+1}$ and $\mcU^{\mcR,k}_t\subset\mcU^{\mcR,k+1}_t$ for all $(t,x,k,n)\in[0,T]\times\R^d\times\bbN\times\bbN$.\qed\\

From \cite{dbl-obst-qvi} we know that there are maps $\underline v^{\mcR,k},\bar v^{\mcR,n}, v\in\Pi^g_c$ (that are all independent of the choice of probability space $\bbA$) such that for each $(t,x,k,n)\in[0,T]\times\R^d\times\bbN\times\bbN$,
\begin{align}\label{ekv:dual-bnds}
  \underline v^{\mcR,k}(t,x) \leq Y^{t,x,k,n}_t \leq \bar v^{\mcR,n}(t,x)
\end{align}
and
\begin{align}\label{ekv:dbl-obst}
  \lim_{k\to\infty} \underline v^{\mcR,k}(t,x) =  \lim_{n\to\infty} \bar v^{\mcR,n}(t,x)=v(t,x).
\end{align}
In \cite{dbl-obst-qvi}, the map $v$ is further characterized as the unique viscosity solution in $\Pi^g$ to \eqref{ekv:var-ineq}. Since any non-anticipative map in $\mcV^{S}$ is uniformly bounded, taking the limit as $n\to\infty$ in \eqref{ekv:trunk-game-value} and using \eqref{ekv:dual-bnds} gives
\begin{align}\label{ekv:dual-lower-bnd}
  \underline v^{\mcR,k}(t,x)\leq \essinf_{\nu\in\mcV^{S}}\esssup_{u^\mcR\in\mcU^{\mcR,k}_{t}}J^{\mcR}(t,x;u^\mcR,\nu^S(u^\mcR))
\end{align}
We will not make use of the precise form of the bounds $\underline v^{\mcR,k}$ and $\bar v^{\mcR,n}$ in the present work. Instead, value of the original game along with the characterization in terms of viscosity solutions is proved by combining \eqref{ekv:dbl-obst} with the following proposition, the proof of which relies on \eqref{ekv:trunk-game-value}, \eqref{ekv:dual-bnds} and \eqref{ekv:dual-lower-bnd}.
\begin{prop}\label{prop:rel-rand-to-orig}
We have
\begin{enumerate}[(A)]
  \item for each $(t,x)\in[0,T]\times\R^d$ and $n\in\bbN$ we have $\lim_{l\to\infty}\bar v^{\infty,l}(t,x)\leq \bar v^{\mcR,n}(t,x)$; and
  \item for each $(t,x)\in[0,T]\times\R^d$ and $k\in\bbN$ we have $\underline v^{\mcR,k}(t,x)\leq \underline v^{k,\infty}(t,x)$.
\end{enumerate}
\end{prop}
The remainder of the paper is devoted to proving Proposition~\ref{prop:rel-rand-to-orig}.

\begin{rem}
Note that Proposition~\ref{prop:rel-rand-to-orig} together with Corollary~\ref{cor:vf-bnds} implies that for each $k,n\in\bbN$, we have
\begin{align*}
  \underline v^{\mcR,k}(t,x)\leq \underline v(t,x) , \bar v(t,x) \leq \bar v^{\mcR,n}(t,x)
\end{align*}
from which value of the game follows by \eqref{ekv:dbl-obst}.
\end{rem}

\subsection{Discretization of the impulse control for the randomized problem}
%Arguing as above, we note that there is a map $\underline v^{\mcR,k}_{\eps}:[0,T]\times\R^d\to\R$ such that
%\begin{align*}
%  \underline v^{\mcR,k}_{\eps}(t,x)=\essinf_{\nu^S\in\mcV^S}\esssup_{u\in\mcU^{B,\pi^\eps_A(\mu),k,\eps}_t}J^\mcR(t,x;u,\nu^S(u))
%\end{align*}
%for each $(t,x)\in [0,T]\times\R^d$. The aim of this part is to find an approximately optimal control for the randomized problem that lies in the set $\mcU^{B,\pi^\eps_A(\mu),k,\eps}$ for some $\eps>0$.

To compare the original and dual formulations, we make use of the discretized version of the original game. Consequently, we also require a corresponding discretization of the dual game. To construct such a discretization in the randomized setting, we define the discretized random measure $\mu^{t,\eps}$ as $\mu^{t,\eps}=(\sigma^{t,\eps}_j,\zeta^{t,\eps}_j)_{j=1}^{M^{\mu^t}_{T}}:=\pi^\eps((\sigma_j,\zeta_j)_{j=M^\mu_{t}}^{M^\mu_{T}}):= \pi^\eps(\mu^t)$ and set $R^{t,x;u,\eps}:=\lim_{j\to\infty}R^{t,x;u,\eps,j}$, the pointwise limit of the solutions to
\begin{align}\nonumber
R^{t,x;u,\eps,j}_s&=x+\int_t^s a(r,R^{t,x;u^\mcR,\eps,j}_r)dr+\int_t^s\sigma(r,R^{t,x;u^\mcR,\eps,j}_r)dB_r+\sum_{i=1}^{N^\mcR}\ett_{[\tau^\mcR_i\leq s]}\jmp(\tau^\mcR_i,R^{t,x;[u^\mcR]_{i-1},\eps,j}_{\tau^\mcR_i},\beta^\mcR_i)
\\
&\quad +\sum_{i=1}^{M^{\mu^t}_{T}\wedge j}\ett_{[\sigma^{t,\eps}_i\leq s]}\gamma(\sigma^{t,\eps}_i,R^{t,x;u^\mcR,\eps,i-1}_{\sigma^{t,\eps}_i},\zeta^{t,\eps}_i),\quad\forall s\in [t,T],\label{ekv:fwd-sde-dual-eps}
\end{align}
so that $R^{t,x;u^\mcR,\eps}$ corresponds to delaying the Poisson induced jumps in $X^{t,x;u^\mcR}$ to the next point in $\bbT^\eps$. We then let
\begin{align*}
  X^{t,x;u^\mcR,\eps}_s&=R^{t,x;u^\mcR,\eps}_{t^\eps_i}+\int_{t^\eps_i}^s a(r,X^{t,x;u^\mcR,\eps}_r)dr+\int_{t^\eps_i}^s\sigma(r,X^{t,x;u^\mcR,\eps}_r)dB_r
\\
&\quad + \int_{t^\eps_i}^s\int_A\gamma(r,X^{t,x;u^\mcR,\eps}_{r-},e)\mu(dr,de),\quad\forall s\in [{t^\eps_i},{t^\eps_{i+1}})
\end{align*}
and define
\begin{align*}
  J^{\mcR,\eps}(t,x;u^\mcR,\nu)&:=\E^\nu\Big[\psi(R^{t,x;u^\mcR,\eps}_T)+\int_t^{T}f(s,X^{t,x;u^\mcR,\eps}_s)ds - \sum_{j=1}^{N^\mcR}\ell(\tau^\mcR_j,R^{t,x;[u^\mcR]_{j-1},\eps}_{\tau^\mcR_j},\beta^\mcR_j)
  \\
  &\quad + \int_{t}^T\int_A\chi(s,X^{t,x;u^\mcR,\eps}_{s},e)\mu(ds,de)\Big|\mcF^{\mcR}_{t}\Big].
\end{align*}
We then let $\underline v^{\mcR,k}_{\eps}:[0,T]\times\R^d\to\R$ be the unique solution to the to the dynamic programming relation
\begin{align}\label{ekv:dynP-eps-dual}
  \begin{cases}
    \underline v^{\mcR,k}_{\eps}(T,x)=\psi(x),\quad\forall x\in\R^d,
    \\
    \underline v^{\mcR,k}_{\eps}(t^\eps_i,x)
    =\underline v^{\mcR,k}_{\eps}(t^\eps_i+,x)\vee\supOP_\eps \underline v^{\mcR,k-1}_{\eps}(t^\eps_i,x),\quad\forall (i,x)\in \{1,\ldots,n^\eps_\bbT-1\}\times\R^d
    \\
     \underline v^{\mcR,k}_{\eps}(t,x)
    =\essinf_{\nu\in\mcV} \E^\nu\Big[\int_{t}^{\pi^\eps_1(t)}f(s,X^{t,x;\emptyset,\eps}_s)ds +\int_{t}^{\pi^\eps_1(t)}\int_A\chi(s,X^{t,x;\emptyset,\eps}_{s},e)\mu(ds,de)\\
    \quad + \underline v^{\mcR,k}_{\eps}(\pi^\eps_1(t),R^{t,x;\emptyset,\eps}_{\pi^\eps_1(t)})\Big|\mcF^{B,\mu}_t\Big],\quad\forall (t,x)\in ([0,T)\setminus\bbT^\eps)\times\R^d,
  \end{cases}
\end{align}
where continuity of the map $x\mapsto \underline v^{\mcR,k}_{\eps}(t^\eps_i,x)$ follows by standard arguments for BSDEs with constrained jumps. %Continuity of $x\mapsto \underline v^{\mcR,k}_{\eps}(t^\eps_i,x)$ also gives that for each $t\in [0,T]$ we have $\essinf_{\nu\in\mcV} \E^\nu\Big[\int_{t}^{\pi^\eps_1(t)}f(s,X^{t,x;\emptyset,\eps}_s)ds +\int_{t}^{\pi^\eps_1(t)}\int_A\chi(s,X^{t,x;u,\eps}_{s},e)\mu(ds,de) + \underline v^{\mcR,k}_{\eps}(\pi^\eps_1(t),R^{t,x;\emptyset,\eps}_{\pi^\eps_1(t)})\Big|\mcF^{B,\mu}_t\Big]=w(x)$ for some continuous function $w:\R^d\to\R$.

We let $\bbF^{B,\mu^{t,\eps}}:=(\mcF^{B,\mu^{t,\eps}}_{s})_{0\leq s\leq T}$ be the augmented natural filtration generated by $B$ and $\mu^{t,\eps}$ and let $\mcU^{B,\mu^t,k,\eps}_t$ be the subset of $u^\mcR=(\tau^\mcR_j,\beta^\mcR_j)_{j=1}^{N^\mcR}\in \mcU^{\mcR,k,\eps}_t$ such that $\tau^\mcR_j$ is a $\bbF^{B,\mu^{t,\eps}}$-stopping time and $\beta_j^\mcR$ is $\mcF^{B,\mu^{t,\eps}}_{\tau^\mcR_j}$-measurable. Then, $R^{t,x;u^\mcR,\eps}$ is $\bbF^{B,\mu^{t,\eps}}$-adapted whenever $u^\mcR\in\mcU^{B,\mu^{t,\eps},k,\eps}_t$. Moreover, we let

\begin{lem}
Let $u^{\mcR,k,\eps}:=(\tau^{\mcR,k,\eps}_j,\beta^{\mcR,k,\eps}_j)_{j=1}^{N^{\mcR,k,\eps}}\in\mcU^{B,\mu^{t,\eps},k,\eps}_t$ be defined as
\begin{itemize}
  \item $\tau^{\mcR,k,\eps}_j:=\inf\{s\in\bbT^\eps_{\tau^{\mcR,k,\eps}_{j-1}}:\underline v^{\mcR,k-j+1}_\eps(s,R^{t,x;[u^{\mcR,k,\eps}]_{j-1},\eps}_s)=\supOP \underline v^{\mcR,k-j}_\eps(s,R^{t,x;[u^{\mcR,k,\eps}]_{j-1},\eps}_s)\}\wedge T$
  \item $\beta^{\mcR,k,\eps}_j\in\argmax_{b\in \bar U^\eps}\{\underline v^{\mcR,k-j}_\eps(\tau^{\mcR,k,\eps}_j,R^{t,x;[u^{\mcR,k,\eps}]_{j-1},\eps}_{\tau^{\mcR,k,\eps}_j}) + \jmp(\tau^{\mcR,k,\eps}_j,R^{t,x;[u^{\mcR,k,\eps}]_{j-1},\eps}_{\tau^{\mcR,k,\eps}_j},b))$\\$ - \ell(\tau^{\mcR,k,\eps}_j,R^{t,x;[u^{\mcR,k,\eps}]_{j-1},\eps}_{\tau^{\mcR,k,\eps}_j},b)\}$
\end{itemize}
with $\tau^{\mcR,k,\eps}_0:=t$ and $N^{\mcR,k,\eps}:=\max\{j:\tau^{\mcR,k,\eps}_j<T\}$. Then,
\begin{align}\label{ekv:dual-eps-rep}
  \underline v^{\mcR,k}_{\eps}(t,x) = \essinf_{\nu^S\in\mcV^S}\esssup_{u^\mcR\in\mcU^{\mcR,k,\eps}_t}J^{\mcR,\eps}(t,x;u^\mcR,\nu^S(u^\mcR)) = \essinf_{\nu\in\mcV}J^{\mcR,\eps}(t,x;u^{\mcR,k,\eps},\nu).
\end{align}
\end{lem}

\noindent\emph{Proof.} \textbf{Step 1.} We first prove that
\begin{align}\label{ekv:dual-eps-1}
  \underline v^{\mcR,k}_{\eps}(t,x)\leq \essinf_{\nu\in\mcV}J^{\mcR,\eps}(t,x;u^{\mcR,k,\eps},\nu).
\end{align}
Arguing as in Proposition 4.1 of \cite{dbl-obst-qvi} we find that there is a quadruple $(Y,Z,V,K)$, with
\begin{align*}
  Y_s:=\ett_{[t,\tau^{\mcR,k,\eps}_1]}\underline v^{\mcR,k}_{\eps}(s,X^{t,x;\emptyset,\eps}_s)+\sum_{j=1}^k\ett_{(\tau^{\mcR,k,\eps}_j,\tau^{\mcR,k,\eps}_{j+1}]}(s)\underline v^{\mcR,k-j}_{\eps}(s,X^{t,x;[u^{w}]_j,\eps}_s)
\end{align*}
and $(Z,V,K)\in\mcH^2(B)\times\mcH^2(\mu)\times\mcK^2$ such that
\begin{align}\label{ekv:bsde-c-jmp-eps}
  \begin{cases}
  Y_{s}=\psi(R^{t,x;u^{\mcR,k,\eps}}_T)+\int_{s}^{T} f\big(r,X^{t,x;u^{\mcR,k,\eps},\eps}_r\big)dr-\int_{s}^{T}Z_rdB_r - \int_{{s'}}^{T}\!\!\!\int_A V_{r}(e)\mu(dr,de)
  \\
  \quad -\sum_{j=1}^{N^{\mcR,k,\eps}}\ett_{[{s}\leq \tau^{\mcR,k,\eps}_j]}\ell({\tau_j^{\mcR,k,\eps}},R^{t,x;[u^{\mcR,k,\eps}]_{j-1},\eps}_{\tau_j^{\mcR,k,\eps}},\beta_j^{\mcR,k,\eps})-(K_T-K_{s}), \quad\forall {s}\in[t,T],
  \\
  V_s\geq -\chi(s,X^{t,x;u^{\mcR,k,\eps},\eps}_{s-},e),\quad d\Prob\otimes ds\otimes \lambda(de)-a.e.
  \end{cases}
\end{align}
On the other hand, by repeating the argument in \cite{Kharroubi2010}, we obtain that \eqref{ekv:bsde-c-jmp-eps} admits a unique maximal solution $(\tilde Y,\tilde Z,\tilde V,\tilde K)$, which has the following representation:
\begin{align*}
  \tilde Y_t=\essinf_{\nu\in\mcV}J^{\mcR,\eps}(t,x;u^{\mcR,k,\eps},\nu).
\end{align*}
Since $\tilde Y$ is the maximal solution, it follows that $Y_t\leq \tilde Y_t$, and thus \eqref{ekv:dual-eps-1} holds.\\
%In particular, $\underline v^{\mcR,k}_{\eps}(t,x)=Y_t$. On the other hand, for any $u\in\mcU_t^{B,\mu,k,\eps}$, a straightforward application of the theory in \cite{Kharroubi2010} gives that $\underline v^{\mcR,k}_\eps(t,x)$ dominates $\tilde Y_t$, where $(\tilde Y,\tilde Z,\tilde V,\tilde K)$ is the maximal solution to
%\begin{align*}
%  \begin{cases}
%  \tilde Y_{s}=\psi(R^{t,x;u}_T)+\int_{s}^{T} f\big(r,X^{t,x;u,\eps}_r\big)dr-\int_{s}^{T}\tilde Z_rdB_r - \int_{{s'}}^{T}\!\!\!\int_A \tilde V_{r}(e)\mu(dr,de)
%  \\
%  \quad -\sum_{j=1}^{N}\ett_{[{s}\leq \tau_j]}\ell({\tau_j},R^{t,x;[u]_{j-1},\eps}_{\tau_j},\beta_j)-(\tilde K_T-\tilde K_{s}),\quad\forall {s}\in[t,T],
%  \\
%  \tilde V_s\geq -\chi(s,X^{t,x;u,\eps}_{s-},e),\quad d\Prob\otimes ds\otimes \lambda(de)-a.e.
%  \end{cases}
%\end{align*}
%Consequently, \eqref{ekv:dual-eps-1} holds.\\

\textbf{Step 2.} We move on to prove that
\begin{align}\label{ekv:dual-eps-2}
  \underline v^{\mcR,k}_{\eps}(t,x)\geq \essinf_{\nu^S\in\mcV^S}\esssup_{u^\mcR\in\mcU^{\mcR,k,\eps}_t}J^{\mcR,\eps}(t,x;u^\mcR,\nu^S(u)).
\end{align}
The proof is based on backward induction and we assume that for some $i\in \{2,\ldots,n_\bbT^\eps\}$, Equation \eqref{ekv:dual-eps-2} holds for $t=t^\eps_i$ and all $x\in\R^d$ and consider $t=t^\eps_{i-1}$ (the argument for $t\in (t^\eps_{i-1},t^\eps_i)$ is simpler and will be left out). Step 1 and continuity gives that
\begin{align*}
  \underline v^{\mcR,k}_{\eps}({t^\eps_{i-1}},x)& = \max_{(b_j)_{j=1}^\kappa\in \cup_{j=1}^k (\bar U^\eps)^j}\essinf_{\nu\in\mcV}\E^\nu\Big[\int_{t^\eps_{i-1}}^{t^\eps_i}f(s,X^{{t^\eps_{i-1}},x;(t^\eps_{i-1},b_j)_{j=1}^\kappa,\eps}_s)ds
  \\
  &\quad+\int_{t^\eps_{i-1}}^{t^\eps_i}\int_A\chi(s,X^{{t^\eps_{i-1}},x;(t^\eps_{i-1},b_j)_{j=1}^\kappa,\eps}_{s},e)\mu(ds,de)
  \\
    &\quad + \underline v^{\mcR,k-\kappa}_{\eps}(t^\eps_i,R^{{t^\eps_{i-1}},x;(t^\eps_{i-1},b_j)_{j=1}^\kappa,\eps}_{t^\eps_i}) + \sum_{j=1}^N\ell({t^\eps_{i-1}},R^{t,x;(t^\eps_{i-1},b_{j'})_{j'=1}^j,\eps},b_j)\Big|\mcF^{B,\mu}_{t^\eps_{i-1}}\Big]
    \\
    & = \essinf_{\nu^S\in\mcV^S}\esssup_{u^\mcR\in\mcU^{\mcR,k,\eps}_{\{t^\eps_{i-1}\}}}\E^{\nu^S(u^\mcR)}\Big[\int_{t^\eps_{i-1}}^{t^\eps_i} f(s,X^{{t^\eps_{i-1}},x;u^\mcR,\eps}_s)ds + \int_{t^\eps_{i-1}}^{t^\eps_i}\int_A\chi(s,X^{{t^\eps_{i-1}},x;u^\mcR,\eps}_{s},e)\mu(ds,de)
  \\
    &\quad + \underline v^{\mcR,k-N^\mcR}_{\eps}(t^\eps_i,R^{{t^\eps_{i-1}},x;u^\mcR,\eps}_{t^\eps_i}) + \sum_{j=1}^{N^\mcR}\ell(\tau^\mcR_j,R^{t,x;[u^\mcR,\eps}]_{j-1},\beta^\mcR_j)\Big|\mcF^{B,\mu}_{t^\eps_{i-1}}\Big]
\end{align*}
We pick $\varrho>0$, and note that there is a $\nu^{S,\varrho}_1\in\mcV^S$ such that for any $u=(\tau_j,\beta_j)_{j=1}^N$, we have
\begin{align*}
  \underline v^{\mcR,k}_{\eps}({t^\eps_{i-1}},x)&\geq \E^{\nu^{S,\varrho}_1(u^1)}\Big[\int_{t^\eps_{i-1}}^{t^\eps_i}f(s,X^{{t^\eps_{i-1}},x; u^1,\eps}_s)ds+\int_{t^\eps_{i-1}}^{t^\eps_i}\int_A\chi(s,X^{{t^\eps_{i-1}},x;u^1,\eps}_{s},e)\mu(ds,de)
  \\
    &\quad + \underline v^{\mcR,k-N_1}_{\eps}(t^\eps_i,R^{{t^\eps_{i-1}},x;(t^\eps_{i-1},b_j)_{j=1}^{N_1},\eps}_{t^\eps_i}) + \sum_{j=1}^{N_1}\ell({t^\eps_{i-1}},R^{t,x;[u^1]_{j-1},\eps},\beta_j)\Big]-\varrho,
\end{align*}
where $u^1=(\tau_j,\beta_j)_{j=1}^{N_1}$ with $N_1:=\max\{j\leq 0:\tau_j= t^\eps_{i-1}\}$. Now, by continuity of $x\mapsto \underline v^{\mcR,j}_{\eps}(t^\eps_i,x)$ and since $R^{{t^\eps_{i-1}},x;u^1,\eps}_{t^\eps_i}$ has moments of all orders under $\Prob^{\nu^{S,\varrho}_1(u^1)}$ there is a finite collection of sets $\{\mcO_m\}_{m=1}^{K}$ and points $x_{m}\in \mcO_m\subset\R^d$ such that $|v^{\mcR,k-\kappa}_{\eps}(t^\eps_i,x)-v^{\mcR,k-\kappa}_{\eps}(t^\eps_i,x_m)|\leq \varrho$ for all $x\in \mcO_m$ for $m=1,\ldots,K$ and $\kappa=1,\ldots,k$ and, moreover, $\Prob^{\nu^{S,\varrho}_1(u^1)}[R^{{t^\eps_{i-1}},x; u^1,\eps}_{t^\eps_i}\in(\mcO^c_1\cap\cdots\cap\mcO^c_K)]\leq\varrho$. Hence,
\begin{align*}
  \tilde v^{\mcR,k}_{\eps}({t^\eps_{i-1}},x)&\geq \E^{\nu^{S,\varrho}_1(u^1)}\Big[\int_{t^\eps_{i-1}}^{t^\eps_i}f(s,X^{{t^\eps_{i-1}},x;u^1,\eps}_s)ds +\int_{t^\eps_{i-1}}^{t^\eps_i}\int_A\chi(s,X^{{t^\eps_{i-1}},x;u^1,\eps}_{s},e)\mu(ds,de)
  \\
    &\quad + \sum_{m=1}^K\ett_{[R^{{t^\eps_{i-1}},x; u^1,\eps}_{t^\eps_i}\in\mcO_m]}\sum_{k'=1}^k\ett_{[k-N_1=k']}\underline v^{\mcR,k'}_{\eps}(t^\eps_i,x_m) + \sum_{j=1}^{N_1}\ell({t^\eps_{i-1}},R^{t,x;[u^1]_{j-1},\eps},b_j)\Big]-2\varrho-C\sqrt{\varrho}.
\end{align*}
On the other hand, for each $m\in\{1,\ldots,K\}$ and $k'\in\{1,\ldots,k\}$, there is a $\nu^{S,\varrho}_{2,m,k'}\in\mcV^S$ such that
\begin{align*}
  \underline v^{\mcR,k'}_{\eps}(t^\eps_i,x_m)&\geq\sup_{\tilde u\in\mcU^{\mcR,k',\eps}_{t^\eps_i}}J^{\mcR,\eps}(t^\eps_i,x_m;\tilde u,\nu^{S,\varrho}_{2,m,k'}(u^1\oplus_{t^\eps_{i}} \tilde u))-\varrho
  \\
  &\geq J^{\mcR,\eps}(t^\eps_i,x_m;u^2,\nu^{S,\varrho}_{2,m,k'}(u^1\oplus_{t^\eps_{i}} u^2))-\varrho
\end{align*}
with $u^2=(\tau_j,\beta_j)_{j=N_1+1}^{N}$. Now,
\begin{align*}
  u\mapsto \{(r,e)\mapsto\nu^{S,\varrho}[u](r,e)\}:=(r,e)\mapsto\ett_{[0\leq r\leq t^\eps_{i}]}\nu^{S,\varrho}_1[u_1](r,e)
  \\
  +\ett_{[ t^\eps_{i}<r\leq T]}\sum_{m=1}^K\ett_{[R^{{t^\eps_{i-1}},x; u^1,\eps}_{t^\eps_i}\in\mcO_m]}\sum_{k'=1}^k\ett_{[k-N_1=k']}\nu^{S,\varrho}_{2,m,k'}[u_2](r,e)&\in\mcV^S
\end{align*}
and using the tower property we conclude that $\underline v^{\mcR,k}_{\eps}({t^\eps_{i-1}},x)\geq J^{\mcR,\eps}(t^\eps_{i-1},x;u,\nu^{S,\varrho}(u))-C(\varrho+\sqrt{\varrho})$ where $C>0$ does not depend on $\varrho$, $u$ or $\nu^{S,\varrho}$. As $u$ was arbitrary we may take the supremum over all $u\in\mcU^{\mcR,k,\eps}_{t^\eps_{i-1}}$ and \eqref{ekv:dual-eps-2} follows.\\

Combining the two above steps and using the trivial inequality
\begin{align*}
  \essinf_{\nu\in\mcV}J^{\mcR,\eps}(t,x;u^{\mcR,k,\eps},\nu)\leq \essinf_{\nu^S\in\mcV^S}\esssup_{u^\mcR\in\mcU^{\mcR,k,\eps}_t}J^{\mcR,\eps}(t,x;u^\mcR,\nu^S(u^\mcR))
\end{align*}
proves \eqref{ekv:dual-eps-rep}.\qed\\

\begin{lem}
For each $(t,x,k)\in[0,T]\times\R^d\times \bbN$, we have
\begin{align*}
  \lim_{\eps\searrow 0}\underline v^{\mcR,k}_\eps(t,x)\geq\underline v^{\mcR,k}(t,x).
\end{align*}
\end{lem}

\noindent\emph{Proof.} We let $J^{\mcR}_l(t,x;u^\mcR,\nu)$ (resp.~$J^{\mcR,\eps}_l(t,x;u^\mcR,\nu)$) be the cost/reward obtained by limiting the number of points of $\mu$ (resp.~$\mu^t$) on the interval $[t,T]\times A$ to $l$ in the definitions of $J^{\mcR}_l(t,x;u^\mcR,\nu)$ (resp.~$J^{\mcR,\eps}_l(t,x;u^\mcR,\nu)$), so that
\begin{align*}
  J^{\mcR}_l(t,x;u^\mcR,\nu)&:=\E^\nu\Big[\psi(X^{t,x;u^\mcR,\eps}_T)+\int_t^{T}f(s,X^{t,x;u^\mcR,l}_s)ds - \sum_{j=1}^{N^\mcR}\ell(\tau^\mcR_j,X^{t,x;[u^\mcR]_{j-1},l}_{\tau^\mcR_j},\beta^\mcR_j)
  \\
  &\quad + \int_t^T\ett_{[\mu([t,s),A)<l]}\chi(s,X^{t,x;u^\mcR,j-1}_{s-},e)d\mu(ds,de)\Big|\mcF^{B,\mu}_{t}\Big],
\end{align*}
where
\begin{align*}
X^{t,x;u^\mcR,l}_s&=x+\int_t^s a(r,X^{t,x;u^\mcR,j}_r)dr+\int_t^s\sigma(r,X^{t,x;u^\mcR,l}_r)dB_r+\sum_{j=1}^N\ett_{[\tau_j\leq s]}\jmp(\tau^\mcR_j,X^{t,x;[u^\mcR]_{j-1},l}_{\tau^\mcR_j},\beta^\mcR_j)
\\
&\quad +\int_t^s\ett_{[\mu([t,r),A)<l]}\gamma(r,X^{t,x;u^\mcR}_{r-},e)d\mu(dr,de),\quad\forall s\in[t,T]
\end{align*}
(resp.
\begin{align*}
  J^{\mcR,\eps}_l(t,x;u^\mcR,\nu)&:=\E^\nu\Big[\psi(X^{t,x;u^\mcR,\eps,l}_T)+\int_t^{T}f(s,X^{t,x;u^\mcR,\eps,l}_s)ds - \sum_{j=1}^{N^\mcR}\ell(\tau^\mcR_j,X^{t,x;[u^\mcR]_{j-1},\eps,l}_{\tau^\mcR_j},\beta^\mcR_j)
  \\
  &\quad + \sum_{j=1}^{M^{\mu}_{t,T}\wedge l}\chi(\sigma^{t,\eps}_j,X^{t,x;u^\mcR,\eps,j-1}_{\sigma^{t,\eps}_j},\zeta^{t,\eps}_j)\Big|\mcF^{B,\mu}_{t}\Big]).
\end{align*}
Then, arguing as in Proposition~\ref{prop:trunk-vfs-approx} gives that
\begin{align*}
  \essinf_{\nu^S\in\mcV^S}\esssup_{u\in\mcU^{\mcR,k}}J^\mcR_l(t,x;u,\nu^S(u))\to \essinf_{\nu^S\in\mcV^S}\esssup_{u\in\mcU^{\mcR,k}}J^\mcR(t,x;u,\nu^S(u))
\end{align*}
and
\begin{align*}
  \essinf_{\nu^S\in\mcV^S}\esssup_{u\in\mcU^{\mcR,k,\eps}}J^{\mcR,\eps}_l(t,x;u,\nu^S(u))\to \underline v^{\mcR,k}_\eps(t,x)
\end{align*}
as $l\to\infty$, where the latter convergence holds uniformly in $\eps$. On the other hand, Lemma~\ref{app:lem:SDEflow} is independent of the underlying probability space and we conclude that
\begin{align*}
  \essinf_{\nu^S\in\mcV^S}\esssup_{u^\mcR\in\mcU^{\mcR,k,\eps}}J^{\mcR,\eps}_l(t,x;u^\mcR,\nu^S(u^\mcR))\to \essinf_{\nu^S\in\mcV^S}\esssup_{u^\mcR\in\mcU^{\mcR,k}}J^\mcR_l(t,x;u^\mcR,\nu^S(u^\mcR))
\end{align*}
as $\eps\searrow 0$, proving that
\begin{align*}
  \lim_{\eps\searrow 0}\underline v^{\mcR,k}_\eps(t,x)=\essinf_{\nu^S\in\mcV^S}\esssup_{u\in\mcU^{\mcR,k}}J^\mcR(t,x;u,\nu^S(u))
\end{align*}
and the assertion follows from \eqref{ekv:dual-lower-bnd}.\qed\\

\begin{cor}\label{cor:app-OC-dual}
For every $(t,x)\in[0,T]\times\R^d$, $k\in\bbN$ and $\varrho> 0$, there is an $\eps>0$ such that
\begin{align}\label{ekv:app-OC-dual}
  \underline v^{\mcR,k}(t,x)\leq\essinf_{\nu^S\in\mcV^{S}}J^{\mcR}(t,x;u^{\mcR,k,\eps},\nu^S(u^{\mcR,k,\eps}))+\varrho.
\end{align}
\end{cor}

\section{Proof of Proposition~\ref{prop:rel-rand-to-orig}\label{sec:proof}}
We introduce the set of extended setups for to original zero-sum game, denoted $\bigSET$, defined as the set of all
\begin{align*}
  \tilde\bbA=(\tilde\Omega,\tilde\mcF,\tilde\bbF,\tilde\Prob,\tilde W),
\end{align*}
where $(\tilde\Omega,\tilde\mcF,\tilde\Prob)$ is a complete probability space, $\tilde\bbF:=(\tilde\mcF_t)_{0\leq t\leq T}$ is a filtration satisfying the usual conditions and $\tilde W$ is a $(\tilde\bbF,\tilde\Prob)$-standard $d$-dimensional Brownian motion. We refer to $\bigSET$ as the set of \emph{extended setups} since $\tilde\bbF$ is not necessarily the augmented natural filtration of $\tilde W$.

For each $t\in[0,T]$, we defined the corresponding control set $\mcU_t^k(\mathbb A)$ as the set of all $u:=(\tau_j,\beta_j)_{j=1}^N$ such that $(\tau_j)_{j\geq 0}$ is a non-decreasing sequence of $\tilde\bbF$-stopping times with $\tau_0=t$, $\beta_j$ is $U$-valued and $\tilde\mcF_{\tau_j}$-measurable and $N:=\max\{j\geq 0:\tau_j<T\}\leq k$, $\tilde\Prob$-a.s.

The sets $\mcU_t(\bbA)$, $\mcU^{k,\eps}_t(\bbA)$, $\mcU^S_t(\bbA)$, $\mcA^k_t(\bbA)$, etcetera, are defined analogously and we let $\mcV(\bbA)$ (resp.~$\mcV^n(\bbA)$) be the set of $\Pred(\tilde\bbF)\otimes\mcB(A)$-measurable, bounded maps $\nu=\nu_t(\omega,e):[0,T]\times\Omega\times U\to [0,\infty)$.

For an arbitrary $\tilde\bbA=(\tilde\Omega,\tilde\mcF,\tilde\bbF,\tilde\Prob,\tilde W)\in\bigSET$ we extend the definition of the cost/reward functional $J$ and let for each $(t,x)\in [0,T)\times\R^d$ and impulse controls $u:=(\tau_j,\beta_j)_{j=1}^N\in\mcU_t(\tilde\bbA)$ and $\alpha:=(\eta_j,\theta_j)_{j=1}^M\in\mcA_t(\tilde\bbA)$, the \cadlag process $\tilde X^{t,x;u,\alpha}$ solve the SDE with impulses,
\begin{align}\nonumber
\tilde X^{t,x;u,\alpha}_s&=x+\int_t^s a(r,X^{t,x;u,\alpha}_r)dr+\int_t^s\sigma(r,X^{t,x;u,\alpha}_r)d\tilde W_r+\sum_{j=1}^N\ett_{[\tau_j\leq s]}\jmp(\tau_j,X^{t,x;[u]_{j-1},\alpha}_{\tau_j},\beta_j)
\\
&\quad +\sum_{j=1}^M\ett_{[\eta_j\leq s]}\gamma(\eta_j,X^{t,x;u,[\alpha]_{j-1}}_{\eta_j},\theta_j)\label{ekv:fwd-sde-tilde}
\end{align}
and arrive at the corresponding cost/reward functional
\begin{align*}
  J(t,x;u,\alpha,\tilde\bbA)&:=\tilde\E\Big[\psi(\tilde X^{t,x;u,\alpha}_T)+\int_t^{T}f(s,\tilde X^{t,x;u,\alpha}_s)ds
  \\
  &\quad- \sum_{j=1}^N\ell(\tau_j,\tilde X^{t,x;[u]_{j-1},\alpha}_{\tau_j},\beta_j)+ \sum_{j=1}^M\chi(\eta_j,\tilde X^{t,x;u,[\alpha]_{j-1}}_{\eta_j},\theta_j)\Big|\tilde\mcF_t\Big].
\end{align*}

\subsection{Part A\label{subsec:Prop-A}}
The novel work in \cite{Fuhrman15} considered a weak formulation of an optimal control problem where, in addition to a supremum over controls, the value function was obtained by taking the supremum over all conceivable probability spaces. In particular, this made it straightforward to prove that the value in the original control problem dominates that of the randomized version. An essential contribution made in \cite{Bandini18} was to consider a strong version of the control problem, where the probability space is fixed. As described below, considering the type of zero-sum games that we analyze does not lead to a significant increase of the complexity compared to the analysis in \cite{Bandini18}.

For the purpose of refining the argument in \cite{Bandini18} to apply in our setup, we augment the probability space $\bbA:=(\Omega,\mcF,\bbF,\Prob)$ from Section~\ref{sec:zsg} with the probability space $(\Omega^\mu,\mcF^\mu,\Prob^\mu)$ carrying a Poisson random measure $\mu$ with $\Prob^\mu$-compensator $\lambda$. We let $\tilde{\mathbb A}=(\tilde\Omega,\tilde\mcF,\tilde\bbF,\tilde\Prob,\tilde W,\tilde\mu)\in\bigSET^\mcR$, where $\tilde\Omega:=\Omega\times\Omega^\mu$, $\tilde \mcF$ is the $\Prob\otimes\Prob^\mu$ completion of $\mcF\otimes\mcF^\mu$, $\tilde\Prob$ denotes the extension of $\Prob\otimes\Prob^\mu$ to $\tilde\mcF$, while $\tilde W$ and $\tilde\mu$ are the canonical extensions of $W$ and $\mu$ to $\tilde\Omega$ and $\tilde\bbF$ is the augmented natural filtration generated by $\tilde W$.\\

\emph{Proof of Proposition~\ref{prop:rel-rand-to-orig}-A.} We let $\bar\bbA=(\tilde\Omega,\tilde\mcF,\bar\bbF,\tilde\Prob,\tilde W)\in\bigSET$, where $\bar\bbF$ is the $\tilde\Prob$-augmented natural filtration generated by $\tilde W$, and have, since $\bar\bbF$ is a subfiltration of $\tilde\bbF$, that
\begin{align*}
  \bar v^{\infty,l}(t,x)&= \esssup_{u^S\in\mcU^S_t(\bar{\bbA})}\essinf_{\alpha\in \mcA^l_t(\bar\bbA)}J(t,x; u^{S}(\alpha),\alpha,\tilde\bbA)
  \\
  &\leq \esssup_{\tilde u^S\in\mcU^S_t(\tilde{\bbA})}\essinf_{\alpha\in \mcA^l_t(\bar\bbA)} J(t,x; \tilde u^{S}(\alpha),\alpha,\tilde\bbA).
\end{align*}
For each $\varrho>0$ and $l\in\bbN$, there is thus a strategy $\tilde u^{S,\varrho}\in\mcU^{S}_t(\tilde\bbA)$, such that
\begin{align*}
  \bar v^{\infty,l}(t,x)&\leq \essinf_{\alpha\in \mcA^l_t(\bar{\bbA})}\tilde J(t,x;\tilde u^{S,\varrho}(\alpha),\alpha,\tilde\bbA)+\varrho.
\end{align*}
Taking the expectation on both sides and using that $\bar v^{\infty,l}(t,x)$ is deterministic, we find that
\begin{align*}
  \bar v^{\infty,l}(t,x)&\leq\inf_{\alpha\in \mcA^l_t(\bar{\bbA})}\tilde\E[J(t,x;\tilde u^{S,\varrho}(\alpha),\alpha,\tilde\bbA)]+\varrho.%\E\Big[\psi(X^{t,x;u^{S,\varrho}(\alpha),\alpha}_T)+\int_t^{T}f(s,X^{t,x;u^{S,\varrho}(\alpha),\alpha}_s)ds - \sum_{j=1}^N\ell(\tau_j,X^{t,x;[u^{S,\varrho}(\alpha)]_{j-1},\alpha}_{\tau_j},\beta_j)
  %\\
  %&\quad+ \sum_{j=1}^M\chi(\eta_j,X^{t,x;u^{S,\varrho}(\alpha),[\alpha]_{j-1}}_{\eta_j},\theta_j)\Big]+\varrho.
\end{align*}
The expression on the right-hand side represents an impulse control problem. Notably, this problem deviates from the standard archetype, as the map $\alpha \to \tilde X^{t,x;\tilde u^{S,\varrho}(\alpha),\alpha}$ does not conform to conventional regularity assumptions.

Conversely, the pivotal outcomes leading up to Proposition 4.2 in \cite{Bandini18} hinge exclusively upon the properties of measurability, without necessitating any further regularity constraints on this dataset. Consequently, we are able to replicate the reasoning delineated in the corresponding proofs, culminating in the deduction that for each $n\in\bbN$, we have
\begin{align*}
  \inf_{\alpha\in \mcA^l_t(\bar\bbA)}\tilde\E[J(t,x;\tilde u^{S,\varrho}(\alpha),\alpha,\tilde\bbA)] &=  \inf_{\nu\in\mcV^n_{\inf>0}(\tilde\bbA)}\inf_{\alpha\in \mcA^{l}_t(\tilde\bbA)}\tilde\E^\nu\Big[\psi(\tilde X^{t,x;\tilde u^{S,\varrho}(\alpha),\alpha}_T)+\int_t^{T}f(s,\tilde X^{t,x;\tilde u^{S,\varrho}(\alpha),\alpha}_s)ds
  \\
  &\quad- \sum_{j=1}^N\ell(\tau_j,\tilde X^{t,x;[\tilde u^{S,\varrho}(\alpha)]_{j-1},\alpha}_{\tau_j},\beta_j)+ \sum_{j=1}^M\chi(\eta_j,\tilde X^{t,x;\tilde u^{S,\varrho}(\alpha),[\alpha]_{j-1}}_{\eta_j},\theta_j)\Big].
\end{align*}
Now, as $\alpha^{t,\tilde\mu}:=(\tilde\sigma_{j},\tilde\zeta_{j})_{j=M^{\tilde\mu}_t}^{M^{\tilde\mu}_T\wedge M^{\tilde\mu}_t+l}\in\mcA^{l}_t(\tilde{\bbA})$, where $M^{\tilde\mu}_t:=\tilde\mu((0,t],A)$, and $u^{\mcR,\varrho}:=u^{S,\varrho}(\alpha^{t,\mu})\in\mcU^{\mcR}_t(\tilde{\bbA})$, it follows that
\begin{align*}
  \bar v^{\infty,l}(t,x)\leq \inf_{\nu\in\mcV^n_{\inf>0}(\tilde\bbA)}J^\mcR_l(t,x;u^{\mcR,\varrho},\nu,\tilde\bbA)+\varrho.
\end{align*}
Using that $\tilde\E^\nu[(M^{\tilde \mu}_{t,T})^2]\leq Cn^2$ whenever $\nu\in\mcV^n_{\inf>0}(\tilde\bbA)$ and Remark~\ref{rem:inf-g-0}, we find that the right-hand side converges to $\bar v^{\mcR,n}(t,x)+\varrho$ as $l\to\infty$ and, since $\varrho>0$ was arbitrary, we conclude that Proposition~\ref{prop:rel-rand-to-orig}-A holds.\qed\\

\subsection{Part B}
As in Section 4.3 of \cite{Fuhrman2020}, we consider the probability space $(\Omega,\mcF,\bbF,\Prob,W)$ defined in Section~\ref{sec:prel} and introduce an auxiliary probability space $( \Omega',\mcF',\Prob')$ on which lives real-valued random variables $(U^m_j,S^m_j)_{m,j\geq 1}$ and random measures $(\pi^\iota)_{\iota\geq 1}$ such that
\begin{enumerate}
  \item the $U^m_j$ are all uniformly distributed on $(0,1)$
  \item the probability distribution of $S^m_j$ admits a density $f^m_j$ with respect to the Lebesgue measure, that has support on the interval $((1-2^{1-j})/m,(1-2^{-j})/m)$, so that $0<S^m_1<S^m_2<\cdots<1/m$ for every $m\in\bbN$.
  \item every $\pi^\iota$ is a Poisson random measure on $(0,\infty)\times A$, with compensator $\iota^{-1}\lambda(de)dt$, with respect to its natural filtration;
  \item the random elements $U^m_j,S^{m'}_{j'},\pi^\iota$ are all independent.
\end{enumerate}
We let $\hat\Omega:=\Omega\times\Omega'$, let $\hat \mcF$ be the $\Prob\otimes\Prob'$ completion of $\mcF\otimes\mcF'$ and denote by $\check\Prob$ the extension of $\Prob\otimes\Prob'$ to $\hat\mcF$ (we reserve the notation $\hat\Prob$ for a different probability measure on $(\hat\Omega,\hat\mcF)$ to be defined later). Further, we let $\hat W,\hat\mu,\hat U^m_j,\hat S^{m}_{j}$ and $\hat \pi^\iota$ denote the canonical extensions of $W,\mu,U^m_j,S^{m}_{j}$ and $\pi^\iota$ to $\hat\Omega$.\\

We introduce two filtrations, and get two corresponding extended setups:
\begin{itemize}
  \item $\hat \bbA:=(\hat\Omega,\hat\mcF,\hat\bbF,\check\Prob,\hat W)$, where $\hat\bbF:=(\hat \mcF_t)_{{0\leq t\leq T}}$ is the $\check\Prob$-augmented natural filtration on $(\hat\Omega,\hat\mcF)$ generated by $\hat W$; and
  \item $\check \bbA:=(\hat\Omega,\hat\mcF,\check\bbF,\check\Prob,\hat W)$, where $\check\bbF=(\check \mcF_{t})_{0\leq t\leq T}$ is the $\check\Prob$-completion of the filtration $(\mcF_t\times \mcF')_{0\leq t\leq T}$.
\end{itemize}
Similar to the above, we have
\begin{align*}
  \underline v^{k,l}(t,x)&= \essinf_{\alpha^S\in \mcA^{S,l}_t(\hat\bbA)}\esssup_{u\in\mcU^k_t(\hat\bbA)}J(t,x; u,\alpha^S(u),\hat\bbA).
\end{align*}
The above notation is slightly clumsy. To improve it, we add a hat to objects defined over $\hat\bbA$, while check is used to denote objects living in $\check\bbA$ so that, for example, $\hat J(t,x;u,\alpha)=J(t,x;u,\alpha,\hat \bbA)$.\\

To facilitate the extension of $\hat\bbA$ to a setting for the dual game, we need the following definitions:
\begin{itemize}
\item For each $\alpha\in\check\mcA_t$, $\hat\bbF^{\alpha}=(\hat\mcF^{\alpha}_s)_{0\leq s\leq T}$ is the filtration on $(\hat\Omega,\hat\mcF)$ defined as $\hat\mcF^{\alpha}_s:=\sigma((\eta_j,\theta_j)_{j=1}^{M_s})$.
\item For each $t\in[0,T]$, $k\in\bbN$ and $\alpha\in\check\mcA_t$, we use the notation $\hat\mcU^{\hat W,\alpha,k}_t=\mcU^k_t((\hat\Omega,\hat\mcF,\hat\bbF\vee \hat\bbF^{\alpha},\check\Prob,\hat W))$ and let $\hat\mcU^{\hat W,\alpha,k,\eps}_t$ be the restriction to impulse controls where $(\eta_j,\theta_j)$ is valued in $\bbT^\eps\times\bar A^\eps$.
\end{itemize}

As we proceed, we fix $(t,x)\in[0,T]\times\R^d$, $k,l\in\bbN$ and $\varrho>0$ and note that by Proposition~\ref{prop:disc-vfs} there is for any sufficiently small $\eps>0$ a $\hat\alpha^{S,\varrho}\in\hat\mcA^{S,l,\eps}_t$ such that
\begin{align}\label{ekv:alpha-varrho}
  \esssup_{\hat u\in \hat\mcU^{k,\eps}_{t}}\hat J(t,x;\hat u,\hat\alpha^{S,\varrho}(\hat u))\leq \underline v^{k,l}(t,x)+\varrho,\quad\hat\Prob{\rm-a.s.}
\end{align}
The idea is to, given $\eps>0$ and $\iota\in\bbN$, use the sequences $(\hat U^m_j)_{j\in\bbN}$ and $(\hat S^{m}_{j})_{j\in\bbN}$ to randomize the strategy $\hat\alpha^{S,\varrho}$ and then add $\hat\pi^\iota$ to get a new strategy $\check \alpha^{S,\eps,\iota}:=(\check\eta^{S,\eps,\iota}_j,\check\theta^{S,\eps,\iota}_j)_{j=1}^{M^{S,\eps,\iota}}\in \check\mcA^{S}_t$ such that the sum of dirac measures $\check \alpha^{S,\eps,\iota}(u):=\sum_{j=1}^k\delta_{(\check\eta^{S,\eps,\iota}_j(u),\check\theta_j^{S,\eps,\iota}(u))}$ has $\check \Prob$-compensator that has a density $\check\nu^{S,\eps,\iota}(u)$ with respect to $\lambda(da)dt$ where $\check\nu^{S,\eps,\iota}$ is a strictly positive (not necessarily bounded) density strategy. Our proof of Proposition~\ref{prop:rel-rand-to-orig}-B hinges on the fact that this can be done in a fashion such that for small $\eps>0$ and large $\iota$, $\check \alpha^{S,\eps,\iota}(\hat u)$ is sufficiently ``close'' to $\hat\alpha^{S,\varrho}(\hat u)$ for any $\hat u\in\hat \mcU^k_t$, which we address in the next lemma.

To ease notation we introduce the time shift of $\alpha=(\eta_j,\theta_j)_{j=1}^M\in\check\mcA$ defined as $\alpha+\Delta t^\eps:=(\eta_j+\Delta t^\eps,\theta_j)_{j=1}^M$.

\begin{lem}\label{lem:will-conv}
For each $\eps>0$ and $\iota\in\bbN$, there is a $\check \alpha^{S,\eps,\iota}\in\check\mcA^{S}_t$ and a corresponding non-anticipative map $\check\nu^{S,\eps,\iota}:\hat\mcU^k\to ([t,T]\times\Omega\times A\to [0,\infty))$ such that for any $\hat u\in \hat\mcU^k_t$, we have
\begin{align}\label{ekv:proximity}
  \check\Prob[\hat\alpha^{S,\varrho}(\hat u)+\Delta t^\eps = \pi^\eps_A(\check \alpha^{S,\eps,\iota}(\hat u))]\geq 1-1/\iota
\end{align}
and the random measure on $[t,T]\times A$ corresponding to $\check \alpha^{S,\eps,\iota}(\hat u)$ has a $\check\Prob$-compensator with respect to the filtration $\hat\bbF\vee \hat\bbF^{\check \alpha^{S,\eps,\iota}(\hat u)}$ that is absolutely continuous with respect to $\lambda$ and takes the form
\begin{align*}
\check\nu^{S,\eps,\iota}_s[\hat u](\hat\omega,e)\lambda(de)ds
\end{align*}
where $\check\nu^{S,\eps,\iota}_t(\hat u)$ is $\Pred(\hat\bbF\vee \hat\bbF^{\check \alpha^{S,\eps,\iota}(\hat u)}) \otimes\mcB(A)$-measurable. Finally, for any $\hat u\in \hat\mcU^k_t$ the number of interventions of the impulse control $\check \alpha^{S,\eps,\iota}(\check u)$ (\ie the $\check \alpha^{S,\eps,\iota}(\check u)$-measure of the set $[t,T]\times A$), denoted $\check M^{S,\eps,\iota}(\hat u)$, has moments of all orders under $\check\Prob$.
\end{lem}

\noindent\emph{Proof.} As mentioned above, we first show existence of a non-negative (not necessarily bounded away from zero) map $\nu$ satisfying the first part of the lemma. For each $m\geq 1$, define the map $q^m:U\times [0,1]\to U:(b,da)\mapsto \frac{1}{\lambda({\bold B}(b,1/m))}\ett_{{\bold B}(b,1/m)}(a)\lambda(da)$ (where ${\bold B}(b,1/m)$ is the closed ball of radius $1/m$, centered at $b$) as in the proof of Lemma 4.4 of \cite{Fuhrman2020}. Now, define the sequence of strategies $(\check\alpha^{S,m}=(\check\eta^{S,m}_j,\check\theta^{S,m}_j)_{j=1}^{\hat M^{S,m}})_{m\geq 1}$ in $\hat\mcA^{S,l}_t$ as
\begin{align*}
  \check \eta^{S,m}_j(u):=\hat \eta^{S,\varrho}_{j}(\check u)+\hat S^m_j,\qquad \check\beta^{S,m}_j(\check u):=q^m(\hat\theta^{S,\varrho}_j(\check u),\hat U^m_j),\qquad \check M^{S,m}(\check u):=\inf\{j\geq 0:\check \eta^{S,m}_j(\check u)<T\}.
\end{align*}
According to Lemma A.11 in~\cite{Fuhrman15} the corresponding $\check\Prob$-compensator with respect to $\hat\bbF\vee \hat\bbF^{\check\alpha^{S,m}(u)}$ is given by the explicit formula
\begin{align*}
  \sum_{j=1}^l\ett_{(\hat \eta^{S,\varrho}_{j}(\check u)\vee \check \eta^{S,\varrho}_{j-1}(\check u) ,\check\eta^{S,\varrho}_{j}(u)]}(t)q^m(\hat\theta^{S,\varrho}_j(u),da)\frac{f^m_j(s-\hat\eta^{S,\varrho}_j(u))}{1-F^m_j(s-\hat\eta^{S,\varrho}_j(u))}ds,
\end{align*}
with $F^m_j(s):=\int_{-\infty}^s f^m_j(r)dr$. For each $m\in\bbN$, this compensator is clearly $\Pred(\hat\bbF\vee \hat\bbF^{\check\alpha^{S,m}(\check u)}) \otimes\mcB(A)$-measurable.

Moreover, there is a $m'\in\bbN$ such that the densities of the random variables in the sequence $(S^{m}_j)_{j\in\bbN}$ all have support in $(0,\Delta t^\eps)$ and $|\check\theta^{S,m}_j(\hat u)-\hat\theta^{S,\varrho}_j(\hat u)|\leq \eps$ for all $\hat u\in \hat\mcU^k_t$ and $j=1,\ldots,l$, whenever $m\geq m'$. We find that
\begin{align*}
\hat\alpha^{S,\varrho}(\hat u)+\Delta t^\eps = \pi^\eps_A(\check \alpha^S(\hat u)),\quad \hat\Prob-\text{a.s.}
\end{align*}
whenever $m\geq m'$.

It remains to modify $\check \alpha^{S,m}$ so that the corresponding density with respect to $\lambda$ is bounded away from $0$ on $[0,T]$. We, therefore, consider the strategy $\check \alpha^{S,\eps,\iota}:=(\check\eta_j^{S,\eps,\iota},\check\theta^{S,\eps,\iota}_j)_{j=1}^{\check M^{S,\eps,\iota}}$ corresponding to the random measure $\check \alpha^{S,m}(\hat u)+\hat \pi^\iota$ and note that the number of interventions of $\check \alpha^{S,\eps,\iota}(\hat u)$ on $[0,T]$ is bounded by $l+\hat M^{\hat \pi^\iota}$, where $\hat M^{\hat \pi^\iota}:=\hat \pi^\iota([0,T],A)$ is Poisson distributed with parameter $\lambda(A)T/\iota$ under $\hat \Prob$. In particular, this gives that $\hat M^{\hat \pi^\iota}$ and then also $l+\hat M^{\hat \pi^\iota}$ has moments of all orders under $\check\Prob$.

Moreover, $\check \alpha^{S,\eps,\iota}(\hat u)$ has a $\check\Prob$-compensator with respect to the filtration $\hat\bbF\vee \hat\bbF^{\check\alpha^{S,\eps,\iota}(\hat u)}$ that has a density with respect to $\lambda$ which is bounded from below by $1/\iota$.\qed\\

Given $\hat u\in \hat\mcU^{k}_{t}$ we now construct a new auxiliary probability measure, under which the corresponding $\check\alpha^{S,\eps,\iota}(\hat u)$ is a Poisson random measure with intensity $\lambda$. For any $\eps>0$, $\iota\in\bbN$ and any impulse control $\hat u\in \hat\mcU^{k,\eps}_{t}$, the response $\check\alpha^{S,\eps,\iota}(\hat u)$ induces a random measure on $[0,T]\times A$ with corresponding control $(\check\sigma_j,\check\zeta_j)_{j\geq 1}$ and density $\check\nu=\check\nu^{S,\eps,\iota}(\hat u)$, with respect to $\lambda$. Since $\check\nu$ is bounded from below we find that
\begin{align*}
\hat G_s:=\exp\Big(\int_{0}^s\int_A(1-(\check\nu_r(e))^{-1})\lambda(de)dr\Big)\prod_{\check\sigma_j\leq s}(\check\nu_{\check\sigma_j}(\check\zeta_j))^{-1}
\end{align*}
is a strictly positive martingale with respect to the filtration $\hat\bbF\vee\hat\bbF^{\check\alpha^{S,\eps,\iota}(\hat u)}$ under $\check\Prob$. We define the equivalent probability measure $\hat\Prob$ on $(\hat\Omega,\hat\mcF)$ as $d\hat\Prob=\hat G_Td\check\Prob$. By the Girsanov theorem, $\check\alpha^{S,\eps,\iota}(\hat u)$ has $\check\Prob$-compensator $\lambda(de)ds$ with respect to the filtration $\hat\bbF\vee \hat\bbF^{\check\alpha^{S,\eps,\iota}(\hat u)}$.

We then consider the probability space
\begin{align*}
  \bbA^{\eps,\iota}(\hat u):=(\hat\Omega,\hat\mcF,\hat\bbF\vee \hat\bbF^{\check\alpha^{S,\eps,\iota}(\hat u)},\hat\Prob,\hat W,\check\alpha^{S,\eps,\iota}(\hat u))\in\bigSET^\mcR.
\end{align*}
Despite the fact that $\check\nu$ is generally not bounded (and we cannot assume that $\check\nu$ belongs to $\check\mcV^n(\bbA^{\eps,\iota}(\hat u))$ for some $n\in\bbN$) we still have a Dol{\'e}ans-Dade exponential
\begin{align*}
\hat\kappa^{\check\nu}_s:=\exp\Big(\int_{0}^s\int_A(1-\check\nu_r(e))\lambda(de)dr\Big)\prod_{\hat\sigma_j\leq s}\check\nu_{\hat\sigma_j}(\hat\zeta_j)
\end{align*}
for which $\hat \E[\hat\kappa^{\check\nu}_T]=\check \E[\hat G_T\hat\kappa^{\check\nu}_T]=1$, proving that $\hat\kappa^{\check\nu}$ is a $\hat\Prob$-martingale and we can define a corresponding probability measure $\hat\Prob^{\check\nu}$ on $(\hat\Omega,\hat\mcF)$ as $d\hat\Prob^{\check\nu}:=\hat\kappa^{\check\nu}_Td\hat\Prob$ and since $\hat M_T\hat\kappa^{\check\nu}_T\equiv 1$, we conclude that $\hat\Prob^{\check\nu}=\check\Prob$ on $(\hat\Omega,\hat\mcF)$. In particular, %letting
%\begin{align*}
%\check X^{t,x;u}_s&=x+\int_t^s a(r,\check X^{t,x;u}_r)dr+\int_t^s\sigma(r,\check X^{t,x;u}_r)d\hat W_r+\int_t^s\!\!\!\int_A\gamma(r,\check X^{t,x;u}_{r-},e)\check\alpha^{S,\eps,\iota}[\hat u](dr,de)
%\\
%&\quad +\sum_{j=1}^N\ett_{s\leq\tau_j}\jmp(\tau_j,\check X^{t,x;[u]_{j-1}}_{\tau_j},\beta_j),\quad\forall s\in[t,T],
%\end{align*}
we get that
\begin{align}\label{ekv:prim-dual-eq}
  \hat J(t,x;\hat u,\check\alpha^{S,\eps,\iota}(\hat u))= J^{\mcR}(t,x;\hat u,\check \nu^{S,\eps,\iota}(\hat u),\bbA^{\eps,\iota}(\hat u)),\quad \hat\Prob-\text{a.s.}
\end{align}
%, where
%\begin{align*}
%\hat J^{\mcR,\check\alpha^{S,\eps,\iota}(\hat u)}(t,x;u,\nu) & := \E^\nu\Big[\psi(\check X_T^{t,x;u})+\int_t^{T}f(s,\check X_s^{t,x;u})ds - \sum_{j=1}^N\ell(\tau_j,\check X^{t,x;[u]_{j-1}}_{\tau_j},\beta_j)
%\\
%&\quad + \int_t^T\int_A\chi(s,\check X^{t,x;u}_{s-},e)\check\alpha^{S,\eps,\iota}[\hat u](ds,de)\Big].
%\end{align*}
%We have the following corollary result to Lemma~\ref{lem:will-conv}:
%\begin{cor}
%For $\eps>0$ sufficiently small and $\iota\in\bbN$ sufficiently large, we have
%\begin{align}\label{ekv:below-underline-Y}
%  \hat J^{\mcR}(t,x;\hat u,\check \nu^{S,\eps,\iota}(\hat u),\bbA^{\eps,\iota}(\hat u))\leq \underline v^{k,l}(t,x)  + 2\varrho,
%\end{align}
%for any $\hat u \in \hat\mcU^{k,\eps}_t$.
%\end{cor}
%
%\emph{Proof.}

There are two main issues that have to be dealt with before we can apply \eqref{ekv:proximity} along with continuity to finish the proof of Proposition~\ref{prop:rel-rand-to-orig}. First, the statement in Lemma~\ref{lem:will-conv} only holds for $\hat u \in \hat\mcU^{k,\eps}_t$ whereas in the present framework $\underline v^{\mcR,k}(t,x)$ is obtained by maximising over $\check u^\mcR \in \hat\mcU^{\mcR,k,\eps}_t(\bbA^{\eps,\iota}(\hat u))$. Moreover, $\check \nu^{S,\eps,\iota}(\hat u)$ is not necessarily bounded and needs to be approximated by a sequence of bounded densities. We start with the latter issue.
\begin{lem}\label{lem:approx-dens}
For any $(t,x)\in[0,T]\times\R^d$, $\eps>0$, $\iota\in\bbN$ and $\hat u \in \hat\mcU^{k,\eps}_t$, let $\check\alpha^{S,\eps,\iota}(\hat u)$ and $\nu^{S,\eps,\iota}(\hat u)$ be as in the statement of Lemma~\ref{lem:will-conv}. Then,
\begin{align}\label{ekv:approx-dens}
   J^{\mcR}(t,x;\hat u,\check \nu^{S,\eps,\iota}(\hat u),\bbA^{\eps,\iota}(\hat u))=\lim_{n\to\infty} J^{\mcR}(t,x;\hat u,\check \nu^{S,\eps,\iota,n}(\hat u),\bbA^{\eps,\iota}(\hat u))
\end{align}
in $L^1(\hat\Omega,\hat\mcF,\hat\Prob)$, where $\check\nu^{S,\eps,\iota,n}:=n\wedge\check\nu^{S,\eps,\iota}$.
\end{lem}

\emph{Proof.} To simplify notation, we let $\check\bbA:=\bbA^{\eps,\iota}(\hat u)$, $\check\alpha:=\check\alpha^{S,\eps,\iota}(\hat u)$, $\check \nu^n:=\check \nu^{S,\eps,\iota,n}(\hat u)$ and $\check \nu:=\check \nu^{S,\eps,\iota}(\hat u)$. Moreover, we introduce
\begin{align*}
  \Phi(t,x):=\psi(\check X_T^{t,x;\hat u})+\int_t^{T}f(s,\check X_s^{t,x;\hat u})ds - \sum_{j=1}^N\ell(\hat\tau_j,\check X^{t,x;[\hat u]_{j-1}}_{\hat\tau_j},\hat\beta_j) + \int_t^T\int_A\chi(s,\check X^{t,x;\hat u}_{s-},e)\check\alpha(ds,de)
\end{align*}
that satisfies the bound
\begin{align*}
  |\Phi(t,x)|\leq C(1+k+\hat M^{\check \alpha}_{t,T})(1+\|\check X^{t,x;\hat u}\|^\rho_T)=:\bar\Phi(t,x)
\end{align*}
%Lemma \ref{lem:stop-disc} implies that
and get that
\begin{align*}
  \hat \E\big[|J^{\mcR}(t,x;\hat u,\check \nu^{n},\check \bbA)  - J^{\mcR}(t,x;\hat u,\check \nu,\check \bbA)|\big]&\leq \hat\E\big[|\hat\kappa^{\check\nu^n}_T-\hat\kappa_T^{\check\nu}|\bar\Phi(t,x)\big].
\end{align*}
On the other hand, with $E_K:=\{\omega:\|\check X^{t,x;\hat u}\|_T+M^{\check \alpha}_{t,T}\leq K\}$ we get
\begin{align*}
  \hat\E\big[|\hat\kappa^{\check\nu^n}_T-\hat\kappa_T^{\check\nu}|\bar\Phi(t,x)\big]&\leq \hat\E\big[\ett_{E_K}|\hat\kappa^{\check\nu^n}_T-\hat\kappa_T^{\check\nu}|\bar\Phi(t,x)\big]
  \\
  &\quad + \hat\E^{\check\nu^n}\big[\ett_{E_K^c}\bar\Phi(t,x)\big]+ \hat\E^{\check\nu}\big[\ett_{E_K^c}\bar\Phi(t,x)\big].
\end{align*}
Concerning the last two terms we have
\begin{align*}
  \hat\E^{\check\nu^n}\big[\ett_{E_K^c}\bar\Phi(t,x)\big]&= C\hat\E^{\check\nu^n}\big[\ett_{E_K^c}(1+\hat M^{\check \alpha}_{t,T})(1+\|\check X^{t,x;\hat u}\|^\rho_T)\big]
  \\
  &\leq \frac{C}{K}\hat\E^{\check\nu^n}\big[(\|\check X^{t,x;\hat u}\|_T+\hat M^{\check \alpha}_{t,T})(1+\hat M^{\check \alpha}_{t,T})(1+\|\check X^{t,x;\hat u}\|^\rho_T)\big]
  \\
  &\leq \frac{C}{K},
\end{align*}
where $C>0$ can be chosen independent of $n$ and $K$, and similarly for the last term. Concerning the first term, we have
\begin{align*}
\hat\E\big[\ett_{E_K}|\hat\kappa^{\check\nu^n}_T-\hat\kappa^{\check\nu}_T|\bar\Phi(t,x)\big] &\leq C\hat\E\big[|\hat\kappa^{\check\nu^n}_T-\hat\kappa^{\check\nu}_T|(1+K^{\rho+1})\big].
\end{align*}
On the other hand, by dominated convergence we have
\begin{align*}
\int_{0}^T\int_A(1-\check\nu^n_r(e))\lambda(de)dr \to \int_{0}^T\int_A(1-\check\nu_r(e))\lambda(de)dr,
\end{align*}
$\hat\Prob$-a.s.~and
\begin{align*}
  \prod_{\hat\sigma_j\leq s}\check\nu^n_{\hat\sigma_j}(\hat\zeta_j)\to \prod_{\hat\sigma_j\leq s}\check\nu^n_{\hat\sigma_j}(\hat\zeta_j),
\end{align*}
$\hat\Prob$-a.s.~as the number of terms in the product is $\hat\Prob$-a.s.~finite, effectively implying that $\hat\kappa^{\check\nu^n}_T\to \hat\kappa^{\check\nu}_T$, $\hat\Prob$-a.s. Since $K>0$ was arbitrary, we conclude that
\begin{align*}
  J^{\mcR}(t,x;\hat u,\check \nu^{n},\check \bbA)  \to J^{\mcR}(t,x;\hat u,\check \nu,\check \bbA)
\end{align*}
in $L^1(\hat\Omega,\hat\mcF,\hat\Prob)$ as $n\to\infty$.\qed\\

We now tend to the first issue and have the following:
\begin{lem}\label{lem:approx-OC}
Let $\hat u^{k,\eps}:=(\hat\tau^{k,\eps}_j,\hat\beta^{k,\eps}_j)_{j=1}^{\hat N^{k,\eps}}\in\hat\mcU^{k,\eps}_t$ be defined as
\begin{itemize}
  \item $\hat\tau^{k,\eps}_j:=\inf\{s\in\bbT^\eps_{\hat\tau^{k,\eps}_{j-1}}: \underline v^{\mcR,k-j+1}_{\eps}(s,X^{t,x;[\hat u^{k,\eps}]_{j-1},\hat\alpha^{S,\varrho}(\hat u^{k,\eps})+\Delta t^\eps})$\\$=\supOP_\eps \underline v^{\mcR,k-j}_{\eps}(s,X^{t,x;[\hat u^{k,\eps}]_{j-1},\hat\alpha^{S,\varrho}(\hat u^{k,\eps})+\Delta t^\eps})\}$,
  \item $\hat\beta^{k,\eps}_j\in\argmax_{b\in \bar U^\eps}\{ \underline v^{\mcR,k-j}_{\eps}(\hat\tau^{k,\eps}_j,X^{t,x;[\hat u^{k,\eps}]_{j-1},\hat\alpha^{S,\varrho}(\hat u^{k,\eps})+\Delta t^\eps}_{\hat\tau^{k,\eps}_j}$\\$ + \jmp(\hat\tau^{k,\eps}_j,X^{t,x;[\hat u^{k,\eps}]_{j-1},\hat\alpha^{S,\varrho}(\hat u^{k,\eps})+\Delta t^\eps}_{\hat\tau^{k,\eps}_j},b)) - \ell(\hat\tau^{k,\eps}_j,X^{t,x;[\hat u^{k,\eps}]_{j-1},\hat\alpha^{S,\varrho}(\hat u^{k,\eps})+\Delta t^\eps}_{\hat\tau^{k,\eps}_j},b)\}$
\end{itemize}
with $\hat\tau^{k,\eps}_0:=t$ and $\hat N^{k,\eps}:=\max\{j:\tau^{k,\eps}_j<T\}$. Then, there is a $\eps_0>0$ and a $\iota_0\in\bbN$ such that
\begin{align*}
  \underline v^{\mcR,k}(t,x)\leq J^{\mcR}(t,x;\hat u^{k,\eps},\check\nu^{S,\eps,\iota,n}(\hat u^{k,\eps}),\bbA^{\eps,\iota}(\hat u^{k,\eps}))+\varrho,\quad \hat\Prob-\text{a.s.}
\end{align*}
for all $n\in\bbN$, whenever $\eps\in (0,\eps_0)$ and $\iota\geq\iota_0$.
\end{lem}

\emph{Proof.} Since $\hat u^{k,\eps}\in\hat\mcU^{k,\eps}_t$, we can define a random measure $\check\alpha:=\check\alpha^{S,\eps,\iota}(\hat u^{k,\eps})=\sum_{j\geq 1}\delta_{(\check\eta_j,\check\theta_j)}$ and let $\check\alpha^{\eps}=(\check\eta^\eps_j,\check\theta^\eps_j)_{j=1}^{\check M}=\pi^\eps_A((\check\eta_j,\check\theta_j)_{j=1}^{\check M})$ and set $\check R^{t,x;u,\eps}:=\lim_{j\to\infty}\check R^{t,x;u,\eps,j}$, the pointwise limit of the solutions to
\begin{align*}
\check R^{t,x;u,\eps,j}_s&=x+\int_t^s a(r,\check R^{t,x;u,\eps,j}_r)dr+\int_t^s\sigma(r,\check R^{t,x;u,\eps,j}_r)d\hat W_r+\sum_{i=1}^N\ett_{[\tau_i\leq s]}\jmp(\tau_i,\check R^{t,x;[u]_{i-1},\eps,j}_{\tau_i},\beta_i)
\\
&\quad +\sum_{i=1}^{\check M\wedge j}\ett_{[\check\eta^{\eps}_i\leq s]}\gamma(\check\eta^{\eps}_i,R^{t,x;u,\eps,i-1}_{\check\eta^{\eps}_i},\check\theta^{\eps}_i),\quad\forall s\in [t,T].
\end{align*}
Corollary~\ref{cor:app-OC-dual} gives that if $\eps>0$ is chosen sufficiently small, the impulse control $\check u^{k,\eps}:=(\check \tau^{k,\eps}_j,\check \beta^{k,\eps}_j)_{j=1}^{\check N^{k,\eps}}\in\hat\mcU^{\hat W,\pi^\eps_A(\check\alpha),k,\eps}_t$ defined as
\begin{itemize}
  \item $\check \tau^{k,\eps}_j:=\inf\{s\in\bbT^\eps_{\check \tau^{k,\eps}_{j-1}}:\underline v^{\mcR,k-j+1}_\eps(s,\check R^{t,x;[\check u^{k,\eps}]_{j-1},\eps}_s)=\supOP \underline v^{\mcR,k-j}_\eps(s,\check R^{t,x;[\check u^{k,\eps}]_{j-1},\eps}_s)\}$
  \item $\check \beta^{k,\eps}_j\in\argmax_{b\in U}\{\underline v^{\mcR,k-j}_\eps(\check \tau^{k,\eps}_j,\check R^{t,x;[\check u^{k,\eps}]_{j-1},\eps}_{\check \tau^{k,\eps}_j} + \jmp(\check \tau^{k,\eps}_j,\check R^{t,x;[\check u^{k,\eps}]_{j-1},\eps}_{\check \tau^{k,\eps}_j},b)) - \ell(\check \tau^{k,\eps}_j,R^{t,x;[\check u^{k,\eps}]_{j-1},\eps}_{\check \tau^{k,\eps}_j},b)\}$
\end{itemize}
with $\check \tau^{k,\eps}_0:=t$ and $\check N^{k,\eps}:=\max\{j:\check\tau^{k,\eps}_j<T\}$ satisfies
\begin{align*}
  \underline v^{\mcR,k}(t,x)\leq J^{\mcR}(t,x;\check u^{k,\eps},\check\nu^{S}(\hat u^{k,\eps}),\bbA^{\eps,\iota}(\hat u^{k,\eps}))+\varrho/2,\quad \hat\Prob-\text{a.s.}
\end{align*}
for all $\check\nu^{S}\in\mcV^S(\bbA^{\eps,\iota}(\hat u^{k,\eps}))$. On the other hand, the control $\hat\alpha^{S,\varrho}(\hat u^{k,\eps})+\Delta t^\eps$ and the sequence $(\check\eta^{\eps}_j,\check\theta^{\eps}_j)_{j=1}^{\check M}$ with $\check M:=\check\alpha([t,T],A)$, and then by definition also $\hat u^{k,\eps}$ and $\check u^{k,\eps}$, coincide on the set $\Gamma:=\{\omega:\hat M^{\hat \pi^\iota}_{t,T}=0\}$ and we have, with $\check X^{t,x;u}$ now the unique solution to
\begin{align*}
\check X^{t,x;u}_s&=x+\int_t^s a(r,\check X^{t,x;u}_r)dr+\int_t^s\sigma(r,\check X^{t,x;u}_r)d\hat W_r+\int_t^s\!\!\!\int_A\gamma(r,\check X^{t,x;u}_{r-},e)\check\alpha(dr,de)
\\
&\quad +\sum_{j=1}^N\ett_{[s\leq\tau_j]}\jmp(\tau_j,\check X^{t,x;[u]_{j-1}}_{\tau_j},\beta_j),\quad\forall s\in[t,T],
\end{align*}
that
\begin{align*}
  &J^{\mcR}(t,x;\check u^{k,\eps},\check\nu^{S}(\hat u^{k,\eps}),\bbA^{\eps,\iota}(\hat u^{k,\eps}))-J^{\mcR}(t,x;\hat u^{k,\eps},\check\nu^{S}(\hat u^{k,\eps}),\bbA^{\eps,\iota}(\hat u^{k,\eps}))
  \\
  &\leq \check\E\Big[\ett_{\Gamma^c}\big(\Psi(\check X^{t,x;\check u^{k,\eps}}_T)-\Psi(\check X^{t,x;\hat u^{k,\eps}}_T)+\int_t^{T} (|f(r,\check X^{t,x;\check u^{k,\eps}}_r)|+|f(r,\check X^{t,x;\hat u^{k,\eps}}_r)|)dr
\\
&\quad + \sum_{j=1}^{\check N^{k,\eps}}\ell(\check\tau^{k,\eps}_j,\check X^{t,x;[\check u^{k,\eps}]_{j-1}}_{\check\tau^{k,\eps}_j},\check\beta^{k,\eps}_j) + \sum_{j=1}^{\hat N^{k,\eps}}\ell(\hat\tau^{k,\eps}_j,\check X^{t,x;[\hat u^{k,\eps}]_{j-1}}_{\hat\tau^{k,\eps}_j},\hat\beta^{k,\eps}_j)
\\
&\quad + \int_t^T\!\!\!\int_A(\chi(r,\check X^{t,x;\check u^{k,\eps}}_{r-},e) - \chi(r,\check X^{t,x;\hat u^{k,\eps}}_{r-},e))\check\alpha(dr,de)\big)\Big|\mcF^{\hat W,\hat\alpha}_t\Big]
\\
&\leq C\check\E\Big[\ett_{\Gamma^c}\check M\big(1+\|\check X^{t,x;\check u^{k,\eps}}\|^\rho_T+\|\check X^{t,x;\hat u^{k,\eps}}\|^\rho_T\big)\Big|\mcF^{\hat W,\hat\alpha}_t\Big]
\\
&\leq C\check\E\Big[\hat M^{\hat \pi^\iota}_{t,T}\big(1+\|\check X^{t,x;\check u^{k,\eps}}\|^\rho_T+\|\check X^{t,x;\hat u^{k,\eps}}\|^\rho_T\big)\Big|\mcF^{\hat W,\hat\alpha}_t\Big]
\\
&\leq C(1+|x|^\rho)\check\E\big[(\hat M^{\hat \pi^\iota}_{t,T})^2\big|\mcF^{\hat W,\hat\alpha}_t\big]^{1/2},
\end{align*}
where the left hand side is less than $\varrho/2$ for $\iota$ sufficiently large as $\check\E\big[(\hat M^{\hat \pi^\iota}_{t,T})^2\Big|\mcF^{\hat W,\hat\alpha}_t\big]=(T-t)\lambda(A)/\iota\to 0$, as $\iota\to\infty$. Since $\check\nu^{S,\eps,\iota,n}\in\mcV^S(\bbA^{\eps,\iota}(\hat u^{k,\eps}))$ for each $n\in\bbN$, the desired result follows.\qed\\

%\begin{center}
%\fbox{%
%  \parbox{0.85\textwidth}{
%  Steps are as follows where we for $(t,x,k,l)\in [0,T]\times\R^d\times\bbN^2$ and $\varrho>0$, have:
%\begin{enumerate}
%  \item There is an $\eps_0>0$ such that for any $\eps\in (0,\eps_0]$, the strategy $\hat\alpha^{S,\varrho}\in\hat\mcA^{S,l,\eps}_t$ defined as the canonical extension of $\alpha^{S,k,l,\eps}$ satisfies
%  \begin{align*}
%    \hat J(t,x;\hat u,\hat\alpha^{S,\varrho}(\hat u))\leq \underline v^{k,l}(t,x)+\varrho.
%  \end{align*}
%  for any $\hat u\in\hat\mcU^{k,\eps}_t$.
%  \item Letting $\check\alpha^{S,m,\iota}$ correspond to $\hat \alpha^S=\hat\alpha^{S,\varrho}$ as in the proof of Lemma~\ref{lem:will-conv} we find that for any $\hat u\in \hat\mcU^k_t$ the impulse control $\check\alpha^{S,m,\iota}(\hat u)$ defines a random measure on $[t,T]\times A$. We pick $\hat u =\hat u^{k,l,1/m}$.
%  \item Then, by Corollary~\ref{cor:app-OC-dual} and approximation result $\check\nu^{S,n}:=\check\nu^S \wedge n$ and $n\to\infty$, there is a $m_0$ such that
%  \begin{align*}
%    \hat J(t,x;\check u^{\mcR,k,1/m},\check\alpha^{S,m,\iota}(\hat u^{k,l,1/m}))\geq \underline v^{\mcR,k}(t,x)-\varrho
%  \end{align*}
%  whenever $m\geq m_0$.
%  \item On the other hand,
%  \begin{align*}
%    \hat J(t,x;\check u^{\mcR,k,1/m},\check\alpha^{S,m,\iota}(\hat u^{k,l,1/m}))\leq \hat J(t,x;\hat u^{k,l,1/m},\hat\alpha^{S,\varrho}(\hat u^{k,l,1/m})) +\varrho
%  \end{align*}
%\end{enumerate}
%}}
%\end{center}

\emph{Proof of Proposition~\ref{prop:rel-rand-to-orig}-B.} Combining Lemmas~\ref{lem:approx-dens} and \ref{lem:approx-OC} gives that for all sufficiently large $\iota\in\bbN$, and all sufficiently small $\eps>0$,
\begin{align*}
  \underline v^{\mcR,k}(t,x)&\leq \lim_{n\to\infty}\hat\E\big[J^{\mcR}(t,x;\hat u^{k,\eps},\check\nu^{S,\eps,\iota,n}(\hat u^{k,\eps}),\bbA^{\eps,\iota}(\hat u^{k,\eps}))\big]+\varrho
  \\
  &=\hat\E\big[J^{\mcR}(t,x;\hat u^{k,\eps},\check\nu^{S,\eps,\iota}(\hat u^{k,\eps}),\bbA^{\eps,\iota}(\hat u^{k,\eps}))\big]+\varrho.
\end{align*}
Moreover, since $\hat u^{k,\eps}\in\hat\mcU^{k}_t$ we get by combining \eqref{ekv:alpha-varrho} and \eqref{ekv:prim-dual-eq} that for any $\iota\in\bbN$ and $\eps>0$,
\begin{align*}
  \underline v^{k,l}(t,x)&\geq \hat J(t,x;\hat u^{k,\eps},\hat\alpha^{S,\varrho}(\hat u^{k,\eps})) - \varrho
  \\
  &= J^{\mcR}(t,x;\hat u^{k,\eps},\check \nu^{S,\eps,\iota}(\hat u^{k,\eps}),\bbA^{\eps,\iota}(\hat u^{k,\eps}))+(\hat J(t,x;\hat u^{k,\eps},\hat\alpha^{S,\varrho}(\hat u^{k,\eps}))-\hat J(t,x;\hat u^{k,\eps},\check\alpha^{S,\eps,\iota}(\hat u^{k,\eps})))-\varrho.
\end{align*}
Hence, whenever $\iota\in\bbN$ is sufficiently large and $\eps>0$ is sufficiently small, we have
\begin{align*}
  \underline v^{\mcR,k}(t,x)-\underline v^{k,l}(t,x)&\leq \hat\E\big[\hat J(t,x;\hat u^{k,\eps},\hat\alpha^{S,\varrho}(\hat u^{k,\eps}))-\hat J(t,x;\hat u^{k,\eps},\check\alpha^{S,\eps,\iota}(\hat u^{k,\eps}))\big]+2\varrho
  \\
  &=\check\E\big[\hat J(t,x;\hat u^{k,\eps},\hat\alpha^{S,\varrho}(\hat u^{k,\eps}))-\hat J(t,x;\hat u^{k,\eps},\check\alpha^{S,\eps,\iota}(\hat u^{k,\eps}))\big]+2\varrho,
\end{align*}
where the last equality follows from the fact that $\hat \E$ and $\check\E$ coincide on $\hat\mcF_t$. In order to adequately deal with the remaining term, we separate the jumps in $\check\eta^{S,\eps,\iota}_i(\hat u^{k,\eps})$ into one part that comes from $\check\eta^{S,\eps}_i(\hat u^{k,\eps})$ (or at least could, from what we observe, originate from this component) and one part containing the remainder. We, thus, introduce the admissible setup $\bbA^{\eps,\iota}=(\hat\Omega,\hat\mcF,\hat \bbF\vee\hat\mcF^{\check \alpha^{S,\eps,\iota}(\hat u^{k,\eps})},\check\Prob,\hat W)\in\bigSET$ and the control $\tilde\alpha^{\eps,\iota}=(\tilde\eta^{\eps,\iota}_j,\tilde\theta^{\eps,\iota}_j)_{j=1}^{\tilde M^{\eps,\iota}}\in\mcA^{l}_t(\bbA^{\eps,\iota})$, where
\begin{itemize}
  \item $\tilde\eta^{\eps,\iota}_j:=\min_{i\geq 1}\{\check\eta^{S,\eps,\iota}_i(\hat u^{k,\eps}):\check\eta^{S,\eps,\iota}_i(\hat u^{k,\eps})\geq \hat\eta^{S,\varrho}_j(\hat u^{k,\eps})\text{ and } \pi^{A,\eps}_2(\check\theta^{S,\eps,\iota}_i(\hat u^{k,\eps}))= \hat\theta^{S,\varrho}_j(\hat u^{k,\eps})\}$ and
  \item $\tilde\theta^{\eps,\iota}_j:=\check\theta^{S,\eps,\iota}_j(\hat u^{k,\eps})$,
\end{itemize}
for $j=1,\ldots,\tilde M^{\eps,\iota}$, where $\tilde M^{\eps,\iota}:=\sup\{j\geq 1:\tilde\eta^{\eps,\iota}_j<T\}\vee 0$. With this definition, we may write
\begin{align*}
  &|\hat J(t,x;\hat u,\hat\alpha^{S,\varrho}(\hat u^{k,\eps}))-\hat J(t,x;\hat u^{k,\eps},\check \alpha^{S,\eps,\iota}(\hat u^{k,\eps}))|
  \\
  &\leq |\hat J(t,x;\hat u^{k,\eps},\check \alpha^{S,\eps,\iota}(\hat u^{k,\eps}))-\hat J(t,x;\hat u^{k,\eps},\tilde\alpha^{\eps,\iota})|+|\hat J(t,x;\hat u^{k,\eps},\tilde\alpha^{\eps,\iota})-\hat J(t,x;\hat u,\hat\alpha^{S,\varrho}(\hat u^{k,\eps}))|.
\end{align*}
Concerning the first term, we let $\Gamma:=\{\omega:\hat M^{\hat \pi^\iota}_{t,T}=0\}$ and have that $\tilde\alpha^{\eps,\iota}\equiv \check \alpha^{S,\eps,\iota}(\hat u^{k,\eps})$ on $\Gamma$. Consequently,
\begin{align*}
  &\check\E\big[|\hat J(t,x;\hat u^{k,\eps},\check \alpha^{S,\eps,\iota}(\hat u^{k,\eps}))-\hat J(t,x;\hat u^{k,\eps},\tilde\alpha^{\eps,\iota})|\big]
  \\
  &\leq \check\E\Big[\ett_{\Gamma^c}\big(|\Psi(\hat X^{t,x;\hat u^{k,\eps},\check \alpha^{S,m,\iota}(\hat u^{k,\eps})}_T)|+|\Psi(\hat X^{t,x;\hat u^{k,\eps},\tilde\alpha^{\eps,\iota}}_T)|
  \\
  &\quad+\int_t^{T} (|f(r,\hat X^{t,x;\check \alpha^{S,m,\iota}(\hat u^{k,\eps}),\hat u^{k,\eps}}_r)|+|f(r,\hat X^{t,x;\tilde\alpha^{\eps,\iota},\hat u^{k,\eps}}_r)|)dr
\\
&\quad + \sum_{j=1}^{\hat N^{k,\eps}}(\ell(\hat\tau_j,\hat X^{t,x;[\hat u^{k,\eps}]_{j-1},\check \alpha^{S,m,\iota}(\hat u^{k,\eps})}_{\hat\tau_j},\hat\beta_j) + \ell(\hat\tau_j,\hat X^{t,x;[\hat u^{k,\eps}]_{j-1},\tilde\alpha^{\eps,\iota}}_{\hat\tau_j},\hat\beta_j))
\\
&\quad + \sum_{j=1}^{\check M^{S,m,\iota}(\hat u^{k,\eps})}\chi(\check\eta_j^{S,m,\iota}(\hat u^{k,\eps}),\hat X^{t,x;\hat u^{k,\eps},[\check \alpha^{S,m,\iota}(\hat u^{k,\eps})]_{j-1}}_{\check\eta_j^{S,m,\iota}(\hat u^{k,\eps})},\check\theta_j^{S,m,\iota}(\hat u^{k,\eps}))
\\
&\quad+ \sum_{j=1}^{\tilde M^{\eps,\iota}}\chi(\tilde\eta^{\eps,\iota}_j,\hat X^{t,x;\hat u^{k,\eps},\tilde\alpha^{\eps,\iota}}_{\tilde\eta^{\eps,\iota}_j},\tilde\theta^{\eps,\iota}_j)\big)\Big]
\\
&\leq C\check\E\Big[\ett_{\Gamma^c}\check M^{S,m,\iota}(\hat u^{k,\eps})\big(1+\|\hat X^{t,x;\hat u^{k,\eps},\check \alpha^{S,m,\iota}(\hat u^{k,\eps})}\|^\rho_T+\|\hat X^{t,x;\hat u^{k,\eps},\check \alpha^{S,m,\iota}(\hat u^{k,\eps})}\|^\rho_T\big)\Big]
\\
&\leq C\check\E\Big[\hat M^{\hat \pi^\iota}\big(1+\|\hat X^{t,x;\hat u^{k,\eps},\check \alpha^{S,m,\iota}(\hat u^{k,\eps})}\|^\rho_T+\|\hat X^{t,x;\hat u^{k,\eps},\check \alpha^{S,m,\iota}(\hat u^{k,\eps})}\|^\rho_T\big)\Big].
\end{align*}
Since $\check\E\big[(\hat M^{\hat \pi^\iota})^2\big]\to 0$, $\hat\Prob$-a.s.,~as $\iota\to\infty$, there is a $\iota_1\in\bbN$ such that
\begin{align*}
  &\check\E\big[|\hat J(t,x;\hat u^{k,\eps},\check \alpha^{S,\eps,\iota}(\hat u^{k,\eps}))-\hat J(t,x;\hat u^{k,\eps},\tilde\alpha^{\eps,\iota})|\big]\leq \varrho,
\end{align*}
whenever $\iota\geq\iota_1$. Since $\bbA^{\eps,\iota}\in\bigSET$ for any $\eps>0$ and $\iota\in\bbN$, the second term can be made arbitrarily small according to Remark~\ref{app:rem:indep-of-filtr} by choosing $\eps>0$ sufficiently small. In particular, we find that
\begin{align*}
  \underline v^{\mcR,k}(t,x)\leq \underline v^{k,l}(t,x) + 4\varrho
\end{align*}
and since $\varrho>0$ was arbitrary we conclude that $\underline v^{\mcR,k}(t,x)\leq \underline v^{k,l}(t,x)$. As $l\in\bbN$ was arbitrary, we can thus take the limit as $l\to\infty$ and use Proposition~\ref{prop:trunk-vfs-approx} to arrive at the inequality $\underline v^{\mcR,k}(t,x)\leq \underline v^{k,\infty}(t,x)$.\qed\\

\appendix

\include{imp-imp-app}

\bibliographystyle{plain}
\bibliography{imp-imp_ref}
\end{document}

%% file: imp-imp-app.tex
% Appendix to imp-imp-game

\section{Some estimates related to the cost/reward functional}
We consider an arbitrary $\bbA=(\tilde\Omega,\tilde\mcF,\tilde\bbF,\tilde\Prob,\tilde W)\in\bigSET$ (or $\bbA=(\tilde\Omega,\tilde\mcF,\tilde\bbF,\tilde\Prob,\tilde W,\tilde\mu)\in\bigSET^\mcR$) and introduce the coupling (see~\cite{HamMM_SWG} for a similar definition in case of switching controls) of controls $u\in\mcU(\bbA)$ and $\alpha\in\mcA(\bbA)$ as $u\circ\alpha=(\phi_j,\zeta_j,\vartheta_j)_{j=1}^{M+N}$, where
\begin{align*}
\phi_j=\tau_{N_j+1}\wedge \eta_{M_j+1},
\end{align*}
where $N_1=M_1=0$ and
\begin{align*}
N_{j}&=N_{j-1}+\zeta_{j}
\\
M_{j}&=M_{j-1}+(1-\zeta_j)
\\
\zeta_j&=\ett_{[\tau_{N_{j-1}+1}\leq\eta_{M_{j-1}+1}]}
\end{align*}
for $j=1,\ldots,N+M$, and
\begin{align*}
\vartheta_j=\beta_{N_j}\ett_{[\zeta_j=1]}+\theta_{M_j}\ett_{[\zeta_j=0]}.
\end{align*}
We can thus define the set of coupled impulse controls $\bar\mcC_t(\bbA)$ as the set of all $\pi=(\phi_j,\zeta_j,\vartheta_j)_{j=1}^{L}$, where $(\phi_j)_{j\geq 0}$ is a non-decreasing sequence of $\tilde\bbF$-stopping times, $(\zeta_j,\vartheta_j)$ is $\tilde\mcF_{\phi_j}$-measurable and valued in $\{0\}\times A \cup \{1\}\times U$ while $L:=\sup\{j:\phi_j<T\}$ is $\tilde\Prob$-a.s.~finite. Finally, $X^{t,x;u,\alpha}=\mathcal X^{t,x;u\circ\alpha}$, where for each $\pi\in\bar\mcC_t(\bbA)$,
\begin{align*}
\mathcal X^{t,x;\pi}_s&=x+\int_t^s a(r,\mathcal X^{t,x;\pi}_r)dr+\int_t^s\sigma(r,\mathcal X^{t,x;\pi}_r)d\tilde W_r+\sum_{j=1}^L\ett_{[\phi_j\leq s]}\varsigma(\phi_j,\mathcal X^{t,x;[\pi]_{j-1}}_{\phi_j},\zeta_j,\vartheta_j)
\end{align*}
using the notation
\begin{align*}
  \varsigma(t,x,c,b):=\ett_{[c=1]}\jmp(t,x,b)+\ett_{[c=0]}\gamma(t,x,b)
\end{align*}
and $[\pi]_j:=(\phi_i,\zeta_i,\vartheta_i)_{i=1}^{L\wedge j}$.

First, we have the following moment estimate for the solution to the controlled SDE:
\begin{lem}\label{app:lem:SDEmoment}
Under Assumption~\ref{ass:onSDE}, the solution the SDE \eqref{ekv:fwd-sde} has moments of all orders, in particular, for each $p\geq 0$, there is a $C>0$, such that\footnote{When $\bbA\in\bigSET^\mcR$ we define $\tilde\E^\nu$ analogously to what was done in Section~\ref{sec:dual-game} and when $\bbA\in\bigSET$ we set $\tilde\E^\nu=\tilde\E$.} for any $\bbA\in\bigSET$,
\begin{align}\label{app:ekv:SDEmoment}
\sup_{\bbA\in\bigSET}\sup_{(u,\alpha,\nu)\in\bar\mcU(\bbA)\times\bar\mcA(\bbA)\times\mcV(\bbA)} \tilde\E^\nu\Big[\sup_{s\in[t,T]}|X^{t,x;u,\alpha}_s|^{p}\Big|\tilde\mcF_t\Big] \leq C(1+|x|^p),
\end{align}
$\tilde\Prob$-a.s.~for all $(t,x)\in[0,T]\times\R^d$.
\end{lem}

\noindent\emph{Proof.} The existence of an injective map from\footnote{Whenever we deem that this will not cause any confusion, we suppress the dependence on $\bbA$.} $\bar\mcU_t\times\bar\mcA_t$ to $\bar\mcC_t$ implies that the left hand side in \eqref{app:ekv:SDEmoment} is bounded by
\begin{align}\label{app:ekv:SDEmoment-pi}
\sup_{\bbA\in\bigSET}\sup_{\nu\in\mcV(\bbA)}\sup_{\pi\in\mcC(\bbA)}\tilde\E^\nu\Big[\sup_{s\in[t,T]}|\mcX^{t,x;\pi}_s|^{p}\Big|\tilde\mcF_t\Big]
\end{align}
We use the shorthand $\mcX^j:=\mcX^{t,x;[\pi]_j}$. By Assumption~\ref{ass:onSDE}.(\ref{ass:onSDE-Gamma}) we get for $s\in [\phi_{j},T]$, using integration by parts, that
\begin{align*}
|\mcX^{j}_s|^2&= |\mcX^{j}_{\phi_{j}}|^2+2\int_{(\phi_{j})+}^s \mcX^{j}_{r}dX^{j}_r+\int_{\phi_{j}+}^s d[\mcX^{j},\mcX^{j}]_r
\\
&\leq K^2_\Gamma\vee |\mcX^{{j-1}}_{\phi_{j}}|^2+2\int_{(\phi_{j})+}^s \mcX^{j}_{r} dX^{j}_r+\int_{\phi_{j}+}^s d[\mcX^{j},\mcX^{j}]_r.
\end{align*}
Now, either $|\mcX^{{j-1}}_{\phi_{j}}|\leq K_\Gamma$ in which case
\begin{align*}
|\mcX^{j}_s|^2\leq |x|^2\vee K^2_\Gamma+2\int_{\phi_{j}+}^s \mcX^{j}_{r} dX^{j}_r+\int_{\phi_{j}+}^s d[\mcX^{j},\mcX^{j}]_r.
\end{align*}
or $|\mcX^{{j-1}}_{\phi_{j}}|> K_\Gamma$ implying that
\begin{align*}
|\mcX^{j}_s|^2&\leq K^2_\Gamma\vee |\mcX^{{j-2}}_{\phi_{j-1}}|^2+2\int_{\phi_{j-1}+}^{\phi_{j}} \mcX^{j-1}_{r} dX^{j-1}_r+\int_{\phi_{j-1}+}^{\phi_{j}} d[\mcX^{j-1},\mcX^{j-1}]_r
\\
&\quad+2\int_{\phi_{j}+}^s \mcX^{j}_{r} dX^{j}_r+\int_{\phi_{j}+}^s d[\mcX^{j},\mcX^{j}]_r.
\end{align*}
In the latter case the same argument can be repeated and we conclude that
\begin{align}\label{ekv:X2-bound}
|\mcX^{j}_s|^2&\leq |x|^2\vee K_\Gamma^2+\sum_{i=j_0}^{j} \Big\{2\int_{\tilde\phi_{i}+}^{s\wedge\tilde\phi_{i+1}} \mcX^{i}_{r}dX^{i}_r+\int_{\tilde\phi_{i}+}^{s\wedge\tilde\phi_{i+1}} d[\mcX^{i},\mcX^{i}]_r\Big\},
\end{align}
where $\tilde\phi_0=t$, $\tilde\phi_i=\phi_i$ for $i=1,\ldots,j$ and $\tilde\phi_{j+1}=\infty$ and $j_0:=\max\{i\in \{1,\ldots,j\}:|\mcX^{{i-1}}_{\phi_{i}}|\leq K_\Gamma\}\vee 0$.

Now, since $\mcX^{i}$ and $\mcX^{j}$ coincide on $[0,\phi_{i+1\wedge j+1})$ we have
\begin{align*}
\sum_{i=j_0}^{j}\int_{\phi_{i}+}^{s\wedge\tilde\phi_{i+1}} \mcX^{i}_{r} dX^{i}_r
&=\int_{\phi_{j_0}}^s \mcX^{j}_{r}a(r,\mcX^{j}_r)dr+\int_{\phi_{j_0}}^{s}\mcX^{j}_{r}\sigma(r,\mcX^{j}_r)d\tilde W_r,
\end{align*}
and
\begin{align*}
\sum_{i=j_0}^{j} \int_{\tilde\phi_{i}+}^{s\wedge\tilde\phi_{i+1}} d[\mcX^{i},\mcX^{i}]_r&=\int_{\phi_{j_0}}^{s} \sigma^2(r,\mcX^{j}_r)dr.
\end{align*}
Inserted in \eqref{ekv:X2-bound} this that
\begin{align}\nonumber
|\mcX^{j}_s|^2&\leq |x|^2\vee K_\Gamma^2+\int_{\phi_{j_0}}^s (2\mcX^{j}_{s}a(r,\mcX^{j}_r)+\sigma^2(r,\mcX^{j}_r))dr+2\int_{\phi_{j_0}}^{s}\mcX^{j}_{r}\sigma(r,\mcX^{j}_r)d\tilde W_r
\\
&\leq |x|^2+C\Big(1+\int_{t}^{s}|\mcX^{j}_{r}|^2dr + \sup_{\eta\in[t,s]}\Big|\int_{t}^{\eta}\mcX^{j}_r\sigma(r,\mcX^{j}_r)d\tilde W_r\Big|\Big)\label{ekv:X-bound1}
\end{align}
for all $s\in [\zeta,T]$. Since $\tilde W$ is a standard Brownian motion under $\tilde\Prob^\nu$ we can apply the Burkholder-Davis-Gundy inequality to get that for $p\geq 2$,
\begin{align*}
\tilde\E^\nu\Big[\sup_{r\in[t,s]}|\mcX^{j}_r|^{p}\Big|\mcF_t \Big]\leq |x|^2+C\big(1+\tilde\E\Big[\int_{t}^{s}|\mcX^{j}_{r}|^{p}dr+\big(\int_{t}^{s}|\mcX^{j}_r|^4 dr\big)^{p/4}\Big|\tilde\mcF_t\Big]\big)
\end{align*}
and Gr\"onwall's lemma gives that for $p\geq 4$,
\begin{align*}
\tilde\E^\nu\Big[\sup_{s\in[t,T]}|\mcX^{j}_s|^{p}\Big|\tilde\mcF_t\Big]&\leq C(1+ |x|^{p}),
\end{align*}
where the constant $C=C(T,p)$ does not depend on $\nu$, $\pi$ or $j$ and \eqref{app:ekv:SDEmoment-pi} follows by letting $j\to\infty$ on both sides and using Fatou's lemma. The result for general $p\geq 1$ follows by Jensen's inequality.\qed\\

We move on to consider the cost/reward function
\begin{align*}
  J(t,x;u,\alpha,\nu,\bbA)&:=\tilde\E^\nu\Big[\psi(X^{t,x;u,\alpha}_T)+\int_t^{T}f(s,X^{t,x;u,\alpha}_s)ds - \sum_{j=1}^N\ell(\tau_j,X^{t,x;[u]_{j-1},\alpha},\beta_j)
  \\
  &\quad + \sum_{j=1}^M\chi(\eta_j,X^{t,x;u,[\alpha]_{j-1}},\theta_j)\Big].
\end{align*}
Our approach to relate the randomized game to the original one relies heavily on the fact that a slight perturbation of the control set does not lead to a dramatic change in the corresponding reward/cost $J$. An important fact in this regard is provided in the statement of the following lemma which shows that a form of continuity holds:
\begin{lem}\label{app:lem:SDEflow}
Let Assumption~\ref{ass:oncoeff} and Assumption~\ref{ass:onSDE} hold. Fix $\bbA\in\bigSET$, $K>0$, $k,l\in\bbN$ and sequences $(u^i)_{i\geq 0}$, $(\hat u^i)_{i\geq 0}$, $(\alpha^i)_{i\geq 0}$ and $(\hat\alpha^i)_{i\geq 0}$, with $u^i,\hat u^i\in\mcU^k_t(\bbA)$ and $\alpha^i,\hat \alpha^i\in\mcA^l_t(\bbA)$, such that
\begin{align*}
  |u^i-\hat u^i|+|\alpha^i-\hat\alpha^i|\to 0,\quad\text{as }i\to\infty
\end{align*}
and initial points $(t,x^i,\hat x^i)\in[0,T]\times\R^d\times\R^d$ with $|x^i|\leq K$ and $\nu\in\mcV(\bbA)$. Then,
\begin{align*}
  \tilde\E\big[|J(t,x^i;u^i,\alpha^i,\nu,\bbA)-J(t,\hat x^i;\hat u^i,\hat\alpha^i,\nu,\bbA)|\big]\to 0,
\end{align*}
if $|x^i-\hat x^i|\to 0$ as $i\to\infty$.
\end{lem}

The proof of the above lemma uses the following intermediate result defined in terms of impulse controls on $\mcC$. Here, $\mcC$ is the set of $\pi\in\bar\mcC$ such that
\begin{align*}
  \tilde\E\Big[\Big(\sum_{j=1}^L|\varphi(\phi_j,\mcX^{t,x;[\pi]_{j-1}}_{\phi_j},\zeta_j,\vartheta_j)|\Big)^2\Big]<\infty
\end{align*}
using the notation
\begin{align*}
  \varphi(t,x,c,b):=\ett_{[c=1]}\ell(t,x,b)-\ett_{[c=0]}\chi(t,x,b).
\end{align*}
For $\pi\in\mcC_t$, we define the corresponding reward as
\begin{align*}
  \mathcal J(t,x;\pi,\nu,\bbA)&:=\tilde\E^\nu\Big[\psi(\mcX^{t,x;\pi}_T)+\int_t^{T}f(s,\mcX^{t,x;\pi}_s)ds - \sum_{j=1}^L\varphi(\phi_j,\mcX^{t,x;[\pi]_{j-1}},\zeta_j,\vartheta_j)\Big]
\end{align*}
\begin{lem}\label{app:lem:SDEflow-pre}
Let Assumption~\ref{ass:oncoeff} and Assumption~\ref{ass:onSDE} hold. Fix $\bbA\in\bigSET$ and $\nu\in\mcV(\bbA)$. For each $t\in[0,T]$, $K>0$, $l\in\bbN$ and sequences $(x_i)_{i\geq 1}$ and $(\hat x_i)_{i\geq 1}$ in $\{x\in\R^d:|x|\leq K\}$ and $(\pi^i)_{i\geq 1}$ and $(\hat\pi^i)_{i\geq 1}$ in $\mcC_t^l(\bbA)$, we have
\begin{align*}
  \E\big[|\mathcal J(t,x_i;\pi^i,\nu,\bbA)-\mathcal J(t,\hat x_i;\hat \pi^i,\nu,\bbA)|\big]\to 0,\quad \text{as }i\to\infty
\end{align*}
whenever $|x_i-\hat x_i|\to 0$ as $i\to\infty$ and the impulse controls $\pi^i$ and $\hat\pi^i$ have the same number of interventions (\ie $\hat L^i=L^i$), $\zeta^i_j=\hat\zeta^i_j$ and $\hat\phi^i_j-\phi^i_j\searrow 0$ and $|\vartheta^i_j-\hat\vartheta^i_j|\to 0$ for $j=1\ldots,L$, $\tilde\Prob$-a.s.,~as $i\to\infty$.
\end{lem}

\noindent\emph{Proof.} The proof borrows from that of Lemma 5.3 in \cite{ElephantIH}. To simplify notation we let $\mcX^{i,j}:=\mcX^{t,x;[\pi^i]_j}$ and $\hat \mcX^{i,j}:=\mcX^{t,x;[\hat \pi^i]_j}$ for $j=0,\ldots,k$. Moreover, we let $\delta \mcX^{i,j}:=\hat \mcX^{i,j}-\mcX^{i,j}$. We have
\begin{align*}
\hat \mcX^{i,j}_{\hat\phi^i_{j}}-\mcX^{i,j}_{\hat\phi^i_{j}}&= \varsigma(\hat\phi^i_{j},\hat \mcX^{i,j-1}_{\hat\phi^i_{j}},\hat\zeta^i_{j},\hat\vartheta^i_{j}) - \varsigma(\phi^i_{j},\mcX^{i,j-1}_{\phi^i_{j}},\zeta^i_{j},\vartheta^i_{j})+\int_{\phi^i_{j}}^{\hat\phi^i_{j}} a(r,\mcX^{i,j}_r)dr+\int_{\phi^i_{j}}^{\hat\phi^i_{j}}\sigma(r,\mcX^{i,j}_r)d\tilde W_r
\\
&\leq|\varsigma(\hat\phi^i_{j},\hat \mcX^{i,j-1}_{\hat\phi^i_{j}},\hat\zeta^i_{j},\hat\vartheta^i_{j}) - \varsigma(\hat\phi^i_{j},\mcX^{i,j-1}_{\hat\phi^i_{j}},\zeta^i_{j},\vartheta^i_{j})|
+|\varsigma(\hat\phi^i_{j},\mcX^{i,j-1}_{\hat\phi^i_{j}},\zeta^i_{j},\vartheta^i_{j})- \varsigma(\phi^i_{j},\mcX^{i,j-1}_{\phi^i_{j}},\zeta^i_{j},\vartheta^i_{j})|
\\
&\quad+\int_{\phi^i_{j}}^{\hat\phi^i_{j}} a(r,\mcX^{i,j}_r)dr+\int_{\phi^i_{j}}^{\hat\phi^i_{j}}\sigma(r,\mcX^{i,j}_r)d\tilde W_r,
\end{align*}
where the last two terms tend to 0 in $L^2(\tilde\Omega,\tilde\mcF,\tilde\Prob)$ as $i\to\infty$ by Lemma~\ref{app:lem:SDEmoment}, while
\begin{align*}
|\varsigma(\hat\phi^i_{j},\mcX^{i,j-1}_{\hat\phi^i_{j}},\zeta^i_{j},\vartheta^i_{j})- \varsigma(\phi^i_{j},\mcX^{i,j-1}_{\phi^i_{j}},\zeta^i_{j},\vartheta^i_{j})|\to 0
\end{align*}
in $\tilde\Prob$-probability as $i\to\infty$, since $\varsigma$ is continuous.  On the other hand,
\begin{align*}
|\delta \mcX^{i,j}_{r}|&\leq |\delta  \mcX^{i,j}_{\hat \phi^i_{j}}|+\int_{\hat \phi^i_{j}}^{r}|\tilde a(s,\mcX^{i,j})-\tilde a(s,\tilde \mcX^{i,j})|ds
+\Big|\int_{\hat \phi^i_{j}}^{r}(\sigma(s,\mcX^{i,j})-\sigma(s,\tilde \mcX^{i,j}))d\tilde W_s\Big|
\end{align*}
The Burkholder-Davis-Gundy inequality now gives that
\begin{align*}
\tilde\E^\nu\Big[\sup_{r\in[\hat \phi^i_{j},s]}|\delta \mcX^{i,j}_{r}|^{2p}\Big] &\leq C\tilde\E\Big[|\delta \mcX^{i,j}_{\hat \phi^i_{j}}|^{2p}+ \Big(\int_{\hat \phi^i_{j}}^s|\tilde a(r,\mcX^{i,j})-\tilde a(r,\tilde \mcX^{i,j})|dr\Big)^{2p}
\\
& \quad+ \Big(\int_{\hat \phi^i_{j}}^{s}|\sigma(r,\tilde \mcX^{i,j})-\sigma(r,\mcX^{i,j})|^2dr\Big)^{p}\Big]
\\
&\leq C\tilde\E^\nu\Big[|\delta \mcX^{i,j}_{\hat \phi^i_{j}}|^{2p}+ \big(\int_{0}^s|\delta \mcX^{i,j}_{r}|^{2}dr\big)^p\Big]
\end{align*}
by the Lipschitz continuity of the coefficients, which by Gr\"onwall's lemma implies that
\begin{align*}
\tilde\E^\nu\Big[\sup_{r\in[\hat \phi^i_{j},T]}|\delta \mcX^{i,j}_{r}|^{2p}\Big] &\leq C\tilde\E\Big[|\delta \mcX^{i,j}_{\hat \phi^i_{j}}|^{2p}\Big].
\end{align*}
Since $\delta \mcX^{i,0}\equiv 0$, we can use induction to conclude that
\begin{align*}
  \tilde\E^\nu\big[\sup_{r\in[\hat \phi^i_{j},T]}|\delta \mcX^{i,j}_{r}|^{2p}\big]\to 0,\quad\text{as }i\to\infty.
\end{align*}
Now, as $\mcX^{t,x;\pi^i}-\mcX^{t,x;\hat\pi^i}=\delta \mcX^{i,j}$ on $[\hat \phi^i_{j},\phi^i_{j+1})$, we find that
\begin{align*}
  \tilde\E^\nu\big[\sup_{r\in[0,T]\setminus \cup_{j=1}^k [\phi^i_{j},\hat \phi^i_{j}]}|\delta \mcX^{i}_{r}|^{2p}\big]\to 0,\quad\text{as }i\to\infty,
\end{align*}
but the sequence $(\sup_{r\in [0,T]}(|\mcX^{t,x;\pi^i}_r|+|\mcX^{t,x;\hat\pi^i}_r|))_{i\in\bbN}$ is uniformly integrable and by continuity of $f$ and $\psi$, we conclude that
\begin{align*}
  \tilde\E^\nu\Big[|\psi(\mcX^{t,x;\pi^i}_T)-\psi(\mcX^{t,x;\hat\pi^i}_T)|+\int_t^{T}|f(s,\mcX^{t,x;\pi^i}_s)-f(s,\mcX^{t,x;\hat\pi^i}_s)|ds\Big]\to 0,\quad\text{as }i\to\infty.
\end{align*}
Similarly, continuity of $\varphi$ gives that
\begin{align*}
  \tilde\E^\nu\Big[\sum_{j=1}^{L^i}|\varphi(\phi^i_j,\mcX^{t,x;[\pi^i]_{j-1}},\zeta^i_j,\vartheta^i_j) - \varphi(\hat \phi^i_j,\mcX^{t,x;[\hat \pi^i]_{j-1}},\hat \zeta^i_j,\hat \vartheta^i_j)\Big]\to 0,\quad\text{as }i\to\infty
\end{align*}
and the desired result follows.\qed\\

\noindent\emph{Proof of Lemma~\ref{app:lem:SDEflow}.} Let $\pi^i=(\phi^i_j,\zeta^i_j,\vartheta^i_j)_{j=1}^{L^i}:=u^i\circ\alpha^i$ and $\hat\pi^i=(\hat\phi^i_j,\hat\zeta^i_j,\hat\vartheta^i_j)_{j=1}^{\hat L^i}:=\hat u^i\circ\hat\alpha^i$. For $\phi\in\mcC$, we let $N_0=M_0=0$ and define
\begin{align*}
N_j:=\min\{\iota>N_{j-1};\zeta_\iota=1\},\quad j=1,\ldots,N,
\\
M_j:=\min\{\iota>M_{j-1};\zeta_\iota=0\},\quad j=1,\ldots,M
\end{align*}
with $N=\sum_{j=1}^L\zeta_j$ and $M=L-N$.
We then introduce the operators $\Psi_{k'},\check\Psi_{k'}:\mcC(\bbA)\times\mcT\to \mcC(\bbA)$ for $k'\in \{1,\ldots,k\}$, defined as
\begin{align*}
  \Psi_{k'}(\pi,\tau)&:=(\ett_{[j\neq N_{k'}]}\phi_j+\ett_{[j=N_{k'}]}(\phi_j\vee(\phi_{j+1}\wedge\tau)) ,\zeta_j,\vartheta_j)_{j=1}^{L},
\end{align*}
whereas,
\begin{align*}
  \check\Psi_{k'}(\pi,\tau)&:=(\phi_{\varsigma^{k'}_j} ,\zeta_{\varsigma^{k'}_j},\vartheta_{\varsigma^{k'}_j})_{j=1}^{L},
\end{align*}
with
\begin{align*}
  \varsigma^{k'}_j:=j\ett_{[\notin \{N_{k'},N_{k'}+1\}]}+\ett_{[j= N_{k'}]}(j+\ett_{[\phi_{j}=\phi_{j+1}]}\ett_{[\phi_{j}<\tau]})+\ett_{[j= N_{k'}+1]}(j-\ett_{[\phi_{j-1}=\phi_{j}]}\ett_{[\phi_{j-1}<\tau]}).
\end{align*}
By first applying $\Psi_{k}$ to $(\pi^i,\hat \tau^i_k)$ we move $\pi^i_{N^i_{k}}$ representing $\tau^i_k$ towards $\hat\tau^i_k$, whenever $\tau^i_k<\hat\tau^i_k$ and $\phi^i_{N^i_k}<\phi^i_{N^i_k+1}$. However, $\Psi_{k}$ does not alter the order of interventions in $\pi^i$ and we may not end up with $\check\phi^i_{N^i_k}\geq \hat \tau^i_k$, where $\check\pi^{i,k,1}=\Psi_{k}(\pi^i,\hat\pi^i)$. To remedy this we apply $\check\Psi_{k}$ to the pair $(\check\pi^{i,k,1},\hat \tau^i_k)$ to switch positions of intervention $N^i_k$ and $N^i_k+1$, whenever intervention $N^i_k+1$ stood in the way of us reaching $\check\phi^i_{N^i_k}\geq \hat \tau^i_k$. We now let $\pi^{i,k,1}=\check\Psi_{k}(\check\pi^{i,k,1},\hat \tau^i_k)$ and repeat this process to produce a sequence $(\check\pi^{i,k,\iota},\pi^{i,k,\iota})_{\iota=1}^l$ through
\begin{align*}
\begin{cases}
  \check\pi^{i,k,\iota}=\Psi_{k}(\pi^{i,k,\iota-1},\hat \tau^i_k)
  \\
  \pi^{i,k,\iota}=\check\Psi_{k}(\check\pi^{i,k,\iota},\hat \tau^i_k)
\end{cases}\quad\iota=1,\ldots,l.
\end{align*}
Since $\alpha^i$ belongs to $\mcA^l$ it has at most $l$ interventions in the interval $[\tau_k,\hat\tau_k]$ and, therefor, $\phi^{i,k,l}_{N^{i,k,l}_k}\geq \hat\tau_k$. Moreover, as $\check\Psi_{k}$ applied to $(\check\pi^{i,k,\iota},\hat \tau^i_k)$ does not shift the position of intervention $\check N^{i,k,\iota}$ with the preceding intervention if not $\check\zeta^{i,k,\iota}_{\check N^{i,k,\iota}+1}=0$, we have by \eqref{ekv:jump-commute}-\eqref{ekv:interv-cost-commute} that
\begin{align*}
  \mathcal J(t,x_i;\check\pi^{i,k,\iota},\nu,\bbA)=\mathcal J(t,x_i;\pi^{i,k,\iota},\nu,\bbA).
\end{align*}
Hence,
\begin{align*}
  |\mathcal J(t,x_i;\pi^{i},\nu,\bbA)-\mathcal J(t,x_i;\pi^{i,k,l},\nu,\bbA)|\leq \sum_{\iota=1}^l|\mathcal J(t,x_i;\check\pi^{i,k,\iota},\nu,\bbA)-\mathcal J(t,x_i;\pi^{i,k,\iota-1},\nu,\bbA)|,
\end{align*}
where the left hand side tends to 0 in $L^1(\bbA)$ as $i\to\infty$ by Lemma~\ref{app:lem:SDEflow-pre}. We proceed by letting $\pi^{i,k-1,0}=\pi^{i,k,l}$ and repeat the above procedure by applying $\Psi_{k-1}$ followed by $\check\Psi_{k-1}$, $l$ times to get $\pi^{i,k-2,0}=\pi^{i,k-1,l}$ and so on. Eventually, we get to $\pi^{i,0,0}$ in which all interventions corresponding to $u$ have been shifted such that the corresponding intervention times are not dominated by the times in $\hat u$ and
\begin{align*}
  \E\big[|\mathcal J(t,x_i;\pi^{i},\nu,\bbA)-\mathcal J(t,x_i;\pi^{i,0,0},\nu,\bbA)|\big]\to 0,\quad\text{as }i\to\infty.
\end{align*}
The same approach can be taken to shift the times of interventions in $\pi^{i,0,0}$ corresponding to $\alpha^i$ so that they do not precede those of $\hat\alpha^i$, lets call the result $\bar\pi^{i}$. By symmetry we can then shift the time of the interventions in $\hat\pi^i$ to align in time with those of $\bar\pi^{i}$ and the assertion follows by again applying Lemma~\ref{app:lem:SDEflow-pre} after possibly shifting the order of some interventions in the resulting impulse control and using \eqref{ekv:jump-commute}-\eqref{ekv:interv-cost-commute}.\qed\\

Similar to what was noted in Remark 4.1 of \cite{Fuhrman15}, Lemma~\ref{app:lem:SDEflow} can be slightly generalized by allowing the impulse controls to be adapted to a sequence of filtrations $\tilde\bbF^i$:
\begin{rem}\label{app:rem:indep-of-filtr}
The assertion in Lemma~\ref{app:lem:SDEflow} remains true if we add a sequence $(\bbA^i)_{i\in\bbN}$ of extended setups, differing from $\bbA$ only in the filtration used, \ie $\mathbb A^i=(\tilde\Omega,\tilde\mcF,\tilde\bbF^i,\tilde\Prob,\tilde W)\in\bigSET$, and have that $u^i\in\mcU^k_t(\bbA^i)$ and $\alpha^i\in \mcA^l_t(\bbA^i)$.
\end{rem}
%%%%%%%%%%%%%%%%%%%%%%%%%%%%%%%%%%%%%%%%%%%%%%%%%%%%%%%%%%%%%%%%%%%%%%%%%%%%%%%%%%%%%%%%%%%%%%%%%%%%%%%%%%%%%%%%%%%%%%%%%%%%%%%%%%%%%%%%%%%%%%%%%%
%%%%%%%%%%%%%%%%%%%%%%%%%%%%%%%%%%%%%%%%%%%%%%%%%%%%%%%%%%%%%%%%%%%%%%%%%%%%%%%%%%%%%%%%%%%%%%%%%%%%%%%%%%%%%%%%%%%%%%%%%%%%%%%%%%%%%%%%%%%%%%%%%%
%%%%%%%%%%%%%%%%%%%%%%%%%%%%%%%%%%%%%%%%%%%%%%%%%%%%%%%%%%%%%%%%%%%%%%%%%%%%%%%%%%%%%%%%%%%%%%%%%%%%%%%%%%%%%%%%%%%%%%%%%%%%%%%%%%%%%%%%%%%%%%%%%%

\section{Connection between BSDEs with jumps and the dual game\label{app:dual-trunc}}

\subsection{Some preliminary result}

Some of the results gathered in this section uses the auxiliary lower barrier $h\in \Pi^g_c$ which is assumed to satisfy $h(T,\cdot)\leq\psi$. Hence, $\Psi(t,x):=\ett_{[t<T]}h(t,x)+\ett_{[t=T]}\psi(x)$ is jointly continuous on $[0,T)\times\R^d$ and $\lim_{(t,x')\to(T,x)}\Psi(t,x')\leq \Psi(T,x)$ for all $x\in\R^d$.

We give the following proposition, strategically formulated to streamline subsequent implementation processes. A proof (based on comparable findings documented in \cite{Dumitrescu15} and \cite{HamMor16}) can be found in \cite{qvi-stop}.

\begin{prop}\label{prop:rbsde-jmp}
For each $(t,x)\in [0,T]\times\R^d$ and $n\in\bbN$, there is a unique quadruple\\ $(Y^{t,x,n},Z^{t,x,n},V^{t,x,n},K^{+,t,x,n})\in\mcS^{2}_t \times \mcH^{2}_t(B) \times \mcH^{2}_t(\mu) \times \mcA^{2}_{t}$ such that
\begin{align}\label{ekv:rbsde-jmp-base}
  \begin{cases}
     Y^{t,x,n}_s=\psi(X^{t,x}_T)+\int_s^T f^n(r,X^{t,x}_r, Y^{t,x,n}_r, Z^{t,x,n}_r, V^{t,x,n}_{r})dr -\int_s^T  Z^{t,x,n}_r dB_r
    \\
    \quad-\int_{s}^T\!\!\!\int_A  V^{t,x,n}_{r}(e)\mu(dr,de)+(K^{+,t,x,n}_T- K^{+,t,x,n}_s),\quad \forall s\in [t,T],
    \\
    Y^{t,x,n}_s\geq h(s,X^{t,x}_s),\, \forall s\in [t,T] \quad\text{and}\quad\int_t^T\!\! \big(Y^{t,x,n}_s-h(s,X^{t,x}_s)\big)dK^{+,t,x,n}_s,
  \end{cases}
\end{align}
where
\begin{align*}
f^n(t,x,y,z,v):= f(t,x,y,z)-n\int_A(v(e)+\chi(t,x,e))^-\lambda(de).
\end{align*}
Moreover, $Y^{t,x,n}_s=\esssup_{\tau\in\mcT_s}Y^{t,x,n;\tau}_s$, where for each $\tau\in\mcT_t$, the triple
$(Y^{t,x,n;\tau},Z^{t,x,n;\tau},V^{t,x,n;\tau})\in\mcS^{2}_{[0,\tau]} \times \mcH^{2}_{[0,\tau]}(B) \times \mcH^{2}_{[0,\tau]}(\mu)$ satisfies
\begin{align}\nonumber
  Y^{t,x,n;\tau}_s&=\Psi(\tau,X^{t,x}_\tau)+\int_s^\tau f^n(r,X^{t,x}_r, Y^{t,x,n;\tau}_r, Z^{t,x,n;\tau}_r, V^{t,x,n;\tau}_{r})dr
  \\
  &\quad -\int_s^\tau  Z^{t,x,n;\tau}_r dB_r-\int_{s}^\tau\!\!\!\int_A  V^{t,x,n;\tau}_{r}(e)\mu(dr,de),\quad \forall s\in [0,\tau],\label{ekv:bsde-jmp}
\end{align}
and for each $\eta\in \mcT_t$, the stopping time
\begin{align*}
  \tau^*:=\inf\{s\geq \eta:Y^{t,x,n}_s=h(s,X^{t,x}_s)\}\wedge T
\end{align*}
is optimal in the sense that $Y^{t,x,n}_\eta=Y^{t,x,n;\tau^*}_\eta$. Finally, there is a function $v_n\in\Pi^g_c$ such that $v_n(s,X^{t,x}_s)=Y^{t,x,n}_s$ for all $s\in[t,T]$. %and $v_n$ is the unique viscosity solution (see Definition~\ref{rem:visc-non-loc} below for a definition) in $\Pi^g_c$ to
\begin{align}\label{ekv:obst-prob-n}
\begin{cases}
  \min\{v_n(t,x)-h(t,x),-\frac{\partial}{\partial t}v_n(t,x)-\mcL v_n(t,x)+\mcK^n v_n(t,x)\\
  \quad- f(t,x,v_n(t,x),\sigma^\top(t,x)\nabla_x v_n(t,x))\}=0,\quad\forall (t,x)\in[0,T)\times \R^d \\
  v_n(T,x)=\psi(x),
\end{cases}
\end{align}
where $\mcK^n \phi(t,x):=n\int_A(\phi(t,x+\gamma(t,x,e)) + \chi(t,x,e)-\phi(t,x))^-\lambda(de)$.
\end{prop}

Note that the variational inequality in \eqref{ekv:obst-prob-n} is non-local in the sense that, by invoking the operator $\mcK^n$, the driver depends on values of $v$ at points other that $(t,x)$. When defining the corresponding viscosity solutions, we adhere to the principle in \cite{HamMor16} (as opposed to the path taken in \cite{Barles97}) and define:

\begin{defn}\label{rem:visc-non-loc}
A locally bounded map $v:[0,T]\times\R^d\to\R$ is a viscosity supersolution (resp. subsolution) to \eqref{ekv:obst-prob-n} if it is l.s.c.~(resp. u.s.c.), satisfies $v(T,x)\geq \psi(x)$ (resp. $v(T,x)\leq \psi(x)$) and if for any $(t,x)\in [0,T)\times\R^d$ and $(p,q,X)\in \bar J^- v(t,x)$ (resp. $\bar J^+ v(t,x)$), we have
\begin{align*}
  \min\big\{&v(t,x)-h(t,x),-p-q^\top a(t,x)
  \\
  &-\frac{1}{2}\trace(\sigma\sigma^\top(t,x)X)+\mcK^n v(t,x)+f(t,x,v(t,x),\sigma^\top(t,x)q)\}\big\}\geq 0\quad (\text{resp. }\leq 0).
\end{align*}
It is a viscosity solution if $v_*$ is a supersolution and $v^*$ is a subsolution.
\end{defn}

\begin{rem}\label{rem:Vtx-cont}
The proof of Proposition~\ref{prop:rbsde-jmp} delineated in \cite{qvi-stop} borrows an important fact from \cite{HamMor16} (see Proposition 3.1 therein), namely that the jump component can be written $V^{t,x,n}_s(e)=v_n(s,X^{t,x}_s+\gamma(s,X^{t,x}_s,e))-v_n(s,X^{t,x}_s)$. Of particular importance for our subsequent analysis is the fact that the map $(t,x,e)\mapsto V^{t,x,n}_t(e)$ can be choosen to be deterministic and continuous.
\end{rem}

\begin{prop}
There is a sequence $(v_{k,n})_{k,n\in\bbN}$, with $v_{k,n}\in\Pi^g_c$ such that $v_{k,n}(s,X^{t,x}_s)=Y^{t,x,k,n}_s$ for all $(t,x)\in[0,T]\times\R^d$ and $s\in[t,T]$. Moreover, for $k>0$, $v_{k,n}$ is a viscosity solution to
\begin{align}\label{ekv:var-ineq-nk}
\begin{cases}
  \min\{v_{k,n}(t,x)-\supOP v_{k-1,n}(t,x),-\frac{\partial}{\partial t}v_{k,n}(t,x)-\mcL v_{k,n}(t,x)+\mcK^n v_{k,n}(t,x)\\
\quad-f(t,x,v_{k,n}(t,x),\sigma^\top(t,x)\nabla_x v_{k,n}(t,x))\}=0,  \quad\forall (t,x)\in[0,T)\times \R^d \\
  v_{k,n}(T,x)=\psi(x),
\end{cases}
\end{align}
where we recall that $\mcK^n \phi(t,x):=n\int_A(\phi(t,x+\gamma(t,x,e))+\chi(t,x,e)-\phi(t,x))^-\lambda(de)$.
\end{prop}

\noindent\emph{Proof.} The result follows by repeated use of Proposition~\ref{prop:rbsde-jmp}.\qed\\

\subsection{Proof of Proposition~\ref{prop:trunk-game}}
For $(u,\nu)\in\mcU^\mcR\times\mcV$, we let $(P^{t,x;u,\nu},Q^{t,x;u,\nu},R^{t,x;u,\nu}_{r})$ be the unique solution with $P^{t,x;u,\nu}-\Lambda^{t,x;u}\in\mcS^2_t$ and $(Q^{t,x;u,\nu},R^{t,x;u,\nu})\in\mcH^2_t(B)\times\mcH^2_t(\mu)$ to
\begin{align}\nonumber
P^{t,x;u,\nu}_s&=\psi(X^{t,x;u}_T)+\int_s^T f(r,X^{t,x;u}_r)dr-\int_s^T Q^{t,x;u,\nu}_r dB_r
\\
&\quad-\int_{s}^T\!\!\!\int_A R^{t,x;u,\nu}_{r}(e)\mu^\nu(dr,de) + \int_{s}^T\!\!\!\int_A\chi(r,X^{t,x;u}_{r},e)\nu_r(e)\lambda(de)dr - \Lambda^{t,x;u}_{T}+\Lambda^{t,x;u}_s, \label{ekv:non-ref-bsde-simp}
\end{align}
where $\mu^\nu(dr,de):=\mu(dr,de)-\nu_r(e)\lambda(de)ds$. Note that $P^{t,x;u,\nu}_t=\E^\nu\big[P^{t,x;u,\nu}_t\big|\mcF^{B,\mu}_t\big]=J^\mcR(t,x;u,\nu)$. The proof of Proposition~\ref{prop:trunk-game} is based on applying the comparison theorem for BSDEs with jumps (see Theorem 4.2 in \cite{QuenSul13}) to $P^{t,x;u,\nu}$.\\

\noindent\emph{Proof of Proposition~\ref{prop:trunk-game}.} We let $\nu^{S,n,*}\in\mcV^{S,n}$ be defined as
\begin{align*}
  &\nu^{S,n,*}(u)(s,e):=n\ett_{[0,\tau_1]}(s) \ett_{\Upsilon^{t,x,k}_s(e)}+n\sum_{j=1}^k\ett_{(\tau_j,\tau_{j+1}]}(s) \ett_{\Upsilon^{t,x,u,k-j}_s(e)}
\end{align*}
for all $(u,s,e)\in \mcU^\mcR_t\times [t,T]\times A$, where
\begin{align*}
  \Upsilon^{t,x,u,j}_s(e):=\{(\omega,s,e): v_{j,n}(s,X^{t,x;u}_{s-}+\jmp(s,X^{t,x;u}_{s-},e))+\chi(s,X^{t,x;u}_{s-},e)-v_{j,n}(s,X^{t,x;u}_{s-} + \jmp(s,X^{t,x;u}_{s-},e))<0\}.
\end{align*}
Remark~\ref{rem:Vtx-cont} then implies that
\begin{align*}
  \nu^{S,n,*}(u)(r,e)=n\ett_{[0,\tau_1]}(r)\ett_{[V^{t,x,k,n}_{r}(e)<-\chi(r,X^{t,x}_{r},e)]}+n\sum_{j=1}^k\ett_{(\tau_j,\tau_{j+1}]}(r) \ett_{[V^{\tau_j,X^{t,x;u}_{\tau_j},k-j,n}_{r}(e)<-\chi(r,X^{t,x;u}_{r},e)]}
\end{align*}
$d\Prob\otimes\lambda(de)\otimes dr$-a.e. Moreover, let $u^{k,*}:=(\tau^*_j,\beta^*_j)_{j=1}^{N^*}\in\mcU^{\mcR,k}_{t}$ be defined through
\begin{itemize}
  \item $\tau^*_{j}:=\inf \big\{r \geq \tau^*_{j-1}:\:v_{k-j+1,n}(r,X^{t,x;[u^*]_{j-1}}_{r})=\supOP v_{k-j,n}(r,X^{t,x;[u^*]_{j-1}}_r)\big\}\wedge T$,
  \item $\beta^*_j\in\mathop{\arg\max}_{b\in U}\{v_{k-j,n}(\tau^*_{j},X^{t,x;[u^*]_{j-1}}_{\tau^*_{j}}+\jmp(\tau^*_{j},X^{t,x;[u^*]_{j-1}}_{\tau^*_{j}},b))-\ell(\tau^*_{j},X^{t,x;[u^*]_{j-1}}_{\tau^*_{j}},b)\}$,
\end{itemize}
for $j=1,\ldots,k$ with $\tau_0^*:=t$, while $N^*:=\max\{j\in \{0,\ldots,k\}:\tau^*_j<T\}$. In the proof, which is divided into two steps, we show that $(u^{k,*},\nu^{S,n,*})\in\mcU_t^{\mcR,k}\times \mcV_t^{S,n}$, is a type of saddle point for the game in \eqref{ekv:trunk-game-value}.\\

\textbf{Step 1:} We first show that $Y^{t,x,k,n}_t=P^{t,x;u^{k,*},\nu^{n,*}}_t$, with $\nu^{n,*}:=\nu^{S,n,*}(u^{k,*})$. We have that, $K^{+,t,x,k,n}_{\tau^*_1}-K^{+,t,x,k,n}_{t}=0$, $\Prob$-a.s., and we find that
\begin{align*}
  Y^{t,x,k,n}_{s}&= \ett_{[\tau^*_1<T]}\supOP v_{k-1,n}(\tau^*_1,X^{t,x}_{\tau_1^*}) +\ett_{[\tau^*_1=T]}\psi(X^{t,x}_T) + \int_{s}^{\tau_1^*} f^n(r,X^{t,x}_r,Y^{t,x,k,n}_r,Z^{t,x,k,n}_r,V^{t,x,k,n}_{r})dr
  \\
  &\quad - \int_{s}^{\tau_1^*}Z^{t,x,k,n}_r dB_r-\int_{{s}}^{\tau_1^*}\!\!\!\int_A V^{t,x,k,n}_{r}(e)\mu(dr,de)
  \\
  &= \ett_{[\tau^*_1<T]}\{v_{k-1,n}(\tau^*_1,X^{t,x;[u^{k,*}]_1}_{\tau_1^*})-\ell({\tau_1^*},X^{t,x}_{\tau_1^*},\beta_1^*)\} +\ett_{[\tau^*_1=T]}\psi(X^{t,x}_T)
  \\
  &\quad +\int_{s}^{\tau_1^*} f^n(r,X^{t,x}_r,Y^{t,x,k,n}_r,Z^{t,x,k,n}_r,V^{t,x,k,n}_{r})dr - \int_{s}^{\tau_1^*}Z^{t,x,k,n}_r dB_r -\int_{{s}}^{\tau_1^*}\!\!\!\int_A V^{t,x,k,n}_{r}(e)\mu(dr,de)
\end{align*}
for all ${s}\in [t,\tau^*_1]$. To simplify notation later, we let $(\mcY^0,\mcZ^0,\mcV^0):=(Y^{t,x,k,n},Z^{t,x,k,n},V^{t,x,k,n})$. Similarly, there is a quadruple $(\mcY^1,\mcZ^1,\mcV^1,\mcK^1)\in\mcS^2_t\times\mcH^2_t(B)\times\mcH^2_t(\mu)\times\mcA^2_{t}$ that solves the reflected bsde
\begin{align*}
\begin{cases}
  \mcY^{1}_{s}=\psi(X^{t,x;[u^{k,*}]_1}_T)+\int_{s}^T f^n\big(r,X^{t,x;[u^{k,*}]_1}_r,\mcY^1_r,\mcZ^1_r,\mcV^1_r\big)dr-\int_{s}^T \mcZ^1_r dB_r-\int_{{s}}^{T}\!\!\!\int_A \mcV^{1}_{r}(e)\mu(dr,de)
  \\
  \quad+ \mcK^1_T-\mcK^1_{s},\quad\forall {s}\in[\tau^*_1,T], \\
  \mcY^{1}_{s}\geq \supOP v_{k-2,n}({s},X^{t,x;[u^{k,*}]_1}_{s}),\:\forall {s}\in[\tau^*_1,T]\quad{\rm and}\quad
  \int_{\tau_1^*}^T(\mcY^{1}_{s}-\supOP v_{k-2,n}({s},X^{t,x;[u^{k,*}]_1}_{s}))d\mcK^{1}_{s}=0.
\end{cases}
\end{align*}
Continuity together with an approximation of $X^{t,x;[u^{k,*}]_1}$ and a stability result for reflected BSDEs with jumps akin to the one in Proposition 3.6 of \cite{ElKaroui1} (see \eg Proposition A.1 in \cite{Dumitrescu15}) implies that $\mcY^{1}_{{s}}=v_{k-1,n}({s},X^{t,x;[u^{k,*}]_1}_{{s}})$, $\Prob$-a.s., for each ${s}\in [\tau^*_1,T]$. In particular, as Proposition~\ref{prop:rbsde-jmp} gives that $\mcK^1_{\tau^*_2}-\mcK^1_{\tau^*_1}=0$, $\Prob$-a.s., we conclude that
\begin{align*}
  \mcY^{1}_{s}&= \ett_{[\tau^*_2<T]}\{v_{k-2,n}(\tau^*_2,X^{t,x;[u^{k,*}]_2}_{\tau_2^*})-\ell({\tau_2^*},X^{t,x;[u^{k,*}]_1}_{\tau_2^*},\beta_2^*)\} +\ett_{[\tau^*_2=T]}\psi(X^{t,x;u^{k,*}}_T)
  \\
  &\quad \int_{s}^{\tau_2^*} f^n(r,X^{t,x}_r,\mcY^{1}_r,\mcZ^{1}_r,\mcV^{1}_r)dr - \int_{s}^{\tau_2^*}\mcZ^{1}_r dB_r - \int_{{s}}^{\tau_2^*}\!\!\!\int_A \mcV^{1}_{r}(e)\mu(dr,de)
\end{align*}
for all ${s}\in [\tau^*_1,\tau^*_2]$.

Repeating this process $k$ times we find that there is a sequence $(\mcY^j,\mcZ^j,\mcV^j)_{j=0}^k\subset\mcS^2\times\mcH^2(B)\times\mcH^2(\mu)$ such that $\mcY:=\ett_{[t,\tau^*_{1}]}\mcY^0+\sum_{j=1}^k\ett_{(\tau^*_{j},\tau^*_{j+1}]}\mcY^j$, $\mcZ:=\sum_{j=0}^k\ett_{(\tau^*_{j},\tau^*_{j+1}]}\mcZ^j$ and $\mcV:=\sum_{j=0}^k\ett_{(\tau^*_{j},\tau^*_{j+1}]}\mcV^j$ satisfies
\begin{align*}
  \mcY_{s}&=\psi(X^{t,x;u^{k,*}}_T)+\int_{s}^{T} f^n\big(r,X^{t,x;u^{k,*}}_r,\mcY_r,\mcZ_r,\mcV_r\big)dr-\int_{s}^{T}\mcZ_rdB_r - \int_{{s}}^{T}\!\!\!\int_A \mcV_{r}(e)\mu(dr,de)
  \\
  &\quad -\sum_{j=1}^{N^*}\ett_{[{s}\leq \tau^*_j]}\ell({\tau_j^*},X^{t,x;[u^{k,*}]_{j-1}}_{\tau_j^*},\beta_j^*)
  \\
  &=\psi(X^{t,x;u}_T)+\int_s^T f(r,X^{t,x;u}_r)dr-\int_s^T \mcZ_r dB_r
\\
&\quad-\int_{s}^T\!\!\!\int_A\mcV_{r}(e)\mu^{\nu^{n,*}}(dr,de) + \int_{s}^T\!\!\!\int_A\chi(r,X^{t,x;u}_{r},e)\nu^{n,*}_r(e)\lambda(de)dr - \Lambda^{t,x;u}_{T}+\Lambda^{t,x;u}_s
\end{align*}
for all ${s}\in[t,T]$. Now, $\mcY-\Lambda^{t,x;u}\in\mcS^2_t$ and $(\mcZ,\mcV)\in\mcH^2_t(B)\times\mcH^2_t(\mu)$ and by uniqueness of solutions to \eqref{ekv:non-ref-bsde-simp} we conclude that $Y^{t,x,k,n}_t=P^{t,x;u^{k,*},\nu^{n,*}}_t$.\\

\textbf{Step 2:} To establish \eqref{ekv:trunk-game-value}, we show that for any $(u,\nu)\in\mcU_t^{\mcR,k}\times\mcV^{n}$ it holds that
\begin{align}\label{ekv:is-saddle-point}
  P^{t,x;u,\nu^{S,n,*}(u)}_t\leq P^{t,x;u^{k,*},\nu^{S,n,*}(u^*)}_t \leq P^{t,x;u^{k,*},\nu}_t.
\end{align}
Concerning the first inequality, suppose that $\hat u:=(\hat\tau_j,\hat\beta_j)_{j=1}^{\hat N}\in \mcU^{\mcR,k}_t$ is another impulse control, then
\begin{align*}
  Y^{t,x,k,n}_{s}&= Y^{t,x,k,n}_{\hat\tau_1}+\int_{s}^{\hat\tau_1} f^n(r,X^{t,x}_r,Y^{t,x,k,n}_r,Z^{t,x,k,n}_r,V^{t,x,k,n}_{r})dr - \int_{s}^{\hat\tau_1}Z^{t,x,k,n}_r dB_r
  \\
  &\quad-\int_{{s}}^{\hat\tau_1}\!\!\!\int_A V^{t,x,k,n}_{r}(e)\mu(dr,de)+K^{+,t,x,k,n}_{\hat\tau_1}-K^{+,t,x,k,n}_{s}
  \\
  &\geq \ett_{[\hat\tau_1<T]}\{v_{k-1,n}(\hat\tau_1,X^{t,x}_{\hat\tau_1}+\jmp(\hat\tau_1,X^{t,x}_{\hat\tau_1},\hat \beta_1)) -\ell({\hat\tau_1},X^{t,x}_{\hat\tau_1},\hat\beta_1)\} +\ett_{[\hat\tau_1=T]}\psi(X^{t,x}_T)
  \\
  &\quad +\int_{s}^{\hat\tau_1} f^n(r,X^{t,x}_r,Y^{t,x,k,n}_r,Z^{t,x,k,n}_r,V^{t,x,k,n}_{r})dr - \int_{s}^{\hat\tau_1}Z^{t,x,k,n}_r dB_r-\int_{{s}}^{\hat\tau_1}\!\!\!\int_A V^{t,x,k,n}_{r}(e)\mu(dr,de)\\
  &\quad+K^{+,t,x,k,n}_{\hat\tau_1}-K^{+,t,x,k,n}_{s}
\end{align*}
for all ${s}\in [t,\hat\tau_1]$. Arguing as above gives that there is a sequence $(\hat\mcY^j,\hat\mcZ^j,\hat\mcV^j,\hat\mcK^j)_{j=0}^k\subset\mcS^2\times\mcH^2(B)\times\mcH^2(\mu)\times\mcA^2$ such that letting $\hat\mcY:=\ett_{[t,\hat\tau_{1}]}\hat\mcY^0+\sum_{j=1}^k\ett_{(\hat\tau_{j},\hat\tau_{j+1}]}\hat\mcY^j$, $\hat\mcZ:=\sum_{j=0}^k\ett_{(\hat\tau_{j},\hat\tau_{j+1}]}\hat\mcZ^j$, $\hat\mcV:=\sum_{j=0}^k\ett_{(\hat\tau_{j},\hat\tau_{j+1}]}\hat\mcV^j$ and $\hat\mcK_{s}:=\sum_{j=0}^{k}\ett_{[\hat \tau_j<{s}]}\{\hat\mcK^j_{{s}\wedge\hat\tau_{j+1}}-\hat\mcK^j_{\hat\tau_j}\}$, with $\hat\tau_0:=t$, implies that $(\hat\mcY,\hat\mcZ,\hat\mcV,\hat\mcK)$ satisfies
\begin{align*}
  \hat\mcY_{s}&\geq \psi(X^{t,x;\hat u}_T)+\int_{s}^{T} f^n\big(r,X^{t,x;\hat u}_r,\hat\mcY_r,\hat\mcZ_r,\hat\mcV_r\big)dr-\int_{s}^{T} \hat\mcZ_rdB_r-\int_{{s}}^{\hat\tau_1}\!\!\!\int_A \hat\mcV_{r}(e)\mu(dr,de)
  \\
  &\quad-\sum_{j=1}^{\hat N}\ett_{[{s}\leq\hat\tau_j]}\ell({\hat \tau_j},X^{t,x;[\hat u]_{j-1}}_{\hat \tau_j},\hat \beta_j)+\hat\mcK_{T}-\hat \mcK_{{s}},\quad\forall {s}\in[t,T]
\end{align*}
and $\hat\mcY_t=Y^{t,x,k,n}_t$. Now, the comparison theorem for BSDEs with jumps gives that $P^{t,x;\hat u,\nu^{S,n,*}(\hat u)}_t\leq\hat\mcY_t$ and we conclude that $P^{t,x;\hat u,\nu^{S,n,*}(\hat u)}_t \leq Y^{t,x,k,n}_t$ proving the first inequality in \eqref{ekv:is-saddle-point}. The second inequality is an immediate consequence of comparison.\qed\\